\documentclass[12pt,  oneside
]{amsart}
\usepackage{amscd,amsmath, wasysym}
\title[]{Regularity for the CR vector bundle problem II}
\author[]{Xianghong Gong$^*$
\and Sidney M. Webster}

 \address{Department of Mathematics,
 University of Wisconsin, Madison, WI 53706, U.S.A.}
\email{gong@math.wisc.edu}
 \address{Department of Mathematics,
 University of Chicago, Chicago, IL 60637, U.S.A.}
 \email{webster@math.uchicago.edu}

  \thanks{$^*$ Research  supported in part by NSF grant DMS-0705426.}

\newcommand{\dist}{\operatorname{dist}}

\newtheorem{thm}{Theorem}[section]
\newtheorem{cor}[thm]{Corollary}
\newtheorem{prop}[thm]{Proposition}
\newtheorem{lemma}[thm]{Lemma}

\theoremstyle{definition}

\newtheorem{rem}[thm]{Remark}

\renewcommand{\th}[1]{\begin{thm}\label{#1}}
\newcommand{\co}[1]{\begin{cor}\label{#1}}
\renewcommand{\le}[1]{\begin{lemma}\label{#1}}
\newcommand{\pr}[1]{\begin{prop}\label{#1}}
\newcommand{\pf}{\begin{proof}}
\newcommand{\epf}{\end{proof}}

\newcommand{\ga}{\begin{gather}}
\newcommand{\gan}{\begin{gather*}}
\newcommand{\al}{\begin{align}}
\newcommand{\aln}{\begin{align*}}
\newcommand{\eq}[1]{\begin{equation}\label{#1}}
\newcommand{\eeq}{\end{equation}}

\newcommand{\ci}{~\cite}

\newcommand{\f}[2]{\frac{#1}{#2}}

\newcommand{\df}{\overset{\text{def}}{=\!\!=}}

\newcommand{\cc}{{\bf C}}

\newcommand{\rr}{{\bf R}}

\newcommand{\ov}{\overline}

\newcommand{\RE}{\operatorname{Re}}
\newcommand{\IM}{\operatorname{Im}}
\newcommand{\dbar}{\overline\partial}
\newcommand{\dbm}{\overline\partial_M}

\newcommand{\tr}[1]{{\operatorname {\mathcal R}_{#1}}}
\newcommand{\trd}{{\operatorname {\mathcal R}}}

\newcommand{\cC}{\mathcal C}

\newcommand{\ps}{the cone property}

\newcommand{\rzd}{r_z\cdot(\zeta-z)}
\newcommand{\rzed}{r_\zeta\cdot(\zeta-z)}

\newcommand{\pd}{\partial}

\newcommand{\zs}{\zeta_*}
\newcommand{\xis}{\xi_*}

\newcommand{\yt}{\frac{1}{2}}

\newcommand{\yf}{\frac{1}{4}}

\newcommand{\re}[1]{(\ref{#1})}
\newcommand{\rea}[1]{$(\ref{#1})$}
\newcommand{\rl}[1]{Lemma~\ref{#1}}

\newcommand{\rp}[1]{Proposition~\ref{#1}}
\newcommand{\rt}[1]{Theorem~\ref{#1}}

\newcommand{\rrem}[1]{Remark~\ref{#1}}

\newcommand{\cont}{continuous}

\newcommand{\supp}{\operatorname{supp}}

\newcommand{\db}{\dbar}
\newcommand{\dbb}{\dbar_b}

\newcommand{\hf}{homotopy formula}

\newcounter{pp}
\newcommand{\bpp}{\begin{list}{$\hspace{-1em}\alph{pp})$}{\usecounter{pp}}}
\newcommand{\epp}{\end{list}}

\newcounter{ppp}
\newcommand{\bppp}{\begin{list}{$\hspace{-1em}(\roman{ppp})$}{\usecounter{ppp}}}
\newcommand{\eppp}{\end{list}}

\begin{document}

 \begin{abstract}  We derive a $\mathcal C^{k+\yt}$ H\"older estimate for
 $P\varphi$, where $P$ is
either of the two solution operators in Henkin's local homotopy
formula for $\overline\partial_b$ on a strongly pseudoconvex real
hypersurface $M$ in $\mathbf C^{n}$, $\varphi$ is a $(0,q)$-form
of class $\mathcal C^{k}$
 on $M$, and $k\geq0$ is
an integer.  We also derive a $\mathcal C^{a}$ estimate for
$P\varphi$, when $\varphi$ is of class $\mathcal C^{a}$  and
$a\geq0$ is a real number.  These estimates require that $M$ be of
class $\mathcal C^{k+\frac{5}{2}}$, or $\mathcal C^{a+2}$,
respectively. The explicit bounds for the constants occurring in
these estimates also considerably improve previously known such
results.

These estimates are then applied to the integrability problem for
CR vector bundles to gain improved regularity.  They also constitute
a major ingredient in a forthcoming  work of the authors on the local
CR embedding problem.
   \end{abstract}

%


 \maketitle

\tableofcontents

\setcounter{section}{0}
\setcounter{thm}{0}\setcounter{equation}{0}
\section{ Introduction}\label{sec1}
 In this paper we will prove  the following. \th{1} Let $n\geq4$.
Let  $M$ be 
a strongly
pseudoconvex real hypersurface  in $\cc^n$ of class $\cC^2$.  Let
$\omega$ be an $r\times r$ matrix of continuous $(0,1)$-forms on $M$.
Assume that 
$\dbar_M\omega=\omega\wedge\omega$.
Near each point of $M$,
there exists a  non-singular matrix $A$ of H\"older class $\cC^{1/2}$
satisfying $\dbar_M
A=-A\omega$ 
and
\bpp
\item
$A\in \cC^{a}(M)$, if $\omega\in \cC^{a}(M)$,
$M\in \cC^{a+2}$ and   $a>0$ is a   real number;
\item
$A\in \cC^{k+\yt}(M)$, if  $\omega\in \cC^{k}(M)$, $M\in \cC^{k+\f{5}{2}}$
and   $k>0$ is an integer.
\epp
\end{thm}

If $\omega$ and $A$ are of class $\cC^0$, 
the identities $\dbar_M\omega=\omega\wedge\omega$
and $\dbar_M A=-A\omega$ are in the sense of currents; see section~\ref{sec3}.
This work is a continuation of~\ci{GWzese}.
For earlier results  see~\ci{Wenion} and Ma-Michel~\ci{MMnifi}.

We now describe the above result in terms of an integrability
problem for CR vector bundles (\cite{Wenion}).
Let $E$ be a   complex vector bundle of rank $r$ over $M$ with
a connection $D$.  For a
local frame $e=(e_1,\dots, e_r)$,   $De_i=\omega_i^je_j$,
where $ \omega=(\omega_i^j) $ are connection $1$-forms on $M$;
by a frame change $\tilde
e=Ae$,   $D\tilde e=\tilde\omega\tilde
e$ with $\tilde \omega=(dA+A\omega)A^{-1}$.
  The integrability problem   is to
find an $A$ such that the new connection forms
$\tilde\omega=(\tilde\omega_i^j)$ belong to  the ideal $\mathcal J(M)$
generated by $(1,0)$-forms on $M$.
The integrability  condition is that   the curvature
2-forms $d\omega -\omega\wedge\omega$ belong to $\mathcal J(M)$.

\

We want to mention a few ingredients in the proof.  The integrability problem is  local.
As in~\ci{Weeini},~\ci{GWzese}, we will
use   the Henkin local homotopy formula. Let  $M\subset\cc^n$ be a graph
over a domain
$D\subset\rr^{2n-1}$,
 given by $y^n=|z'|^2+\hat r(z',x^n)$
with $\hat r(0)=\pd \hat r(0)=0$, where $z=(z',z^n)$ are standard coordinates.
Set
$M_\rho=M\cap\{(x^n)^2+y^n<\rho^2\}$ and $D_\rho=\pi(M_\rho)$.
Suppose that $\ov D_{\rho_0}\subset
D$ and
the $\cC^2$-norm  $ \|\hat r\|_{\rho_0,2}$
  of $\hat r$ on $D_{\rho_0}$ is sufficiently small.
For
$0<\rho\leq\rho_0$ and $n\geq4$,
 we have the   Henkin   homotopy formula
 \eq{Hhf}
 \varphi=
 \dbm P\varphi+Q\dbm\varphi
 \end{equation}
for $(0,q)$-forms $\varphi$ on $\ov{M_\rho}$ with $0<q<n-2$.
 We will prove  the following estimates.
 \bppp
\item Let $a\geq0$ be a real number. Then
 $$ 
 \|P\varphi\|_{ (1-\sigma)\rho,a}\leq
C_a\rho^{-s_*}\sigma^{-s} \bigl(\|\varphi\|_{ \rho,a}+
\|\hat r\|_{\rho,a+2}\|\varphi\|_{\rho,0}\bigr).
$$ 
\item Let $k\geq0$ be an integer. Then
\gan
 \|P\varphi\|_{  (1-\sigma)\rho,k+\f{1}{2}}\leq
C_k\rho^{-s_*}\sigma^{-s}\bigl((1+\|\hat r\|_{\rho,\f{5}{2}})
\|\varphi\|_{\rho,k}+
\| \hat r\|_{\rho, k+\f{5}{2}}\|\varphi\|_{\rho,0}\bigr);\\
\|P\varphi\|_{ (1-\sigma)\rho, 1/2}\leq
{C}  \rho^{-1}{\sigma^{1-2n}} \|\varphi\|_{\rho,0},\qquad  k=0,\quad
q=1.
\end{gather*}
\eppp
(See \re{sval} for $s,s_*$.) We emphasize that the estimates hold for all
 $0<\rho\leq\rho_0\leq3$ and $0<\sigma<1$.
Under the coordinates $(z',x^n)$ of $M$,
    $\|\cdot\|_{\rho,a}$ denotes the standard $\cC^a$-norm
on the domain $D_\rho\subset\rr^{2n-1}$.
The same estimates hold for $Q$; however,
the second estimate in (ii), based on a special property of the kernels for $(0,1)$
forms,
is not applicable
to $Q$  when it operates on $(0,q+1)$ form with $q>0$.
See Romero~\ci{Ronion} for
estimates in H\"older
norms  for the Heisenberg group case
and an example showing necessity of blow-up constants.

  The estimate $(i)$ is proved
  in  \rp{dkck} and  $(ii)$  is in \rp{dkck+}.
The   above theorem and  two estimates
 are our main results.
 With the   estimates, we will
prove \rt{1} by using a KAM rapid iteration argument as in~\ci{GWzese}, which
avoids  the Nash-Moser smoothing techniques.

\

We now describe some ideas to derive the estimates.
The
integral operators $P,Q$ are estimated in the same way.
Let us
focus on $P =P_0+P_1$, where
 $P_1$ is an integral operator
 over  $\partial M_\rho$
 and    $P_0$ is  over $M_\rho$. Since we need  estimates
 only on shrinking domains, the boundary integral $P_1$ can be treated
 easily. For the interior integral $P_0$, via  cutoff, the difficulties
 lie in the case where the $(0,q)$-form
 $\varphi=\sum\varphi_{\ov J} d\ov {z^J}$
 has compact support. We will see that
  coefficients  of $(0,q-1)$-form $P_0\varphi$ are  sums of
    $\mathcal Kf(z)=\int_{M_\rho}f(\zeta)k(\zeta,z)\, dV(\zeta)$, where
  $f(\zeta)= \varphi_{\ov J}(\zeta) U_{\ov J}(\zeta)$
  and $U_{\ov J}$ depend on   derivatives
  of  $|\zeta'|^2+\hat r(\zeta',\xi^n)$ of order at most two, and
$$
k(\zeta,z)= \f{r_{\zeta^j}-r_{z^j}}
{(r_\zeta\cdot(\zeta-z))^{a}(r_z\cdot(\zeta-z))^{b}},\quad a=n-q,\   b=q
$$
for $r(z)=-y^n+|z'|^2+\hat r(z',x^n)$; see section~\ref{sec3} for details.

To deal with  the singularity of the
 kernel along $\zeta=z$, we fix $z\in M$ and
 apply the approximate Heisenberg transformation $\zeta\to\zs$ defined by
$$
\psi_z\colon \zs'=\zeta'-z',\quad \zs^n=-2ir_z\cdot (\zeta-z).
$$
We will show that   $\zs',\xi_*^n=2\IM (r_z\cdot(\zeta-z))$ form
coordinates of
$\psi_z(M)$, and under this coordinate system  the integral  becomes
$$
\mathcal Kf(z)
=\int_{\pi(\psi_z(M_\rho))}
f(\psi_z^{-1}(\zs))  k_*(\zs',\xi_*^n,z',x^n)\, dV(\zs',\xi_*^n).
$$
Here
\gan
k_*(\zs',\xi_*^n,z',x^n)= 
\sum_{|I|=1}E_I(\zs',\xi_*^n,z',x^n)\hat k_{ab}^I(\zs',\xi_*^n),\\
\hat k_{ab}^I(\zs',\xi_*^n)
=\f{(\zs',\ov{\zs'},\xi_*^n)^I}{(|\zs'|^2+i\xi_*^n)^a
(|\zs'|^2-i\xi_*^n)^b},\quad a+b=n,
\end{gather*}
where $I=(i_1,\ldots, i_{2n-1})$ and we use standard multi-index notation. The
coefficients
$E_I(\zs',\xi_*^n,\cdot)$ are of class $\cC^{a}$  if   $(\zs',\xi_*^n) (\neq0)$
is fixed
and
 $\hat r\in \cC^{a+2}$.
  Since $f$ has compact support,
one can take  derivatives of $\mathcal Kf$ directly onto $f$ and onto $E_I$
without disturbing the   kernels $ \hat k_{ab}^I$.
 The transformation
$\psi_z$ has been used by other people. See
Bruna-Burgu\'es~\ci{BBeise} and Ma-Michel~\ci{MMnith}.
To obtain estimate
(ii), we need to return to the original coordinates after
differentiation. This will give us another  formula for the
derivatives of $\mathcal K f$, which allows us to reduce the
$\cC^{k+\yt}$-estimates to the   H\"older ${\yt}$-estimate for   new kernels of the
same type.

\

We would like  to mention some  methods to derive  the  fundamental $\yt$-estimate.
Kerzman~\ci{Keseon}
 obtained    H\"older  $\alpha$-estimates for all $\alpha<\yt$
  for $\db$-solutions, by estimating
a Cauchy-Fantappi\`e form. Folland-Stein~\ci{FSsefo} used non-isotropic balls
and piecewise smooth curves in  complex tangential directions to obtain
estimates in their 
spaces on the real hyperquadric.
Henkin-Romanov~\ci{HRseon} obtained the $\yt$-estimate for $\db$-equations
on strongly pseudoconvex domains  via a type of
  Hardy-Littlewood lemma (see also Henkin~\ci{Hesese} for $\dbb$ on strictly
  convex boundaries). Our estimate,  like the classical
H\"older estimate for the Newtonian potential,
is still  based on a decomposition of domain. However, we   delete
  a {\it cylinder}
about  the pole,
instead of  a (non-isotropic) ball. The radius of the cylinder is optimized
for the   $\yt$-exponent and  is yet so large  that, when
estimating the H\"older  $ \yt$-ratio at two points, we
 can ignore  their non-isotropic distance and
   connect them  with a line segment.

\

In the $5$ dimensional case there is an extra
term added to the right-hand side of \re{Hhf}.  As of this writing, it remains
unclear whether such a more general homotopy formula can be used.  See~\ci{Weeinia}
and Nagel-Rosay~\ci{NReini}.

The $\yt$-estimate for a  solution to $\db$-equations for $(0,1)$-forms
 on strictly pseudoconvex domains with $\cC^2$
boundary is obtained by Henkin-Romanov~\ci{HRseon}. For
$\cC^{k+\yt}$-estimates for solutions of $\db$-equations of degree
$(0,1)$ on   strictly pseudoconvex domains with $\cC^m$ boundary
($m\geq k+4$), see Siu~\ci{Sisefo}; for
$\db$-equations in
higher degree, see Lieb-Range~\ci{LReize}. The
$\cC^0$-estimate for a solution to  $\dbb$-equations on $M_\rho$,
without shrinking $M_\rho$, is given by Henkin~\ci{Hesese}. For
$\cC^k$-estimates for the homotopy formula
for $\dbb$ operator on shrinking domains,
see~\ci{Weeinia}.
 Michel-Ma~\ci{MMnith} also obtain $\cC^k$-estimates
for a modified homotopy formula  without shrinking
domains, by  introducing an extra   derivative via $\dbb$.

\

We want to mention that in  estimating ({\it i}) and ({\it ii})
we need some  H\"older inequalities. For
the convenience of the reader, we present these   inequalities  in Appendix~A,
 following the formulation and proofs of H\"ormander~\ci{Hosesi}.

The estimates ({\it i}) and ({\it ii}) will be used to improve
regularity in the local CR embedding problem in~\ci{GWzenib}. To limit the scope
of this paper, we leave the estimates in Folland-Stein spaces for
future work.

\setcounter{thm}{0}\setcounter{equation}{0}


\section{ Notation and counting derivatives}
\label{sec2}
To simplify notation, set $z'=(z^1,\ldots,z^{n-1})$, $z=(z', z^n)$, and
$$
x=(\RE z, \IM z') =\pi(z).
$$
Analogously,  $\xi=(\RE\zeta,\IM\zeta')$.
Denote  by
$|\cdot|$ the  euclidean norms on $\cc^{n-1}, \cc^{n-1}\times\rr=\rr^{2n-1}$ and $\cc^n$.
 Our real hypersurface $M\subset\cc^n$ is always a
  graph over  a domain in $\cc^{n-1}\times\rr$. Let
$M_\rho= M\cap\{(x^n)^2+y^n<\rho^2\}$ and
$D_\rho=\pi(M_\rho)$. 
For the real hyperquadric $ M\colon y^n=|z'|^2$,
$\pi(M_\rho)$ is exactly the ball $B_\rho
=\{x\in\rr^{2n-1}\colon |x|<\rho 
\}$.
 On $\rr^{2n-1}$,
 we will use the volume-form
 $dV = d\xi^1\wedge    \cdots
   \wedge d\xi^{2n-1}$. On $\pd D_\rho$, we will
   need $(2n-2)$-forms $dV^s =d\xi^1\wedge\cdots\wedge d\hat\xi^{s}
   \wedge \cdots\wedge d\xi^{2n-1}$.

 \

Let $k\geq0$ be an integer. Denote by $\pd^Iu$  a   derivative of $u$ of order $|I|$, where
$I$ is a standard multi-index. Let
  $\pd^ku$  denote  the set of the $k$-th
order derivatives of $u$.
For a function  $u$    on $D\subset \rr^{2n-1}$, define
\gan
\|\partial^ku\|_{D,0}= \sup_{x\in D,|I|=k} |\partial^Iu(x)|,
\quad \|u\|_{D,k}=\max_{0\leq j\leq k}\|\partial^ju\|_{D,0},\\
|u|_{D,\alpha}= \sup_{x,y\in D}\frac{|u(x)-u(y)|}{|x-y|^\alpha}, \quad 0< \alpha
<1,\\
\|u\|_{D,k+\alpha}=\max\bigl\{
 \|u\|_{D,k}, |\partial^Iu|_{D,\alpha}\colon |I|=k\bigr\}, \quad 0<\alpha<1.
\end{gather*}
If $A=(a_i^j)$ is a  matrix of $\cC^k$ functions on $D$, we define $\|A\|_{D,k}
=\max_{i,j}\{\|a_{i}^j\|_{D,k}\}$. Define $\|A\|_{D,k+\alpha}$ analogously.

  To avoid   confusion and in the essence of estimates
  ({\it i})-({\it ii}) in the introduction,
the $\cC^{a}$-norm of a function  on $M_\rho$ is its
 $\cC^{a}$-norm
on   $D_\rho \subset\rr^{2n-1}$, defined  above.
The $\cC^{a}(M_\rho)$ norm
of a $(0,q)$-form $\varphi=\sum\varphi_{\ov I} d{\ov{z'}}^I$
  is   the maximum of   $\cC^a$-norms
 of  $\varphi_{\ov I}$.
The $\cC^a$-norm on $M_\rho$ will be denoted by
$\|\cdot\|_{\cC^{a}(M_\rho)}=\|\cdot\|_{D_\rho,a}$, or simply by
$\|\cdot\|_{\rho,a}$ when there is no confusion about  $M$.

   Throughout the paper,    $C_k$ denotes a
constant   dependent  of $k$ and
 this dependence will not be expressed sometimes.
Constants, such as $C, C_1, C_k,$
 might have different values when they reoccur.
All constants are independent of  $M, \hat r, \rho, \rho_0$.

\

Let  $M$ be   defined by
\eq{bhroo-}
r(z)\df-y^n+|z'|^2+\hat r(x)=0, \quad x\in D,
\end{equation}
where  $\hat r\in \cC^2(D)$. Throughout the paper,
we make the basic assumption
\ga
\label{bhroo}\ov{D_{\rho_0}}\subset D,
\quad \hat r(0)=\hat
r_z(0)=0, \quad
  \epsilon=\|\hat r\|_{ \rho_0,2}
  <C_0^{-1}.
\end{gather}
We emphasize that   $C_0$  is a large constant
to be adjusted several times. However, it will  not depend on   any
quantity other than $n$.

\

\noindent
{\bf Counting derivatives.}
We need   to
count derivatives  efficiently and   use the count
to estimate norms.
Such a counting scheme  is essentially in~\ci{Hosesi} and we specialize it
 for two reasons. First,  the homotopy formula involves two extra derivatives of
 the defining
function $r$ of $M$; second, for each consecutive $x$-derivative
on    $\varphi_{\ov J}\circ\psi_z^{-1}$, the $x$-derivative which falls on $\varphi_{\ov J}$
yields
an extra factor  $\pd^2  r(z)$ via the chain rule and $\psi_z(\zeta)=(\zeta'-z'
,-i2r_z\cdot(\zeta-z))$.
We illustrate   below
how to use the scheme to 
cope with the two extra derivatives on $r$ and the  consecutive derivatives.

Recall that for $l\geq1$, $\pd^lr(z)=\pd^lr(x) $ is the set of
the $l$-th order derivatives of $r$.
Define
\al\nonumber
\pd_*^1r(\xi,x)&=p\left( \xi,x, (1+\hat r_{x^n}\hat r_{\xi^n})^{-1}, r_{\ov{z^n}}^{-1},
 \pd \hat r(\xi),
 \pd \hat r(x) \right),\\
 \pd_*^2r(\xi,x)&=p\left(\xi,x,(1+\hat r_{x^n}\hat r_{\xi^n})^{-1},r_{\ov{z^n}}^{-1},
  \pd \hat r(\xi),\pd^2\hat r(\xi),
 \pd \hat r(x),\pd^2\hat r(x)\right). 
\label{dsr2}
\end{align}
Here and in what follows,  $p$  is a   polynomial with constant coefficients.
Its coefficients and degree are bounded
in absolute values
by a constant  depending only on  fixed quantities, say $k, n$. Also, $p$
 might be different when
it reoccurs. 
In general, define
$$
\pd_*^{2+k}r(\xi,x)=\sum \pd_*^2r(\xi,x)\pd^{I_1}r(\xi)\cdots\pd^{I_j}r(\xi)
\pd^{J_1}r(x)\cdots\pd^{J_l}r(x),
$$
where the sum is over finitely many multi-indices $I_i,J_i$ satisfying
$$
\sum_{i=1}^j(|I_i|-2)+\sum_{i=1}^{l}(|J_i|-2)\leq k, \quad |I_i|\geq2,\quad |J_i|\geq 2.
$$
We  will  write   $\pd_*^{2+k}r(\xi,x)=\pd_*^{2+k}r(x)$
when it   depends only on $x$.
  With this abbreviation, we have simple relations
\eq{csch}
\pd_*^{2+k}r\pd_*^{2+j}r=\pd_*^{2+k+j}r, \quad \pd^J\pd_*^{2+k}r=\pd_*^{2+k+|J|}r.
\end{equation}
From \ci{Hosesi}, Corollary A.6 (see~\rp{vbhi} in Appendix~A), we
know that with $D_\rho^2=D_\rho\times D_\rho$,
 $$\prod_{j=1}^m\|f_j\|_{D_\rho^2,k_j+b_j} \leq C_{|a|+|c|+m}\rho^{-b_1-\cdots-b_m
 }
 \Bigl(\prod_{j=1}^m\|f_j\|_{D_\rho^2,k_j+a_j}+\prod_{j=1}^m
 \|f_j\|_{D_\rho^2,k_j+c_j}\Bigr)
 $$
for any  non-negative integers $k_j$ and
non-negative real numbers $a_j,c_j$ such that $(b_1,\ldots, b_m)$
is in the convex hull
of $(a_1,\ldots, a_m)$, $(c_1,\ldots, c_m)$.
   With the above abbreviation, basic assumption \re{bhroo}, and
   $0<\rho\leq\rho_0\leq3,$ one  obtains
\eq{csch+}
\|\pd^{2+k}_*r\|_{ {\rho}, a}\leq C_{a+k}\rho^{-a-k}\| r\|_{\rho, 2+k+a },\quad
\| r\|_{\rho,2+k+a}\df 1+\| \hat r\|_{\rho,2+k+a}
\end{equation}
for  all  real numbers $a\geq0$ and integers $k\geq0$.

We will also need a chain rule. Recall that
$\psi_z\colon\zs'=\zeta'-z', \zs^n=2ir_z\cdot(\zeta-z)$ and   define
$$
\Psi(\xi,x)=(\pi\psi_{z}(\zeta),x), \qquad z=\pi|_{M}^{-1}(x),
\quad \zeta=\pi|_{M}^{-1}(\xi).
$$
Let $0<\rho\leq\rho_0\leq3$. We will show that
$B_{\rho/2}\subset D_\rho\subset B_{2\rho}$ and
$$W_\rho=\Psi(D_{\rho}\times D_{\rho})\subset B_{9\rho}\times D_{\rho}.$$
(See Lemmas ~\ref{mdb}-\ref{psii}.)  The Jacobean  matrix
of $\Psi$ depends only on derivatives of $r$ of order $\leq2$
and has  determinant $1+\hat r_{x^n}\hat r_{\xi^n}$ (by \re{-nxisn}).
Then the chain rule takes
   simple forms. Let $(\Psi^{-1})^j$ be the $j$-th   component of $\Psi^{-1}$. Then
\ga\nonumber 
\pd^I\{(\Psi^{-1})^j\}=\pd^2_*r\circ\Psi^{-1}, \quad |I|=1; \\
\label{p1P-}
\pd^K(f\circ\Psi^{-1})=\sum_{|L|\leq |K|}(\pd^{L}f\pd_*^{2+|K|-|L|}r)\circ\Psi^{-1}.
\end{gather}
We will
show that the Lipschitz constant of $\Psi^{-1}$  on $W_\rho$
is bounded by $C$ (see \rp{psii}). Then taking H\"older ratio in \re{p1P-}
with $k=[a]$ gives us
\ga\label{csch++}
 \|f\circ\Psi^{-1}\|_{W_\rho, a} \leq C_k\rho^{-a}
(\|f\|_{D_{ \rho}\times D_\rho,a}+\|f\|_{D_{ \rho}\times D_\rho,0}
\|r\|_{\rho,2+a}).\end{gather}
Note that the above mentioned $\varphi_{\ov J}\circ\psi_z^{-1}$ is a special case.

The above counting scheme via \re{csch}-\re{csch++}
 and its variants will be used systematically.
\setcounter{thm}{0}\setcounter{equation}{0}
\section{ The Henkin homotopy formula}\label{sec3}

In this section
we   recall the \hf\ by following the formulation   in~\ci{Weeinia}.
We discuss the formula  for
differential forms of low regularity.



\

Let $M\subset\cc^n\colon r=0$ be given by \re{bhroo-}-\re{bhroo}.
We first  recall a representative $\dbm$ for the $\dbb$-operator.
By definition, a
$(p,q)$-form $\varphi$ of class $\cC^a$ on $M$ is the restriction
of some
$(p,q)$-form  $\tilde\varphi$ of class $\cC^a$   in a  neighborhood
of $M$. 
If $a\geq 1$, we define $\dbb \varphi$   to
be the restriction of $\db\tilde\varphi$ to $M$.
Notice that on $M$, $\theta=-2i\pd r=\ov\theta$
and that $\db_b\varphi$
is well-defined  modulo $\db r$ when $0\leq p<n$, and it is
actually well-defined when $p=n$.
 By  a {\it tangential}
 $(0,q)$-form $\varphi$, we mean a   form
 $\varphi=\sum_{|I|=q}\varphi_{\ov I}d\ov {z'}^{I}$ with $d\ov{z'}
 =(d\ov {z^1},\dots, d\ov {z^{n-1}})$.
A \cont\ $(0,q)$-form $\varphi$   can be written uniquely as
$\varphi'+\varphi''\wedge\ov\theta$ for some tangential forms
$\varphi',\varphi''$. Define
\gan \ov {X_\alpha}=\partial_{\ov{z^\alpha}}-\f{r_{\ov{z^\alpha}}}
{r_{\ov{z^n}}}\partial_{\ov{z^n}}, \quad
\ov\partial_M\varphi=\dbm\varphi'=\sum_{|I|=q}\ \sum_{1\leq\alpha<n} \ov {X_\alpha}
\varphi'_{\ov I}
d\ov{z^\alpha}\wedge d\ov{z'}^{ I}.
\end{gather*}
Then a straightforward computation shows that
$$
\dbb\varphi=
\ov\partial_M\varphi\mod{\theta}.
$$
Each $(n,q)$-form on $M$ can be written as $\varphi''\wedge dz^1\wedge\cdots
\wedge dz^{n-1}\wedge\theta$  where $\varphi''$ is tangential.
If $\varphi''\in\cC^1$, then
\eq{dbnq}
\dbb(\varphi''\wedge dz^1\wedge\cdots
\wedge dz^{n-1}\wedge\theta)=\dbm
\varphi''\wedge dz^1\wedge\cdots
\wedge dz^{n-1}\wedge\theta.
\end{equation}
Let $\varphi$ be a continuous $(0,q)$-form  on a domain $\mathcal U\subset
M$,  and $\phi$ be
 a continuous $(0,q-1)$-form. Suppose that $q\geq1$.
We say that
 $\dbb\phi=\varphi \mod{\theta}$ holds  on  $\mathcal U$  as currents,
  if $\int_M\phi\wedge\dbb\psi=(-1)^{q}\int_M \varphi\wedge
\psi$ for all $\cC^1$-smooth $(n,n-q-1)$-forms $\psi$ with
compact support in $\mathcal U$.

Write
 $r_z=r_z(z)$ and $r_\zeta=r_\zeta(\zeta).$ Set $N_0(\zeta,z)=r_\zeta\cdot(\zeta-z)$,
 $S_0(\zeta,z)=r_z\cdot(\zeta-z)$ and
\ga
\label{omega+-}
\Omega_{0,q}^{+-}(\zeta,z)=\f{\pd_\zeta r\wedge (r_z \cdot d\zeta)\wedge
(\pd_{\ov\zeta}\pd_\zeta r)^{n-2-q}\wedge (\pd_{\ov z}r_z \wedge d\zeta)^{q-1}}{
 N_0^{n-q}(\zeta,z) S_0^{q}(\zeta,z)},\\
\label{omega0+-}
 \Omega_{0,q}^{0+-}(\zeta,z)=\f{d\zeta^n\wedge\pd_\zeta r\wedge (r_z \cdot d\zeta)\wedge
(\pd_{\ov\zeta}\pd_\zeta r)^{n-3-q}\wedge (\pd_{\ov z}r_z \wedge d\zeta)^{q-1}}{
(\zeta^n-z^n) N_0^{n-q-1}(\zeta,z) S_0^{q}(\zeta,z)}.
\end{gather}
Note that $\theta(\zeta)$   annihilates
$\Omega_{0,q-1}^{+-}(\zeta,z)$ and $\Omega_{0,q}^{0+-}(\zeta,z)$.
Let $n\geq4$ and
  $0<q<n-2$. For  a tangential
  $(0,q)$-form  $\varphi$ on $\ov{M_\rho}$, we have the \hf\
\eq{hfc0}
\varphi=\dbm (P_0'+P_1')\varphi+(Q_0'+Q_1')\dbm\varphi,
\quad z\in M_\rho.
\end{equation}
Here $P_0',P_1',Q_0',Q_1'$ are tangential parts of
\eq{fhc0+}\begin{array}{ll}
P_0\varphi(z)=c_1\int_{M_\rho}\varphi\wedge\Omega_{0,q-1}^{+-}(\zeta,z),\quad
&P_1\varphi(z)=c_2\int_{\partial M_\rho}\varphi\wedge\Omega_{0,q-1}^{0+-}(\zeta,z),
\vspace{.75ex}
\\
Q_0\psi(z)=c_3\int_{M_\rho}\psi\wedge\Omega_{0,q}^{+-}(\zeta,z),\quad
&Q_1\psi(z)=c_4\int_{\partial M_\rho}\psi\wedge\Omega_{0,q}^{0+-}(\zeta,z).
\end{array}
\end{equation}

\

From now on, all $(0,q)$-forms
$\varphi$  on $M$ are tangential.  By $\|\varphi\|_{\rho,a}$,
we mean the norm $\|\varphi'\|_{\rho,a}$ as  defined in section~\ref{sec2},
where $\varphi'=\varphi\mod\theta$ and $\varphi'$  is tangential.
 By $\|P\varphi\|_{\rho,a}$
as used in the introduction, where $P$ is either
of solution operators in the homotopy formula,
we mean $\|P'\varphi\|_{\rho,a}$.
By an abuse of notation, $\varphi$  stands for forms on
$M_\rho$ and $D_\rho=\pi(M_\rho)$.

Next, we describe  kernels of $P'_j,Q'_j$ on domain $D_\rho$ via coordinates $x$.
We have
\ga \label{rzdzeta}
r_\zeta\cdot d\zeta\wedge r_z\cdot d\zeta=\sum_{1\leq j,l\leq n}r_{\zeta^l}(r_{\zeta^j}-r_{z^j})
d\zeta^j\wedge
d\zeta^l,\\ \label{rzdzeta+}
r_\zeta\cdot d\zeta\wedge r_z\cdot d\zeta\wedge d\zeta^n 
=\sum_{1\leq\alpha,\beta<n}
{r_{\zeta^\beta}(r_{\zeta^\alpha}-r_{z^\alpha})}
d\zeta^\alpha\wedge
d\zeta^\beta\wedge d\zeta^n.
\end{gather}
Assume now that $\zeta,z\in M_\rho$. 
We compute $d\ov{\zeta^n}$ and $d \ov{z^n}$ in different ways. We keep the latter
a $(0,1)$-form and find its tangential part.   We have
\ga
\label{dovzn}
d\ov{z^n}=-\f{r_{\ov{z'}}}{r_{\ov{z^n}}}\cdot d\ov{z'}\mod\theta(z),\\
\label{dzetan}
d \zeta^n=(1+i\hat r_{\xi^n})d\xi^n+i2\RE\{   r_{\zeta'}\cdot d\zeta'\}
=2i  r_{ \ov {\zeta^n}}d\xi^n+i2\RE\{  r_{\zeta'}\cdot d\zeta'\}.
\end{gather}
Note that in \re{dzetan}, we have used $r(\zeta)=-\eta^n+|\zeta'|^2+\hat r(\xi)=0$.
In \re{rzdzeta}-\re{rzdzeta+} and $\ov\pd_zr_z\wedge d\zeta$,
we use \re{dovzn}-\re{dzetan} to rewrite $d\ov{z^n}$ and $d\zeta^n$, respectively.
In $\pd_{\ov\zeta}\pd_\zeta r$,  we use \re{dzetan} to rewrite
$d\zeta^n,d\ov{\zeta^n}$.
From \re{omega+-}-\re{omega0+-}, we obtain  on $M_\rho$
   \ga\label{hftk}
  P'_0\varphi(x)=\sum_{|I|=q-1}\ \sum_{|J|=q}\ \sum_{1\leq j\leq n} d\ov{z'}^{  I}
 \int_{D_\rho}
A_{\ov I}^{j\ov J}(\xi,x)
\f{\varphi_{\ov J}(\xi) (r_{\zeta^j}-r_{z^j})}{(N_0^{n-q}S_0^{q})(\zeta,z)} \, dV(\xi),\\
P'_1\varphi(x)= \sum_{|I|=q-1}\sum_{|J|=q}\sum_{ \alpha,\beta=1}^{ n-1}\sum_{s=1}^{2n-1}
d\ov {z'}^I
 \int_{\partial D_\rho}
\f{B^{\alpha\beta \ov J}_{\ov Is}(\xi,x) \varphi_{\ov J}(\xi)(r_{\zeta^\alpha}-
r_{z^\alpha})r_{\zeta^\beta}}{(\zeta^n-z^n)(N_0^{n-q-1}S_0^{q})(\zeta,z)} \,
d V^s(\xi).
\label{hftk+}
\end{gather}
Here $A_{\ov I}^{j\ov J}$ and $ B^{\alpha\beta \ov J}_{\ov Is}$ are polynomials in
$(r_{\zeta },   r_{ \ov\zeta},
 r_{\zeta\ov\zeta},   r_{\ov{z'}}, 1/{r_{\ov{z^n}}},   r_{z\ov z})$.  We make a remark
 for the case $q=1$. In this case we need
 to remove $\sum_{|I|=0}$ in \re{hftk}-\re{hftk+}. Also,
$A^{j\ov J}\df A_{\ov I}^{j\ov J}$ and $B^{\alpha\beta \ov J}_{s}\df  B^{\alpha\beta \ov J}_{\ov Is}$ are independent of $z$.
This observation will play a   role in the
$\yt$-estimate of $P'\varphi$ when $\varphi$ is a $(0,1)$-form.

\

See also Chen-Shaw~\ci{CSzeon}  for   homotopy formulae.   In this paper
we replace the strict  convexity of defining function $r$
in~\cite{CSzeon}
 by the condition \re{bhroo} with $0< \rho_0\leq3$;   see Appendix~B for details.
%
%
%
We  remark that the homotopy formula \re{hfc0} holds
 as currents, when   $\varphi$
and   $\dbb\varphi$   on $M_\rho$
admit continuous extensions to $\ov{M_\rho}$.
See Appendix~B for a proof
by using the Friedrichs approximation theorem.
Henkin~\ci{Hesese}   formulated $\dbb$ in
the   sense of currents.
  Shaw~\ci{Shnize}
also used  the Friedrichs approximation for $\dbb$-solutions.

%
%
%
%
%
%
%

\setcounter{thm}{0}\setcounter{equation}{0}\section{
Kernels and the approximate Heisenberg transformation
}\label{sec4}

  In this section, we will describe briefly the new kernels when the approximate
  Heisenberg transformation is applied.
  The
   contents of next few sections are indicated at the end of this section.

Recall that 
in \re{hftk}  functions $A_{\ov I}^{j \ov J}$ have
 the form $\pd_*^2r$. Hence coefficients of $P'_0\varphi$ are
sums over $|J|=q$ and $1\leq j\leq n$ of 
$
\mathcal  K\varphi_{\ov J} (x)=\int_{D_\rho}\varphi_{\ov J} ( \xi)   k_{ab}^j(\xi,x)\,dV(\xi).  
$
Here  $\zeta,z\in M$ and
$$
k(\xi,x)=k_{ab}^j(\xi,x)=\f{\pd_*^2r(\xi,x)
 (r_{\zeta^j}-r_{z^j})}{(r_\zeta\cdot(\zeta-z))^{a}(r_z\cdot(\zeta-z))^{b}},
  \quad  a=n-q, \  b=q.
$$

\

To understand the kernel, let us compute its denominator.
Set   $\hat r(z)=\hat r(x)$.
On $M\times M$ we have
\aln
r_z\cdot(\zeta-z)&=\ov{z'}\cdot(\zeta'-z')+\f{i}{2}(\zeta^n-z^n)
+\hat r_z\cdot(\zeta-z)\\
&=\ov{z'}\cdot(\zeta'-z')-\yt(|\zeta'|^2-|z'|^2)+\f{i}{2}(\xi^n-x^n)\\ &\quad
-\yt(\hat r(\zeta)-\hat r(z))+\hat r_z\cdot(\zeta-z)\\
&=i\IM(r_z\cdot(\zeta-z))-\yt|\zeta'-z'|^2\\ &\quad-
\yt(\hat r(\zeta)-\hat r(z))+\RE(\hat r_z\cdot(\zeta-z)).\intertext{Also} 
r_\zeta\cdot(\zeta-z)&=r_z\cdot(\zeta-z)
+(r_\zeta-r_z)\cdot(\zeta-z)\\
&=i\IM(r_z\cdot(\zeta-z))+\yt|\zeta'-z'|^2
+(\hat r_\zeta-\hat r_z)\cdot(\zeta-z)\\ &\quad -
\yt\bigl\{\hat r(\zeta)-\hat r(z)-2\RE(\hat r_z\cdot(\zeta-z))\bigr\}.
\end{align*}
The first two terms in  $r_z\cdot(\zeta-z)$ and $r_\zeta\cdot(\zeta-z)$
can be simplified simultaneously.
Define
\eq{dfnz}
 N(\zeta,z)\df |\zeta'-z'|^2+2i\IM(r_z\cdot(\zeta-z)).
\end{equation}
Then we arrive at the basic relations
\ga
-2r_z\cdot(\zeta-z)=\ov{ N(\zeta,z)} + A(\zeta,z),
\nonumber 
\quad
\nonumber
2r_\zeta\cdot(\zeta-z)=N(\zeta,z) +B(\zeta,z)\\
\intertext{with}
A(\zeta,z)=
\hat r(\zeta)-\hat r(z)-2\RE(\hat r_z\cdot(\zeta-z)),\quad  
B(\zeta,z)=2(\hat r_\zeta-\hat r_z)\cdot(\zeta-z)- A(\zeta, z).\nonumber
\end{gather}
By   the condition \re{bhroo} on $\hat r$, we  will show that $D_{\rho_0}$
is convex and hence
$$
|A(\zeta,z)|\leq C\epsilon|\zeta-z|^2, \quad |B(\zeta,z)|\leq C\epsilon|\zeta-z|^2,\quad
\zeta,z\in M_{\rho_0}.
$$
We will show
 that $|N(\zeta,z)|\geq C^{-1}|\zeta-z|^2$  on
 $M_{\rho_0}\times M_{\rho_0}$. Consequently, the kernel
of $P_0'\varphi$ 
is factored as
\gan k(\xi,x)=
T_1(\xi,x)^{-a}T_2(\xi,x)^{-b}
\f{\pd_*^2r(\xi,x)(r_{\zeta^j}-r_{z^j})}{ N^{a}(\zeta,z)\ov {N^{b}(\zeta,z)} },
\\
 T_1(\xi,x)= 1+N^{-1}(\zeta,z)B(\zeta,z), \quad  T_2(\xi,x)=1+\ov {N^{-1}(\zeta,z)}A(\zeta,z)
\end{gather*}
for $\zeta, z\in M_{\rho_0}$.
Moreover, $\min\{|T_1(\xi,x)|,|T_2(\xi,x)|\}\geq
1/4$ when $(\xi,x)$ is  on  $D_{\rho_0}\times D_{\rho_0}$ and off
its diagonal.
 Now,   identity \re{dfnz} suggests
the following approximate Heisenberg transformation
$$
\psi_z\colon\zs'=\zeta'-z', \quad \zs^n=-2ir_z\cdot(\zeta-z).
$$
We will show that for fixed $z\in M_{\rho_0}$,
$\zs'=\zeta'-z',\xi_*^n=2\IM(r_z\cdot(\zeta-z))$ indeed  form
  coordinates of $M_{\rho_0}$. Thus, with $\zeta, z\in M_{\rho_0}$
\gan
N(\zeta,z)=|\zs'|^2+i\xi_*^n 
\df N_*(\xis), \\
k(\xi,x)=\sum_{|I|=1}E_{I}(\xi_*,x)
{\xi_*^{I}}{N_*(\xis) ^{-a}\ov {N_*^{-b}(\xis)}},\quad a+b=n.
\end{gather*}
When $\varphi$ has compact support, the decomposition for $k(\xi,x)$
 allows us to take
$x$-derivatives of $K\varphi_{\ov J}(x)$ directly onto  coefficients
$E_I(\xi_*,x)$ and onto $\varphi_{\ov J}(\psi_z^{-1}(\xi_*))$, without
destroying the integrability of the kernels.


\

We now describe the contents of next few sections.

We will carry out the details
of this section  in sections~\ref{sec5},~\ref{sec6}, and~\ref{sec8-}.
In section~\ref{sec5} we will also study
how domains $B_{\rho}$, $D_\rho=\pi(M_\rho)$ and  $\pi\psi_z(M_\rho)$     are nested.
We will need  H\"older inequalities
in Appendix~A for   domains $D_\rho$.
Therefore, we will
 verify that  under the basic assumption \re{bhroo} and $0<\rho\leq{\rho_0}\leq3$,
 $D_\rho$ is convex and $B_{\rho/2}\subset D_\rho\subset B_{2\rho}$.
In section~\ref{sec6},
we will also express the graph $\psi_z(M_{\rho_0})$ for   $z\in M_{\rho_0}$ as
$$
\eta_*^n=\sum_{|I|=2} h_I( \xi_*,x) \xi_*^I, \quad \xi_*=(\RE\zs,\IM\zs').
$$
The latter will be used to show $|r_z\cdot(\zeta-z)|\geq C^{-1}(\rho\sigma)^2$
for $z\in M_{(1-\sigma)\rho}$ and $\zeta\in \pd M_\rho$.
In sections~\ref{sec8} and \ref{sec9}, we derive   the  $\cC^a$-estimates. The
$ \yt$-estimate is  in section~\ref{sec14},  following a reduction for $\cC^{k+\yt}$-estimates
in section~\ref{sec10}.

We conclude this section with a lower bound for $|\zeta^n-z^n|$.
By the basic assumption \re{bhroo}, $D_{\rho_0}$
is relatively  compact in $D$   and hence $(x^n)^2+y^n=\rho^2$
on $\pd M_\rho$ for $0<\rho\leq\rho_0$.
Then the image of the projection of $\pd M_{\rho}$ in the
$z^n$-plane is contained
in the parabola $(x^n)^2+y^n=\rho^2$. Therefore,
for $z\in M_{(1-\sigma)\rho}$ and $\zeta\in\partial M_\rho$
 with $\rho\leq\rho_0 (\leq3)$, we obtain
\eq{nsrr++-}
 |\zeta^n-z^n|\geq C^{-1}\rho^2\sigma.
\end{equation}

\setcounter{thm}{0}\setcounter{equation}{0}\section{ Domains and images
under the transformation}
\label{sec5}

Recall that our real hypersurface $M$ is given by \re{bhroo-}-\re{bhroo}, and
 $$M_\rho=M\cap\{(x^n)^2+y^n<\rho^2\},\quad D_\rho=\pi(M_\rho)=\{x\in
 D\colon|x|^2+\hat r(x)<\rho^2\}.$$

 \le{mdb} 
 Let $M$ satisfy \rea{bhroo}. 
Suppose that $0<\rho\leq\rho_0$. Then $\ov{D_\rho}$ is strictly convex with $\cC^2$ boundary.
Also,   $$
 B_{(1-c\epsilon)\rho}
\subset D_{\rho}\subset
B_{(1+c\epsilon)\rho},\quad C^{-1}\rho\sigma\leq
\dist(\pd D_{(1-\sigma)\rho},\pd D_\rho)\leq C\rho\sigma.$$
 Moreover,  constants $C_0$, $c, C$   are
 independent of $\rho$ and $\hat r$.
 \end{lemma}
\pf  Let $\epsilon=\|\hat r\|_{\rho_0,2}<C_0^{-1}$.
The strict convexity   follows from the positivity
of the   Hessian of $\phi(x)=(x^n)^2+y^n=
 |x|^2+\hat r(x)$ on $\ov {D_\rho}$. Since $0\in D_\rho$,
then
  $|\hat r(x)|\leq C\epsilon|x|^2$ on $D_\rho$. Now,
$$
(1-C\epsilon)|x|^2\leq \phi(x) 
\leq
(1+C\epsilon)|x|^2.
$$
In particular, for a possibly  larger $C$, we have $
B_{(1-C\epsilon)^{1/2}\rho}\subset D_\rho\subset B_{(1+C\epsilon)^{1/2}\rho}$,
 since  $\ov{D_{\rho_0}}\subset D$ implies that
   $\phi=\rho$ on $\pd D_\rho$ for all $0<\rho\leq\rho_0$.

Let   
$x\in\pd D_{(1-\sigma)\rho}$.
Then $\phi(x)=((1-\sigma)\rho)^2$ and $\yt(1-\sigma)\rho<|x|<2(1-\sigma)\rho$.
For $y\in D_\rho$, we have $
|\phi(y)-\phi(x)|\leq  C\|\pd^1\phi\|_{\rho,0}|y-x|$.  We get
$$
|\phi(y)|\leq (1-\sigma)^2\rho^2+C\|\pd^1\phi\|_{\rho,0}|y-x|\leq (1-\sigma)^2\rho^2
+C\rho|y-x|.
$$
Hence,  $\phi(y)<\rho$ for $|y-x|\leq C^{-1}\rho\sigma$. This shows that
$\dist(\pd D_{(1-\sigma)\rho},\pd D_\rho)\geq C^{-1}\rho\sigma$. On the other hand,
if $y=(1+t)x\in D_\rho$ and $t>0$,   applying
the mean-value-theorem to $\phi((1+s)x)$ for $0\leq s\leq t$ yields
\aln
\phi(y)&\geq ((1-\sigma)\rho)^2+  x\cdot (y-x)-C\epsilon|y||y-x|
\\
&\geq \rho^2(1-\sigma)^2+t|x|^2-C'\epsilon\rho t|x|\\
&\geq
\rho^2\{(1-\sigma)^2+\f{1}{4} t(1-\sigma)^2 -2C''\epsilon t\}.
\end{align*}
This show that $\phi((1+t)x)>\rho^2$ if $(1+t)x\in D_\rho$ and $t>C\sigma$, a
contradiction. Therefore, $\dist(\pd D_{(1-\sigma)\rho},\pd D_\rho)\leq C\rho\sigma$.
\end{proof}


Recall the approximate Heisenberg transformation
$$
\psi_z\colon\zs'=\zeta'-z', \quad \zs^n=-2ir_z\cdot(\zeta-z).
$$
Before applying $\psi_z$ to kernels $k_{a b}^j(\zeta,z)$, we need to know
 how it  transforms  
 $M_\rho$.
Recall   our notation  $x\in\rr^{2n-1}$ and $\xi\in\rr^{2n-1}$.
 Define a map $\tilde\psi_{x}$
by relations
$$\xi_*=\tilde\psi_{x}(\xi)=\pi\psi_z(\zeta), \quad
 \zeta, z\in M_{\rho_0}.
$$
Set $\Psi(\xi,x)=(\tilde\psi_{x}(\xi),x)$.
We have  the following.
\le{psii}
Let $M\colon y^n=|z'|^2+\hat r(z',x^n)$ satisfy \rea{bhroo} with $0<\rho_0\leq3$.
There exist constants
$C, C_m$, independent of   $\hat r$, $\rho$ and $\rho_0$,
 such that the following hold.
\bppp
 \item If $x\in D_\rho$ and $0< \rho\leq\rho_0$ then
$$
 \tilde\psi_{x}(D_{\rho})\subset   B_{9\rho},\quad
  \Psi(D_{\rho}\times D_\rho)\subset   B_{ 9\rho}\times D_\rho.
$$
\item For $u,v\in D_{\rho_0}\times D_{\rho_0}$,
$$
C^{-1}|v-u|\leq|\Psi(v)-\Psi(u)|\leq C|v-u|.
$$
If particular, if $z\in M_{\rho_0}$ then
 $\psi_z(M_{\rho_0})$ is  a graph over $\pi\psi_z(M_{\rho_0})$.
 \eppp
\end{lemma}
\pf   Let  $R(x)=|z'|^2+\hat r(x)$. By \re{bhroo},
$\ov {D_{\rho_0}}\subset D$
and $\epsilon=\|\hat r\|_{\rho_0,2}<C_0^{-1}$.
For brevity, set $M=M_{\rho_0}$.


(i). Assume that $0<\rho\leq\rho_0$. Since  $D_{\rho_0}$ is convex,
we have $|\hat r(x)|\leq C\epsilon|x|^2$ on $D_{\rho_0}$.
The map $\tilde\psi_{x}$ is defined by $\zs'=\zeta'-z'$ and
\al\label{xisn}
 \xi_*^n&=
\xi^n-x^n+\hat r_{x^n}(R(\xi)-R(x))+2\IM [R_{z'}\cdot(\zeta'-z')]\\
&=\xi^n-x^n+2\IM{(\ov{z'}\cdot \zeta')}+\hat r_{x^n}(R(\xi)-R(x))+\nonumber
2\IM [\hat r_{z'}\cdot(\zeta'-z')].
\end{align}
 Let
$\xi,x$ be in $D_{\rho}$.
Recall that $D_{\rho}\subset B_{(1+C\epsilon)\rho}$. We have
\gan
|x|<(1+C\epsilon)\rho, \quad |\xi|<
(1+C\epsilon){\rho},\quad
|\hat r_{x^n}|\leq \epsilon,  \quad|\hat r_{z'}|\leq\epsilon,\\
|R(\xi)|\leq |\zeta'|^2+C\epsilon|\xi|^2
<(1+C'\epsilon)\rho^2.
\end{gather*}
Thus, $|\zs'|^2=|\zeta'-z'|^2\leq 2|\zeta'|^2+2|z'|^2<4(1+C\epsilon)\rho^2$; by \re{xisn},
$
|\xi_*^n|^2\leq   (2\rho+2\rho^2+C\epsilon\rho)^2$. Since $2\rho^2\leq 6\rho$ then
$$|\zs'|^2+|\xi_*^n|^2\leq4(1+C\epsilon)\rho^2+(2\rho+2\rho^2+C\epsilon\rho)^2
<(  9\rho)^2,$$
if $\epsilon$ is sufficiently small.
We get $\tilde\psi_{x}(D_{\rho})
\subset  B_{9\rho}$.

 (ii).  
 It is obvious that the Lipschitz constant of
 $\Psi$ on the convex domain
 $D_{\rho_0}\times D_{\rho_0}$ is bounded by some $C$.
Let $\xi,x,\tilde\xi,\tilde x$ be in $D_{\rho_0}$.  We need
 to show that $|\Psi(\tilde\xi,\tilde x)-\Psi(\xi,x)|\geq C^{-1}|
 (\tilde\xi,\tilde x)-(\xi,x)|$.
Write $\xi=(\zeta',\xi^n), x=(z',x^n)$. Then
$$
|\Psi( \tilde\xi,\tilde x )-\Psi(\xi,x)|\geq|(\tilde\zeta'-
\tilde z'-\zeta'+z',\tilde x-x)|
\geq|(\tilde\zeta'-\zeta',\tilde x-x)|/C.
$$
Set $\tilde\xi_*=\tilde\psi_{ \tilde x}(\tilde
 \xi)$ and
$\xis=\tilde\psi_{x}(
\xi)$.
It suffices to show that
$$
|\tilde\xi_*^n-\xi_*^n|\geq  |\tilde\xi^n-\xi^n|/C
$$
if $|(\tilde\zeta'-\zeta',\tilde x-x)|<\f{1}{48}
|\tilde\xi^n-\xi^n|$. Assume that the latter holds.
Recall that $0<\rho_0\leq3$ and
 $D_{\rho_0}\subset B_{2\rho_0}$. We have $\max\{|z'|,|\tilde
\zeta'|\}<2\rho_0$. Now the second identity in \re{xisn} implies that
\aln
|\tilde\xi^n_*-\xi_*^n|&\geq|\tilde\xi^n-\xi^n|-
|\tilde x^n-x^n|-2|\tilde\zeta'||\tilde z'-z'|
-2|z'||\tilde\zeta'-\zeta'| 
-C\epsilon|(\tilde\xi-\xi,\tilde x-x)|\\
&\geq (1-\f{1}{48})|\tilde \xi^n-\xi^n|- {8\rho_0}(|\tilde\zeta'-\zeta'|+|\tilde z'-z'|)
-C'\epsilon|\tilde \xi^n-\xi^n|.
\end{align*}
Thus, $|\tilde\xi^n_*-\xi_*^n|\geq (1-\f{1}{48}-\f{8\cdot 3\cdot\sqrt2}{48}-C'\epsilon)
|\tilde \xi^n-\xi^n|\geq |\tilde \xi^n-\xi^n|/4$.

That $\psi_z(M_{\rho_0})$ is a graph
 follows from the injectivity of $\Psi$.
 \end{proof}

\setcounter{thm}{0}\setcounter{equation}{0}\section{ 
$r_z\cdot(\zeta-z), r_\zeta-r_z$ and $\psi_z(M)$   via   Taylor's theorem
}
\label{sec6}

To transform  $k(\zeta,z)$ via $\psi_z$,  we need
  expansions of $r_z\cdot(\zeta-z),r_{\zeta}\cdot
(\zeta-z)$ and $r_{\zeta^j}-r_{z^j}$ in   new variables
$ \xi_* $. We will find these expansions via Taylor's theorem.

Let us recall    Taylor's theorem.
Assume that $0<\rho\leq\rho_0\leq3$ and
$\hat r$ satisfies  \re{bhroo}.
 So $D_{\rho}$ is strictly convex.
If $f$ is a complex-valued function  on the convex set $D_{\rho}$, 
we define $\tr{k}f$ and $\tr{I}f$ on $ D_\rho\times D_\rho$ by
\aln
\tr{k} f(y,x)&\equiv f(y)-
\sum_{0\leq j\leq k-1}\frac{1}{j!}\pd_t^j\vert_{t=0}f(x+t(y-x))\\
&= \frac{1}{(k-1)!}\int_0^1(1-t)^{k-1}\pd_t^kf(x+t(y-x))\, dt=\sum_{|I|=k}\tr{I}f(y,x)(y-x)^I.
\end{align*}
Denote by $\mathcal R^kf$ the set of remainder
coefficients $\tr{I}f$ with $|I|=k$.
For any  real number $a\geq0$, 
\eq{estr}
\|\tr{I}f\|_{ D_\rho\times D_\rho,a}\leq C_{a+|I|}\|f\|_{\rho,a+|I| }.
\end{equation}
We will also need to express the remainders in  $\xi_*=(\RE \zs,\IM\zs')$.
Let $z\in M_{\rho_0}$ and $\pi(z)=x\in D_{\rho_0}$.
By  \rl{psii}, $\psi_z(M_{\rho_0})$ is  a graph
 $$
 \eta_*^n=h(\xis,x),\qquad  \xi_* \in\tilde\psi_{x}(D_{\rho_0}).
 $$
By $\zeta'-z'=\zs'$
and $\zs^n=-2i r_z\cdot(\zeta-z)$, we have
 \aln
\zs^n&=\zeta^n-z^n-i\hat r_{x^n}(\zeta^n-z^n) -2i R_{z'}\cdot\zs'.
\end{align*}
 Computing the real and
imaginary parts, we get
 \al
\label{-nxisn}
 \xi_*^n&=\xi^n-x^n+\hat
r_{x^n}(R( \xi )-R(x))+2\IM
(R_{z'}\cdot\zs'),\\
\label{-nxisn+}
\eta_*^n
&=|\zs'|^2+\hat r( \xi )-\hat r(x)-\hat r_{x^n}(\xi^n-x^n)-2\RE
(\hat r_{z'}\cdot\zs').
\end{align}
In \re{-nxisn}, replace   $R(\xi)-R(x)$ by $$
\sum_{|I|=1} \tr{I}R(\xi,x)(\xi-x)^I=\sum_{|I|=1} \tr{I}R(\xi,x)(\xi_*',\xi^n-x^n)^I$$ and  then solve
for  $\xi^n-x^n$.  We get
\eq{xinx}
\xi^n-x^n=
\sum_{|I|=1}p_I(\xi,x)\xi_*^I,
\end{equation}
where $p_I(\xi,x)$  are   of the form
$$
p\Bigl(\xi,x, \f{1}{1+\hat r_{x^n}\tr{(0',1)}\hat r(\xi,x)},\pd\hat r(x),
\mathcal R^1\hat r(\xi,x)\Bigr)
$$
for some polynomials $p$.
With these   polynomials $p$, introduce notation
\eq{cont}
\trd_*^i g=\sum_{|I |=i}p\Bigl(\xi,x, \f{1}{1+\hat r_{x^n}\tr{(0',1)}
\hat r(\xi,x)},\pd\hat r(x),
\mathcal R^1\hat r(\xi,x)\Bigr)\tr{I }g(\xi,x).
\end{equation}
Note that 
 reappearing
$p,\trd_*^ig$ may
be different.
We now   express Taylor remainders in variables $\zs',\xi_*^n$ as follows
\ga\label{pa}
\tr{k}f(\xi,x)=\sum_{|L|=k}\tr{L}
f(\xi,x)(\xi-x)^L
=\sum_{|L|=k}\trd_*^kf(\xi,x)\xi_*^L,\\
|\trd_*^l\hat r(\xi,x)|\leq C\|
\hat   r\|_{\rho_0,2},\quad l\leq2,\quad \xi,x\in D_{\rho_0}.
\nonumber\end{gather}

We now apply notation \re{cont}.
Using  \re{xinx} in \re{-nxisn+}
 we get 
 \ga\label{ahsord1}
\eta_*^n= |\zs'|^2+\sum_{|I|=1}
  \trd_*^1\hat r(\xi,x) \xi_* ^I,\\
 \label{ahsord2}
\eta_*^n 
=|\zs'|^2+\tr{2}\hat
r(\xi,x)
=|\zs'|^2+
\sum_{|I|=2}\trd_*^2\hat r(\xi,x) \xi_* ^I.
\end{gather}
Set $
r(z)=-y^n+|z'|^2+\hat r(z',x^n).$ 
We have defined
\gan 
N(\zeta,z)=|\zeta'-z'|^2+2i\IM(r_z\cdot(\zeta-z)),\quad
N_*( \xi_* )=|\zs'|^2+i\xi_*^n. 
\end{gather*}
Recall that   for $(\zeta,z)\in M_{\rho_0} \times M_{\rho_0}  $ we have
 \ga
-2r_z\cdot(\zeta-z)=\ov{N(\zeta,z)}+  \tr{2}\hat r(\xi,x), \label{2rz}\quad
2r_\zeta\cdot(\zeta-z)=N(\zeta,z)+B(\zeta,z),\\
B(\zeta,z)=2(\hat r_\zeta-\hat r_z)\cdot(\zeta-z)- 
 \tr{2}\hat r(\xi,x).\label{bzz2}
\end{gather}
Define 
\eq{trpdr}
 \trd_*^i\pd \hat r=\sum_{1\leq j\leq n} \trd_*^i\hat r_{x_j}+
 \sum_{1\leq j<n} \trd_*^i\hat r_{y_j}.
  \end{equation}
We  have
\ga (\hat r_\zeta-\hat r_z)\cdot(\zeta-z)=\sum_{|K|=2 } \trd_*^1
\pd\hat r(\xi,x)  \xi_*^{K},\nonumber\\
    \label{Bzz} B(\zeta,z)= \sum_{|L|=2}
(\trd_*^2\hat r(\xi,x)+
  \trd_*^1 \pd \hat r(\xi,x))  \xi_*^L.
\end{gather}
For the numerator  of the kernel, we have
\al\label{rzzj}
r_{\zeta}-r_{z}=(\zs',0)  +\sum_{|I|=1} (\trd_*^1\pd \hat r, \ldots,\trd_*^1\pd \hat r)(\xi,x)
 \xi_* ^I.
\end{align}


In summary, we have proved the following expansions.
\le{m2h}
Let $M\colon y^n=|z'|^2+\hat r(z',x^n)$ satisfy \rea{bhroo} with $0< \rho_0
\leq3$.
Suppose that $\zeta, z\in M_{\rho_0}$ and
$\zs=\psi_z(\zeta)$.
 Then $\psi_z(M_{\rho_0})$ is given by  $ \eta_*^n
=|\zs'|^2+h( \xi_* ,x)$. Moreover,
\aln
h( \xi_*,x)&=\sum_{|I|=1}\trd_*^1\hat r(\xi,x)
 \xi_* ^I
= \sum_{|I|=2}\trd_*^2\hat r(\xi,x) \xi_*^I,\\
-2\rzd&=|\zs'|^2-i\xi_*^n+
\sum_{|I|=2}\trd_*^2\hat r(\xi,x) \xi_* ^I, \\
2\rzed&=|\zs'|^2+i\xi_*^n
+\sum_{|I|=2 }(\trd_*^2\hat r(\xi,x) +\trd_*^1\pd \hat r(\xi,x))\xi_* ^I,
\\
r_{\zeta}- r_{z}
&=(\zs',0)+\sum_{|I|=1}(\trd_*^1\hat r_{z^1},\ldots, \trd_*^1\hat r_{z^n})
(\xi,x)
  \xi_*^I.
\end{align*}
\end{lemma}
We emphasize that $\trd_*^ig, \trd_*^i\pd\hat r$,  defined by  \re{cont} and \re{trpdr}, might be
different when they reoccur.

\begin{rem}\label{62}
Notice that \re{-nxisn}-\re{ahsord2} are valid  if we fix $\zeta\in M_{\rho_0}$
and vary $z\in M_{\rho_0}$. Therefore, the image of $M_{\rho_0}$ under the map
$z\to\psi_z(\zeta)$ is still given by \re{ahsord1} and \re{ahsord2}, where
 $z$ varies in $M_{\rho_0}$.
\end{rem}

Here are immediate consequences of the above
expansions: 
\ga\label{nr}
\Bigl|\f{2r_\zeta\cdot(\zeta-z)}{N(\zeta,z)}-1\Bigr|<\yt, \quad
\Bigl|\f{2r_z\cdot(\zeta-z)}{\ov{N(\zeta,z)}}+1\Bigr|<\yt,
\quad \f{|N(z,\zeta)|}{|N(\zeta,z)|}<4
\end{gather}
for $\zeta,z\in M_{\rho_0}$ with $\zeta\neq z$ and $\rho_0\leq3$.

\le{distb} Let $M$ satisfy \rea{bhroo} and let $0<\rho\leq\rho_0\leq3$.
Set $d(\zeta,z)=|r_z\cdot(\zeta-z)|$.
Then
 \gan
C^{-1}|\zeta-z|^2\leq d(\zeta,z) \leq C(d(\zeta,v)
+d(v,z)),\qquad \zeta,z,v\in M_{\rho_0};\\
d(\zeta,z)\geq C^{-1}\rho^2\sigma^2, \quad
z\in M_{(1-\sigma)\rho},\  \zeta\in \partial M_{\rho},\ 0<\sigma<1.
\end{gather*}
\end{lemma}
\pf  $M$ is defined by $r=- y^n+|z'|^2+\hat r(z',x^n)=0$ and
$\epsilon=\|\hat r\|_{\rho_0,2}<C_0^{-1}$.
Let $\zeta,z, v$ be in $M_{\rho_0}$. Set
$\zs'=\zeta'-z',
\zs^n=-2ir_z\cdot(\zeta-z)$. Recall that $D_{\rho_0}$ is convex.

(i). By definition,
$$
N(\zeta,z)=|\zeta'-z'|^2+2i\IM(r_z\cdot(\zeta-z))=|\zs'|^2+i\xi_*^n.
$$
 By \rl{m2h}, we have
$$
\psi_zM_{\rho_0}\colon \eta_*^n=|\zs'|^2+
\sum_{|I|=2} \trd_*^2\hat r(\xi,x) \xi_* ^I.
$$
  On $D_{\rho_0}\times D_{\rho_0}$, $|\trd_*^2\hat r|\leq
C\epsilon$, and $| \xi_* |<C$.  Then
$|\zs^n|^2=|\xi_*^n|^2+ |\eta_*^n|^2\leq
C_0' | \xi_* |^2$. Therefore,
$|\zeta^n-z^n|^2=|\f{i}{2r_{z^n}}\zs^n-\f{r_z'}{r_{z^n}}\cdot
\zs'|^2\leq C|N(\zeta,z) |$. Also, $|\zeta'-z'|^2\leq |N(\zeta,z) |$.  We
conclude that $|N(\zeta,z) |\geq C^{-1}|\zeta-z|^2$ for
$\zeta,z\in M_{\rho_0}$. By \re{nr} we have $d(\zeta,z)\geq|\zeta-z|^2/C$.
 Now
\aln |r_\zeta\cdot(\zeta-z)|&\leq
|(r_\zeta-r_v)\cdot (z-v)|+
|r_v\cdot(z-v)|+|r_\zeta\cdot(v-\zeta)|\\
&\leq|r_\zeta-r_v|^2+|z-v|^2+d(z,v)+d(v,\zeta) \leq C(d(\zeta,v)+d(v,z)).
\end{align*}
Thus $d(\zeta,z)\leq C'(d(\zeta,v)+d(z,v))$.


(ii).  By \rl{mdb}, $\dist(\pd
D_{(1-\sigma)\rho},\pd D_\rho)\geq C^{-1}\rho\sigma$.
Then   $d(\zeta,z)\geq C^{-1}(\rho\sigma)^2$ follows from
$d(\zeta,z)\geq
C^{-1}|\zeta-z|^2$.
\end{proof}

\setcounter{thm}{0}\setcounter{equation}{0}
\section{ Outline of $\cC^a$ 
estimates}
\label{sec7}

Let $M$ satisfy \re{bhroo-}-\re{bhroo} and let $0<\rho\leq\rho_0\leq3$.
Recall from \re{hftk}-\re{hftk+} that
  \gan
  P'_0\varphi(x)=\sum_{|I|=q-1}\ \ \sum_{|J|=q}\ \ \sum_{1\leq j\leq n} d\ov{z'}^I
 \int_{D_\rho}
A_{\ov I}^{j\ov J}(\xi,x)
\f{\varphi_{\ov J}(\xi) (r_{\zeta^j}-r_{z^j})}{(N_0^{n-q}S_0^{q})(\zeta,z)} \, dV(\xi),\\
P'_1\varphi(x)= \sum_{|I|=q-1}\, \sum_{|J|=q}\, \sum_{ \alpha,\beta=1}^{n-1 }\,  \sum_{s=1}^{2n-1}
d\ov{z'}^I
 \int_{\partial D_\rho}
\f{B^{\alpha\beta \ov J}_{\ov Is}(\xi,x) \varphi_{\ov J}(\xi)(r_{\zeta^\alpha}-
r_{z^\alpha})r_{\zeta^\beta}}{(\zeta^n-z^n)(N_0^{n-q-1}S_0^{q})(\zeta,z)} \,
d V^s(\xi).
\end{gather*}
Here $\zeta, z\in M_\rho$, and
$A_{\ov I}^{j\ov J}$, $ B^{\alpha\beta \ov J}_{\ov Is}$ are polynomials in
$  (r_{\zeta },   r_{ \ov\zeta},
 r_{\zeta\ov\zeta},   r_{\ov{z'}}, r_{\ov{z^n}}^{-1},   r_{z\ov z})$.

\

We emphasize
that norms are defined on $D_\rho=\pi(M_\rho)$
via coordinates $x$. Let us indicate how to obtain  $\cC^a$
estimates. We will give estimates on shrinking domains $M_{(1-\sigma)\rho}$.
For $\zeta\in\partial M_\rho$ and $z\in M_{(1-\sigma)\rho }$,
we   obtain $\min\{|r_\zeta\cdot(\zeta-z)|,
|r_z\cdot(\zeta-z)|\}\geq C^{-1}(\rho\sigma)^2$
by \rl{distb}
and
$|\zeta^n-z^n|\geq C^{-1}\rho^2\sigma$ by~\re{nsrr++-}.
 This will allow us to estimate the $\cC^a$-norm of
 boundary integral
$P'_1\varphi$ by passing derivatives over the integral sign and
differentiating  the kernels directly.

We now deal with the interior integrals $P'_0\varphi$.
Using
  a partition of unity, we can find  a smooth function $\chi=\chi_{\sigma,\rho}$,
  which is $1$  on  $D_{\rho(1-\sigma/2)}$
and zero  off $D_{\rho (1-\sigma/4)}$, such that $\|\chi\|_{\rho,a}\leq
C_a(\rho\sigma)^{-a}$.
On $D_{\rho}$,  decompose
\gan
\varphi_0=\chi \varphi, \quad
\varphi_1=\varphi-\varphi_0,\qquad
P'_0\varphi=P'_0\varphi_0+P'_0\varphi_1.
\end{gather*}
Now $P'_0\varphi_1$ can be estimated
on $M_{(1-\sigma)\rho}$ by differentiating the kernels directly, since $\varphi_1$
 is supported in $M_\rho\setminus M_{(1-\yt\sigma)\rho}$.
The only non-trivial integral is   $P'_0\varphi_0$, for which $\varphi_0$ has
  compact support
in $M_{(1-\yf\sigma)\rho}$.
To estimate the latter,  
we will apply  the transformation $\psi_z$ for the integral and
then differentiate the new integral.

The estimate for $P'_0\varphi_0$  
is the  most technical part.
We  deal with  this estimate first in sections~\ref{sec8-} and~\ref{sec8}.
The estimates for boundary and cutoff terms are in section~\ref{sec9}.

We now conclude this section with estimates of some integrals.
\le{outl} Let $n\geq2$.
Let $a$ be a real number, and $J=(j_1,\ldots, j_{2n-1})$
be a multiindex of non-negative
integers. Let $\beta= (j_{2n-1}+|J|-2a)+2n-1$ and $0<\rho_1\leq
\rho_0<+\infty$. Then
\gan
\int_{|z'|\leq\rho_1, |x^n|\leq\rho_0}
\f{|(z',\ov{z'},x^n)^J|}{||z'|^2+ix^n|^a}\, dV\leq \begin{cases}
C \rho_1^{1+\beta},& -1<\beta<2n-3,\vspace{.75ex}\\
C(1+|\log\rho_1|), &\beta=-1,\\
C\rho_1^{2n-2},&\beta\geq 2n-3,
\end{cases}\\
\int_{\rho_1\leq |z'|\leq \rho_0, |x^n|\leq\rho_0}
\f{|(z',\ov{z'},x^n)^J|}{||z'|^2+ix^n|^a}\, dV\leq \begin{cases}
C |\rho_1^{1+\beta}-\rho_0^{1+\beta}|,& \beta\neq-1,\beta<2n-3,\vspace{.75ex}\\
C \log(\rho_0/\rho_1), &\beta=-1,\\
C(\rho_0-\rho_1),&\beta\geq 2n-3,
\end{cases}
\end{gather*}
where $C$ depends only on $n, \rho_0$ and $|\beta|$.
\end{lemma}
\begin{proof} Set
$|J|_2=j_1+\cdots+j_{2n-2}+2j_{2n-1}$. When $|z'|^2\leq|x^n|$ and $|x^n|\leq\rho_0$,
we have $|x^n|\leq ||z'|^2+ix^n|\leq 2|x^n|$ and
 \aln
 b(z',x^n)=\f{|(z',\ov{z'},x^n)^J|}{||z'|^2+ix^n|^a}&\leq C
 |x^n|^{\yt|J|_2-a}\leq C' (|z'|^2+|x^n|)^{\yt|J|_2-a}.
 \end{align*}
 When $|z'|^2\geq|x^n|$ and $|z'|\leq\rho_0$, we have
 \aln
b(z',x^n)&\leq
 |z'|^{|J|_2-2a}\leq C (|z'|^2+|x^n|)^{\yt |J|_2-a}= \f{C}{(|z'|^2+|x^n|)^{n-\f{\beta+1}{2}}}.
 \end{align*}
Assume first that $n-\f{\beta+1}{2}>1$, i.e. $\beta<2n-3$.
 Using   polar coordinates, we get
 \aln
  \int_{|z'|\leq\rho_1, |x^n|\leq\rho_0}
b(z',x^n)\, dV  &\leq
C\int_{r=0}^{\rho_1}  dr\int_{x_n=0}^{\rho_0}
 \f{r^{2n-3}\, dx^n }{(r^2+|x^n|)^{n-\f{\beta+1}{2}}}\leq C'\int_0^{\rho_1}r^\beta\,dr.
 \end{align*}
Also, $\int_{\rho_1\leq |z'|\leq\rho_0, |x^n|\leq\rho_0}
b(z',x^n)\, dV  \leq C'\int_{\rho_1}^{\rho_0}r^\beta\,dr$.
 The estimates in this case follow.

 Assume now that $\beta\geq 2n-3$.
 The inequalities hold trivially if $j_{2n-1}\geq a$, in which case
 $b(z',x^n)\leq C$. When $j_{2n-1}<a$, we set $J=(J',j_{2n-1})$ and get $|J'|\geq 2(
 a-j_{2n-1})$.
 Hence
 $$
 b(z',x^n)\leq \f{|(z',\ov{z'})^{J'}|}{||z'|^2+ix^n|^{a-j_{2n-1}}}\leq C,
 $$
 from which the estimates follow.
 \end{proof}
\setcounter{thm}{0}\setcounter{equation}{0}
\section{ New kernels and  two formulae for derivatives
}\label{sec8-}

In this section we   write down the new kernels by using
\rp{m2h} and
  derive two formulae for the derivatives of $P_0\varphi$, where
 $\varphi$ has compact support in $D_\rho$.
The first formula
 will be used for $\cC^a$ estimates and the second is for $\cC^{k+\yt}$ estimates.

 Recall that with $\zeta, z\in M_\rho$  the coefficients of $P'_0\varphi(x)$ are sums of
$$
\mathcal I(x)=\int_{D_\rho}
\varphi_{\ov J}(\xi)
\f{\pd_*^2r(\xi,x)(r_{\zeta^j}-r_{z^j})}{(N_0^{a}S_0^{b})(\zeta,z)} \, dV(\xi)$$
over   $1\leq j\leq n$ and $|J|=q$, where $ a=n-q$ and $b=q$.
Set $
f(\xi,x)=\varphi_{\ov J}( \xi )$.

We have defined $\pd_*^{2+k}r$ in section~\ref{sec2}. Now, CHANGE NOTATION
and let
 \gan \pd_*^2
 r=p\left(\xi,x, (1+\hat r_{x^n}\hat r_{\xi^n})^{-1}, r_{\ov{z^n}}^{-1}, (
 1+\hat r_{x^n}\tr{(0',1)}\hat r(\xi,x))^{-1},
  \pd \hat r(x),  \pd  \hat r(\xi), \mathcal Q(\xi,x)\right),\\
 \mathcal Q(\xi,x)=\left(\pd^2\hat r(x), \pd^2\hat r(\xi),
 \trd_*^1 \pd\hat r(\xi,x),
 \trd_*^1 \hat  r(\xi,x), \trd_*^2 \hat  r(\xi,x)
  \right
 ),
 \end{gather*}
where     $p$  is a  polynomial.
Again,  $ \trd_*^1 \pd\hat r,
\trd_*^2 \hat  r$, defined by \re{cont} and \re{trpdr}, and $p$ might be  different  when they reoccur;
for instance, $(\trd_*^1 r)^2$
may be the product of two different $\trd_*^1 r$'s.
 Define
 \gan
 \pd_*^{2+k}  r(\xi,x)=\sum \pd_*^2  r\cdot
 \pd^{I_1} \mathcal Q_1(\xi,x)\cdots \pd^{I_j} \mathcal Q_j(\xi,x),\quad j\geq0,\\
 \pd_*^{2}  \hat r(\xi,x)=\sum \pd_*^2  r\cdot
  \mathcal Q_1(\xi,x)\cdots \mathcal Q_j(\xi,x), \quad j\geq 1,
 \end{gather*}
 where $ \sum_{l=1}^j |I_l|\leq k$,  $ \mathcal Q_l\in \mathcal Q$
  and both sums
 have finitely many terms.
 Hence, we have simple relations
 \gan
 \pd_*^{2+k}  r\pd_*^{2+j}  r=\pd_*^{2+k+j} r,\qquad
 \pd^J\pd_*^{2+k} r=\pd_*^{2+k+|J|}  r.
\end{gather*}
The chain rule     takes the form
 \ga
 \pd^I_{ \xi_*,x}\Psi^{-1}=\pd^2_*r\circ\Psi^{-1}, \quad|I|=1,\nonumber 
 \\ 
\label{pdjf} \pd^J_{ \xi_* ,x }\{(f\pd_*^{2+k}r)\circ\Psi^{-1}\}=\sum_{|L|\leq |J|}
 (\pd^{L}  f\cdot \pd_*^{2+k+|J|-|L|}r)\circ\Psi^{-1}.
  \end{gather}

\noindent{\bf New kernels}. Set $(\xis,x)=\Psi(\xi,x)$ with $\zeta,z\in M_{\rho_0}$.
Recall that $N_*(\xis)=|\xis'|^2+i\xis^n$. By \rl{m2h},
\ga
\nonumber
N_0(\zeta,z)\equiv 2\rzed=N_*(\xis) \hat T_1( \xi_* ,x ),
\\
\nonumber S_0(\zeta,z)\equiv -2\rzd=\ov {N_*(\xis)}
\hat T_2( \xi_* ,x ),\\
\label{t1=1}
\hat T_1(\xi_* ,x)=T_1\circ\Psi^{-1}(\xi_*,x)=
 1+ \sum_{|J| =2} \pd_*^2\hat r\circ\Psi^{-1}(\xi_* ,x)
N_*^{-1}( \xi_* )   \xi_* ^J,
\\
\label{t2=1}
\hat T_2(\xi_* ,x)=T_2\circ\Psi^{-1}(\xi_*,x)=
 1+\sum_{|J| =2} \pd_*^2\hat r\circ\Psi^{-1}(\xi_* ,x)
 \ov {N_*^{-1}( \xi_* )}\xi_* ^J,\\
\nonumber
\pd_*^2r(\zeta,z)(r_{\zeta^j}-r_{z^j})
=  \sum_{|I|=1} \pd_*^2 r\circ\Psi^{-1}( \xi_* ,x ) \xi_* ^I. 
 \end{gather}
 Note that $|\pd_*^2\hat r|\leq C\epsilon$ and $|\hat T_j(\xis,x)-1|<1/2$ when $\xis\neq0$.
Thus, we obtain
\ga\label{ty12-}
 \f{\pd_*^2r(\zeta,z)(r_{\zeta^j}-r_{z^j})}{(N_0^{a}S_0^{b})(\zeta,z)}=
\sum_{|I|=1}
\left\{\pd_*^2  r\circ\Psi^{-1}  \cdot\hat T_1^{-a}\hat T_2^{-b}\right\}( \xi_*,x)\hat k_{ab}^I(\xi_* ),
 \\  \quad \hat k_{ab}^I(\xis)={\xi_*^I}{
\label{ty12} N_* ^{-a}(\xis) \ov {N_*^{-b}(\xis)}}, \quad a=n-q, \  b=q.
\end{gather}

\noindent{\bf First formula of  derivatives of $\mathcal I$.} Recall that
$\Psi(\xi,x)=(\tilde\psi_x(\xi),x)$,
 $\xi_*=\tilde\psi_{x}(\xi)$ and  $(\tilde\psi_{x}^*dV)(\xi_*)=(\pd_*^1r)
 \circ\Psi^{-1}(\xi_*,x) dV(\xi_*)$.
By \re{ty12-} and $\pd_*^2r\pd_*^1r=\pd_*^2r$,
 we obtain
\eq{firstfa}
\mathcal I(x)=\sum_{|I|=1}\int_{\tilde\psi_x(D_\rho)}
\left\{(f
\pd_*^2  r)\circ\Psi^{-1}\cdot
\hat T_1^{-a}\hat T_2^{-b}\right\}(\xi_*,x)\hat k_{ab}^I( \xi _*) \, dV(\xi_*),
\end{equation}
where $a=n-q, b=q$. By  \rl{outl}, $\hat k_{a b}^I\in L_{loc}^1$.
For each $x\in D_\rho$, the integrand has   compact support in
 $\tilde\psi_x(D_\rho)\subset B_{9\rho}$. To compute $\pd^k \mathcal I(x)$,
 we  extend the integrand of $\mathcal I(x)$ to be zero on
 $B_{9\rho}\setminus \tilde\psi_x(D_\rho)$.
The integral is
 over the fixed domain $B_{9\rho}$.
 So we can interchange the integral sign with $\pd_{x}$, and
    derivatives of $\mathcal I$ have the form 
 \al\label{for1xi}
 &\pd^K\mathcal I(x)=\sum_{j+k'+l+m=|K|}\ \
 \sum_{|J|=j}\ \ \sum_{|K'|=k'}\ \ \sum_{|L|=l}\ \ \sum_{|I|=1}\\
 &\quad\int_{B_{9\rho}} \left\{(\pd^{L} f
 \pd_*^{2+m}r)\circ\Psi^{-1} \cdot  \pd_{x}^{J}{\hat T_1^{-a}}\cdot
 \pd_{x}^{K'}{\hat T_2^{-b}}\right\}(\xis,x)\cdot
\hat k_{ab}^I( \xi _*) \, dV ( \xi _*).
\nonumber\end{align}
By \re{t1=1},  the first-order   $x$-derivatives on $\hat T_1^{-a}$ have the form
\aln
\pd^I_x\hat T_1^{-a}(\xis,x)&=\sum_{|J'|=1}\ \
\sum_{|L'|=2}\{\hat T_1^{-a-1} \pd_x^{J'}(\pd^2_*r\circ\Psi^{-1})\}(\xis,x)N_* ^{-1}(\xis)\xis^{L'}
 \\ & \qquad
  =\sum_{|L'|=2}(\pd^{3}_*r\circ\Psi^{-1}\hat T_1^{-a-1})(\xis,x) N_* ^{-1}(\xis)\xis^{L'}.
  \end{align*}
Take derivatives consecutively and use the product rule. We can write
\al\label{t1a-s}
\pd_x^{J}\hat T_1^{-a}(\xis,x)&= \sum_{s\leq |J|}\
\sum_{|L'|=2s}\left\{(\pd_*^{2+|J|-s}r)\circ\Psi^{-1}
  \cdot\hat T_1^{-a-s}\right\}(\xis,x)\f{\xi_* ^{L'}}{N_*^{s}(\xis)}.
 \end{align}
Analogously,    $\pd_{x}^{K'}\hat T_2^{-b}$ ($|K'|=k'$) can be written as
\al\label{t2b-t}
\pd_x^{K'}\hat T_2^{-b}(\xis,x)=\sum_{t\leq |K'|}\sum_{|L''|=2t} \left\{(\pd_*^{2+|K'|-t}r)\circ\Psi^{-1}
\cdot \hat T_2^{-b-t}\right\}(\xis,x)\f{\xi_* ^{L''}}{\ov{N_*^{t}( \xi_* )}}.
\end{align}
Let $I'=I+L'+L'',  a'=a+s,  b'=b+t$. We have
$$ 
N_*^{-s}(\xis)\xi_* ^{L'}\ov{N_*^{-t}( \xi_* )}
\xi_* ^{L''}\hat k_{ab}^{I}=\hat k_{a'b'}^{I'},\quad 2a'+2b'-|I'|=2n-1.
$$
By   \re{t1a-s}-\re{t2b-t},
    derivatives of $ \mathcal I(x)$ have the form
\al \label{for1}
\pd^K\mathcal I(x)&=
 \sum_{|L|\leq |K|}\ \ \sum_{a\leq a'\leq a+|K|}\ \ \sum_{b\leq b'\leq b+|K|}\ \
 \sum_{2a'+2b'-|I'|=2n-1}\\
&\qquad\int_{B_{9\rho}} \f{(\pd^Lf
\pd^{2+|K|-|L|}_*r)\circ\Psi^{-1}( \xi_* ,x )}{(\hat T_1^{a'}\hat T_2^{b'})(\xi_*,x)}
\hat k_{a'b'}^{I'}
( \xi_* )
 \, dV(\xis).
 \nonumber
 \end{align}%
Recall that $f(\xi,x)=\varphi_{\ov J}(\xi)$ has compact support in $D_\rho$.
Since $|I'|=2a'+2a'-2n+1$, \rl{outl} implies
 $\hat k_{a'b'}^{I'}\in L^1_{loc}$. By the dominated convergence theorem,
we see that $\pd^k\mathcal I$ are continuous. Note that this also implies
that if $\hat r\in\cC^{k+2}(\ov D_{\rho_0})$ and $\varphi\in C^k(\ov
D_{\rho_0})$, then $P'\varphi\in
\cC^k(D_{\rho_0})$.

\

\noindent{\bf Second formula of  derivatives of $\mathcal I$}.
 We return to the original coordinates by letting
$( \xi_*,x)=\Psi(\xi,x)$ and  $\zeta, z\in M_\rho$. So $dV(\xi_*)=\pd_*^1r\, dV(\xi)$. Also
\gan
\hat T_1^{-a'}\circ\Psi(\xi,x)  ={N^{a'}(\zeta,z)} {N_0^{-a'}(\zeta,z)},\quad
\hat T_2^{-b'}\circ\Psi(\xi,x) ={\ov{N^{b'}(\zeta,z)}} {S_0^{-b'}(\zeta,z)},\\
\hat k_{a'  b'}^{I'}\circ\tilde\psi_x(\xi )
=\f{(\RE(\zeta'-z'),\IM(\zeta'-z'),|\zeta'-z'|^2+i2\IM(r_z\cdot(\zeta-z)))^{I'}}{
N(\zeta,z)^{a'}\ov {N^{b'} (\zeta,z)}}.
\end{gather*}
Multiply the same sides of the three identities and expand the last
numerator.  Then
$\{\hat T_1^{-a'}\hat T_2^{-b'}\hat k_{a'b'}^{I'}\}\circ\Psi$ is a  linear combination of
\ga
k_{a'b'}^{I''}(\xi,x)\df
\f{(\zeta'-z',\ov{\zeta' - z'},  \IM(r_z\cdot(\zeta-z)))^{I''}}
{N_0^{a'}(\zeta,z)S_0^{b'}(\zeta,z)}.
\nonumber
\end{gather}
Since $|I'|=2a'+2b'+1-2n$, then
\eq{ipij}
 |I''|-2a'-2b'+1-2n\geq0, \quad
  |I''|\leq 2(2|K|+1).
  \end{equation}
(In fact $|I''|-|I'|=i_{2n-1}'-i_{2n-1}''\geq0$.)   By  \re{for1},
    derivatives of $ \mathcal I(x)$ have the form
\al\label{for2}
\pd^K\mathcal I(x)&=\sum_{|L|\leq |K|}
\ \ \sum_{a\leq a'\leq a+|K|}\ \ \sum_{b\leq b'\leq b+|K|}\ \
\sum_{1+2a'+2b'-2n\leq |I''|\leq 4|K|+2}
 \\
&\qquad
\int_{D_{\rho}}
 \pd^L\varphi_{\ov J}( \xi )
\pd^{2+|K|-|L|}_*r(\xi,x)    k_{a'b'}^{I''}
(\xi,x)
 \, dV(\xi).
 \nonumber
\end{align}
Here $\varphi_{\ov J}$ has compact support
in $D_\rho$. Obviously,  the  $\pd^L\varphi_{\ov J}$ in \re{for2}
 do  not depend  on $x$. This simple observation
  will be crucial for the $\yt$-estimate.

The reader might want to acquaint with the counting
scheme in section~\ref{sec2}
and H\"older inequalities on domains $D_\rho$
 in Appendix~A;
 see \rp{vbhi}.

\setcounter{thm}{0}\setcounter{equation}{0}
\section{\hspace{1ex}$\cC^a$-estimates,
case of compact support}\label{sec8}
In this section, we derive the $\cC^a$-estimate for $P_0\varphi_0$
where
\eq{v0uj}
\varphi_0=\chi\varphi, \quad \|\chi\|_{\rho,a}\leq C_a(\rho\sigma)^{-a}
\end{equation}
and $\chi$  is supported in $D_\rho$. We also derive an estimate for $P_0\varphi$
when $\varphi$ itself has compact support in $D_\rho$.
\pr{cptv} Let $k\geq0$ be an integer and $0\leq\alpha<1$.
Let $M\colon y^n= |z'|^2+\hat r(z',x^n)$
satisfy \rea{bhroo} and $0<\rho\leq\rho_0\leq3$.
  Let $\varphi_0$ be a tangential form as  in \rea{v0uj}. Then
\ga
\label{fk0}
\| P_0'\varphi_0 \|_{ \rho,k+\alpha}\leq C_k \rho^{1-k-\alpha}\sigma^{-k-\alpha}
\bigl(\|\varphi\|_{\rho,k+\alpha}
 +\|\varphi\|_{\rho,0}\|\hat r\|_{\rho,k+2+\alpha}\bigr).
\end{gather}
If $\varphi$
  has  compact support in $D_\rho$ and is tangential, then
\ga\label{fk0+}
\| P_0'\varphi \|_{ \rho,k+\alpha}\leq C_k
\rho^{1-k-\alpha}
\bigl(\|\varphi\|_{\rho,k+\alpha}
 +\|\varphi\|_{\rho,0}\|\hat r\|_{\rho,k+2+\alpha}\bigr).
\end{gather}
  The same estimate holds
 for $Q_0'$.
\end{prop}
\pf
Recall that by applying \re{for1}  to $f=\chi\varphi_{\ov L}$,    $k$-th derivatives  of a coefficient of
$P_0'\varphi_0$ are the  sum of finitely many
\gan
  \mathcal I_k(x)=\int_{B_{9\rho}} \f{u
\circ\Psi^{-1}(\xis,x)}{\hat T_1^{a'}\hat T_2^{b'}(\xis,x)}
 \hat k_{a'b'}^{I'}
(\xis)
 \, dV(\xis).
 \end{gather*}
 Here $a=q$, $b=n-q$, $a\leq a'\leq a+k,b\leq b'\leq b+k$  and
 \gan
 u(\xi,x)\df\pd^I\chi(\xi)\pd^J\varphi_{\ov L}(\xi)
\pd^{2+l}_*r(\xi,x),\quad |I|=i,\ |J|=j,\quad i+j+l=k,
\\
 \hat k_{a'b'}^{I'}(\xis)=
 N^{-a'}( \xi_* )\ov{N^{-b'}( \xis )}\xis^{I'},
 \quad 2a'+2b'=|I'|+2n-1.
\end{gather*}
To obtain  \re{fk0},
we estimate the $\cC^\alpha$-norm of $  \mathcal I_k$.

By the definition of $\pd_*^{2+l}r$, we have
$$
\|\pd^{2+l}_*r\|_{\rho,\alpha}\leq C\sum_{l_1+\cdots+l_t\leq l}
\|r\|_{2+l_1+\alpha}\cdots\|r\|_{2+l_t}.
$$
 Thus
\aln
&\|\pd^i\chi\|_{\rho,0}\|\pd^j\varphi_{\ov L}\|_{\rho,0}\|\pd^{2+l}_*r\|_{\rho,\alpha}
\leq C(\rho\sigma)^{-i}\sum_{l_1+\cdots+l_t\leq l} \|\varphi_{\ov L}\|_{\rho,j}
\|r\|_{2+l_1+\alpha}\cdots\|r\|_{2+l_t}\\
&\hspace{10em}\leq C'(\rho\sigma)^{-i}\rho^{-l-j-\alpha}
(\|\varphi_{\ov L}\|_{\rho,i+j+\alpha} +\|r\|_{2+i+j+\alpha}\|\varphi_{\ov L}\|_{\rho,0}).
\end{align*}
Here the last inequality is obtained by  \rp{vbhi}.
Also
\aln
\|\pd^i\chi\|_{\rho,\alpha}\|\pd^j\varphi_{\ov L}\|_{\rho,0}\|\pd^{2+l}_*r\|_{\rho,0}
&\leq C(\rho\sigma)^{-i-\alpha}\rho^{-l-j}
(\|\varphi_{\ov L}\|_{\rho,i+j} +\|r\|_{2+i+j}\|\varphi_{\ov L}\|_{\rho,0}),\\
\|\pd^i\chi\|_{\rho,0}\|\pd^j\varphi_{\ov L}\|_{\rho,\alpha}\|\pd^{2+l}_*r\|_{\rho,0}
&\leq \f{C}{(\rho\sigma)^{i}\rho^{l+j+\alpha}}
(\|\varphi_{\ov L}\|_{\rho,i+j+\alpha} +\|r\|_{2+i+j+\alpha}\|\varphi_{\ov L}\|_{\rho,0}).
\end{align*}
Therefore,
\ga\label{asim}
\|u\|_{D_\rho^2,\alpha}
 \leq C_k  \rho^{-k-\alpha} \sigma^{-k-\alpha}
(\|\varphi\|_{\rho,k+\alpha}+
 \| r\|_{\rho,k+2+\alpha}\|\varphi\|_{\rho,0}),\\
 \label{asim+}
\|u\|_{D_\rho^2,0}\|r\|_{\rho,2+\alpha}
 \leq C_k  \rho^{-k-\alpha} \sigma^{-k }
(\|\varphi\|_{\rho,k+\alpha}+
 \| r\|_{\rho,k+2+\alpha}\|\varphi\|_{\rho,0}).
\end{gather}
  By \rl{psii}, we know that  $ W_\rho=\Psi(D_\rho\times D_\rho)\subset
 B_{9\rho}\times D_\rho$ and
\eq{pbdp}
 |\Psi^{-1}(v)-\Psi^{-1}(u)|\leq C|v-u|, \quad u,v\in W_\rho.
\end{equation}
Assume that $0\leq\alpha<1$. Fix $ \xis\in B_{9\rho}\setminus\{0\}$
and $x_1,x_2\in D_\rho$.
Assume first that
 $  \xi_j=\tilde\psi_{x_j}^{-1}(\xis)$ are in $D_\rho$
 for $j=1,2$. First, by \re{pbdp} $$|\xi_2 -\xi_1|\leq C
 |x_2
 -x_1| .$$
 Now by \re{t1=1}-\re{t2=1}, we obtain  $|\hat T_j(\xis,x)|\geq1/4$ and
 $$|\hat T_j(\xi_*,x_2)-\hat T_j(\xi_*,x_1)|\leq C|\pd_*^2r(\xi_2,x_2)-
 \pd_*^2r(\xi_1,x_1)|\leq\|r\|_{\rho,2+\alpha}|x_2-x_1|^\alpha.$$
  Thus
\aln
 \Delta&=|\Delta(x_2)-\Delta(x_1)|
 \df\left|\f{u\circ\Psi^{-1}(\xis,x_2)}{\hat T_1^{a'}\hat T_2^{b'}(\xis,x_2)}-
 \f{u\circ\Psi^{-1}(\xis,x_1)}{\hat T_1^{a'}\hat T_2^{b'}(\xis,x_1)}\right|
\\
&  \leq C\|u\|_{\rho,\alpha }   |x_2  -x_1|^\alpha
 +C|u (\xi_2,x_2)
((\pd_*^2  r)(\xi_2,x_2) -
(\pd_*^2  r)(\xi_1,x_1) )|  \\
&   \leq C  (
\|u \|_{\rho,\alpha }
+\|u \|_{\rho,0 }
\|\hat r\|_{\rho,2+\alpha})  |x_2  -x_1|^\alpha.
\end{align*}
By \re{asim}-\re{asim+} we get
\eq{eDel}
\Delta\leq C_k  \rho^{-k-\alpha} \sigma^{-k-\alpha}
(\|\varphi\|_{\rho,k+\alpha}+
 \| r\|_{\rho,k+2+\alpha}\|\varphi\|_{\rho,0})|x_2-x_1|^\alpha.
 \end{equation}
The above holds trivially if $\xi_1,
\xi_2$ are both not in $D_\rho$,
in which case $\Delta=0$.
If $\xi_2\in D_\rho$ and $\xi_1\not\in D_\rho$,
 we replace
$x_1$ by a point $x_3$
in the line segment $[x_1,x_2]$,
for which $\xi_3=\tilde\psi_{x_3}^{-1}(\xis)\in \pd
D_\rho$.  Then $\Delta=|\Delta(x_2)|=|\Delta(x_2)-\Delta(x_3)|$ and
\re{eDel} still holds.
By \rl{outl} (with $\rho_1=\rho_0=9\rho$ and $\beta\geq0$),
$
\int_{B_{9\rho}}|\hat k_{a',b'}^{I'}|dV \leq C \rho.
$
Combining the above estimates yields  \re{fk0}.

The case that $\varphi_0=\varphi$ is simpler, and it does
not involve $\sigma$. So we can remove all powers of $\sigma$
in \re{asim}-\re{asim+}, \re{eDel}, and   \re{fk0}.
The latter becomes \re{fk0+}.
\end{proof}

We compute the $\cC^\yt$ norm of
 $\pd^I\chi\pd^J\varphi_{\ov L}\pd^{2+l}_*r$ for a later use. Here it is crucial to avoid the
 H\"older $\yt$-norm on $\pd^J\varphi_{\ov L}$.
\pr{cptv+}  Let $u(\xi,x)=\pd^I\chi(\xi)\pd^J\varphi_{\ov L}(\xi)\pd^{2+l}_*r(\xi,x)
$ where $\chi$ has compact support in $D_\rho$
and satisfies $\|\chi\|_{\rho,a}\leq C(\rho\sigma)^{-a}$.  Let $|I|+|J|+l=k$. Then
\ga
\label{fk012}
 \|u(
\xi,\cdot)\|_{\rho,\yt} 
\leq C_k (\rho \sigma)^{-k-\yt}
\bigl(\|\varphi_{\ov L}\|_{\rho,k}\|  r\|_{\rho,\f{5}{2}}
 +\|\varphi_{\ov L}\|_{\rho,0}\|  \hat r\|_{\rho,k+\f{5}{2}}\bigr).
\end{gather}
\end{prop}
\begin{proof}
Fix $\xi\in D_\rho$.
The $\varphi_{\ov L}$ appearing in $u(\xi,x)$
 depends only on $ \xi $.
Therefore,  for
 $\|u(\xi,\cdot)\|_{\rho,\alpha}$
($\alpha=1/2$), we only use the sup norm of $\pd^j\varphi_{\ov L}(\xi)$.
Then \aln
\|u(\xi,\cdot)\|_{\rho,1/2}
&\leq C_k( (\rho\sigma)^{-i-\yt}\|\varphi_{\ov L}\|_{\rho,j}\|r\|_{\rho,
l+2}+ (\rho\sigma)^{-i}\|\varphi_{\ov L}\|_{\rho,j}\|r\|_{\rho,
l+\f{5}{2}}).
\\
&\leq C_k'(\rho\sigma)^{-i-\yt} \rho^{-j-l}(\|\varphi_{\ov L}\|_{\rho,j+l}+\|\varphi_{\ov L}\|_{\rho,
0}\|r\|_{\rho,2+j+l})
\\
&\quad +C_k'(\rho\sigma)^{-i}\rho^{-j-l-\f{1}{2}}
(\|\varphi_{\ov L}\|_{\rho,0}\|r\|_{\rho, j+l+\f{5}{2}}+\|\varphi_{\ov L}\|_{\rho,j+l}\|r\|_{\rho, \f{5}{2}}),
\end{align*}
where the last two terms
are
obtained by \rp{vbhi} in which we take  $d_1=0$ and $d_2=\f{5}{2}$. Simplifying  yields
\re{fk012}.\end{proof}

\setcounter{thm}{0}\setcounter{equation}{0}\section{Boundary integrals,
end of $\cC^a$-estimates}
\label{sec9}
In this section we will estimate the boundary integrals $P'_1\varphi$  and
cutoff term $P'_0\varphi_1$, where $\varphi_1$ vanishes on $D_{(1-\yt\sigma)\rho}$.
  Estimates \re{p1va} and \re{p0v0}
below
will be used again
for the $\cC^{k+\yt}$  estimate.

   Recall \re{hftk}-\re{hftk+} that for $\zeta, z\in M_\rho$
   \gan
  P'_0\varphi(x)=\sum_{|I|=q-1}\ \sum_{|J|=q}\ \sum_{1\leq j\leq n} d\ov{z'}^I
 \int_{D_\rho}
A_{\ov I}^{j\ov J}(\xi,x)
\f{\varphi_{\ov J}(\xi) (r_{\zeta^j}-r_{z^j})}{(N_0^{n-q}S_0^{q})(\zeta,z)} \, dV(\xi),\\
P'_1\varphi(x)= \sum_{|I|=q-1}\, \sum_{|J|=q}\, \sum_{ \alpha,\beta=1}^{n-1 }\, \sum_{s=1}^{2n-1}
d\ov{z'}^I
 \int_{\partial D_\rho}
\f{B^{\alpha\beta \ov J}_{\ov Is}(\xi,x) \varphi_{\ov J}(\xi)(r_{\zeta^\alpha}-
r_{z^\alpha})r_{\zeta^\beta}}{(\zeta^n-z^n)(N_0^{n-q-1}S_0^{q})(\zeta,z)} \,
d V^s(\xi).
\end{gather*}
Here $  A_{\ov I}^{j\ov J} $ and $ B^{\alpha\beta \ov J}_{\ov Is} $ are polynomials in
$  (r_{\zeta },   r_{ \ov\zeta},
 r_{\zeta\ov\zeta},   r_{\ov{z'}}, r_{\ov{z^n}}^{-1},   r_{z\ov z})$.
And
 $N_0(\zeta,z)=r_\zeta\cdot(\zeta-z), S_0(\zeta,z)=r_z\cdot(\zeta-z)$.
  For $\zeta, z\in M_\rho$, set
\gan
k (\xi,x)\df
k_{\ov I
}^{j\ov J} (\xi,x)=A_{\ov I}^{j\ov J}(\xi,x)
\f{ (r_{\zeta^j}-r_{z^j})}{(N_0^{a_0}S_0^{b_0})(\zeta,z)},
\quad a_0+b_0=n, \quad b_0=q,\\
l (\xi,x)\df
l_{ \ov Is}^{\alpha\beta \ov J }(\xi,x)=B^{\alpha\beta \ov J}_{\ov Is}(\xi,x)\f{
r_{\zeta^\beta}(r_{\zeta^\alpha}-r_{z^\alpha}) }
{(\zeta^n-z^n)(N_0^{a_0-1}S_0^{b_0})(\zeta,z)}.
\end{gather*}

Recall that    $\chi=\chi_{\sigma,\rho}$ is a smooth function,
  which is $1$  on  $D_{\rho(1-\sigma/2)}$
and zero  off $D_{\rho (1-\sigma/4)}$, and $\|\chi\|_{\rho,a}\leq
C_a(\rho\sigma)^{-a}$.
On $D_{\rho}$, set
\gan
\varphi_0=\chi \varphi, \quad
\varphi_1=\varphi-\varphi_0.\end{gather*}

Assume that $M$ satisfies \re{bhroo-}-\re{bhroo}. Thus
 $\|\hat r\|_{\rho_0,2}<1/C_0$. Assume that $0<\rho\leq\rho_0\leq3$.
  Set $R=|z'|^2+\hat r(x)$. Let $\pi_s$
  be the projection from $\pd D_\rho$ into
  the  subspace $\xi^s=0$. Since $\ov D_\rho$ is bounded and strictly convex,
  then $\pi_s$ is a 2-to-1 map
from   $\pd D_\rho$ onto $\pi_s(\pd D_\rho)=\pi(\ov D_\rho)$. Actually,
 $\pi_s$ sends   $\pi_s^{-1}(\pd(\pi_s(D_\rho)))$ one-to-one
and onto $\pd(\pi_s(D_\rho))$.
Let $ vol(\pi_s D_\rho)$ be the volume of $\pi_s(D_\rho)$ calculated
via the volume-form $dV^s$.
 Recall that $D_\rho$ is contained in $B_{2\rho}$.
 If $f$ is a continuous function on $D_\rho$, then
\al\label{ares}
\Bigl|\int_{\partial D_\rho} f(\xi) \,
dV^s(\xi)
\Bigr|&\leq 2\, vol(\pi_s D_\rho) \|f\|_{\rho,0}\leq 2\,
vol(B_{2\rho}\cap\rr^{2n-2})
\|f\|_{\rho,0}\\
&\leq C\rho^{2n-2}\|f\|_{\rho,0}.\nonumber
\end{align}
 Since   the projection of $\pd D_\rho$ in any
 coordinate hyperplane   is contained in a ball of radius $2\rho$,
by the Fubini theorem one  can verify that
\eq{arvo}
vol\, (D_\rho\setminus D_{(1-\yt\sigma)\rho})
\leq   (2n-1)\cdot C\rho\sigma\cdot
\rho^{2n-2}\leq C'\rho^{2n-1}\sigma.
\end{equation}

\

Let $k\geq0$ be an integer  and $0\leq\alpha<1$.
Fix $\zeta\in M_\rho\setminus
M_{(1-\yt\sigma)\rho}$ and vary $
z\in M_{(1-\sigma)\rho}$.
By  \rl{distb}, we have
 \ga\label{nsrr--}
 |r_\zeta\cdot(\zeta-z)|\geq C^{-1}(\rho\sigma)^2, \quad |r_z\cdot
 (\zeta-z)|\geq C^{-1}(\rho\sigma)^2.
\end{gather}
 We also have
   $|r_{\zeta^j}-r_{z^j}|\leq C\|\zeta-z\|\leq
 C'|r_\zeta\cdot(\zeta-z)|^{1/2}.$ Hence
 \eq{nsrr-}
 |r_{\zeta^j}-r_{z^j}||N_0|^{-a}|S_0|^{-b}
 \leq C|N_0|^{\yt-a-b}\leq C'(\rho\sigma)^{1-2a-2b}
 \end{equation}
  if $a+b\geq1/2$.
Using $|f(x_2)-f(x_1)|\leq \|f\|_{\rho,1}|x_2-x_1|\leq
\|f\|_{\rho,1}(|x_2|+|x_1|)^{1-\alpha}|x_2-x_1|^\alpha$, we get
 $$\|f\|_{\rho,\alpha}\leq \|f\|_{\rho,0}+
 C\|f\|_{\rho,1}\rho^{1-\alpha}.$$
 Therefore,  for $1\leq j\leq n$ and $0<\alpha<1$,
 \eq{rzbe}
 \|r_{\zeta^j}-r_{z^j}(\cdot)\|_{\rho,\alpha}\leq C\rho^{1-\alpha}, \quad
 |r_{\zeta^\beta}|\leq C\rho.
 \end{equation}
 Also for $x_1,x_2\in D_\rho$,
 \aln
 \Bigl|\f{1}{f(x_2)}-\f{1}{f(x_1)}\Bigr|&=\f{ |f(x_1)-f(x_2)| ^{1-\alpha}}
 {|f(x_2)f(x_1)|}|f(x_1)-f( x_2)|^\alpha\\
 \nonumber
& \leq 2^{1-\alpha}\|1/f\|_{\rho,0}^{1+\alpha}\|f\|_{\rho,1}^\alpha|x_2-x_1|^\alpha.
 \end{align*}
 Combining it with H\"older ratio $|1/{f^a}|_{\rho,\alpha}
 \leq C_a\|1/f\|_{\rho,0}^{a-1}|1/f|_{\rho,\alpha}$ for $a\geq1$, we get
 \eq{1fx2+}
\|1/{f^a}\|_{\rho,\alpha}
 \leq \|1/{f^a}\|_{\rho,0}+
 C_a\|1/f\|_{\rho,0}^{a+\alpha}\|f\|_{\rho,1}^\alpha, \quad a\geq1.
  \end{equation}
 Now, by \re{nsrr--}
 \ga
 \|N_0(\zeta,\cdot)^{-a}\|_{(1-\sigma)\rho,\alpha}+
 \|S_0(\zeta,\cdot)^{-a}\|_{(1-\sigma)\rho,\alpha}\leq C(\rho\sigma)^{-2(a+\alpha)},\quad
 a\geq1.\label{nsrr}
 \end{gather}
 Note that $A_{I}^{jJ}$ has the form $\pd_*^2r$. By \rp{vbhi}, we have
 \eq{nsrr0}
 \|\pd_*^2r(\xi,\cdot)\|_{\rho,a}\cdot\|\pd_*^2r(\xi,\cdot)\|_{\rho,b}\leq C_{a,b}\rho^{-a-b}
 \|r\|_{\rho, 2+a+b}.
 \end{equation}
We have
 \gan
\pd_x^Kk (\xi,x)=\sum_{a+ b+ c+|L|=|K|}\ \ \sum_{|L|=0,1}
\f{ \pd_*^{2+c}r(\xi,x)\pd_x^{L}
(r_{\zeta^j}-r_{z^j}) }
{(N_0^{a_0+a}S_0^{b_0+b})(\zeta,z)},
\end{gather*}
where $ a_0+b_0=n$. Fix $\xi$ and vary $x$.
For  the summand, we estimate the $\cC^\alpha$-norms of two terms
in the numerator by \re{nsrr--}-\re{rzbe},
  and use   \re{nsrr} for the reciprocals of two terms in the denominator. Set
 $ |K|=k$ and $|L|=d$. Recall that $d=0$ or $1$.  We get
\aln
&\Bigl\|\f{ \pd_*^{2+c}r(\xi,\cdot)\pd_{x}^{L}
(r_{\zeta^j}-r_{z^j}(\cdot)) }
{(N_0^{a_0+a}S_0^{b_0+b})(\zeta,\cdot)}\Bigr\|_{(1-\sigma)\rho,\alpha}
\leq \f{C}{(\rho\sigma)^{2(a+a_0+b+b_0)}}\Bigl\{
(\rho\sigma)^{1-d}
\rho^{-c-\alpha}\|r\|_{\rho,2+c+\alpha}\\
&\hspace{6.5em}
+\rho^{(1-d)(1-\alpha)}\rho^{-c-d\alpha}\|r\|_{\rho,2+c+d\alpha}+
2(\rho\sigma)^{1-d }
(\rho\sigma)^{-2\alpha} \rho^{-c}\|r\|_{\rho,2+c}\Bigr\}.
\end{align*}
The worst term  in terms of  powers of
$\rho,\sigma$
occurs when $a+b=k$, $c=d=0$.
We see that
  \al\label{ikabj}
\|k (\xi,\cdot)\|_{(1-\sigma)\rho, k+\alpha}\leq C_{k}
(\rho\sigma)^{-s_1} 
 \|r\|_{ {\rho},k+2+\alpha}, \quad \xi\in D_\rho\setminus D_{(1-\yt\sigma)\rho}
 \end{align}
with $ s_1 =2n-1+2k+2\alpha$.

To estimate  the   boundary term,  by
\re{nsrr++-} we have
\eq{nsrr++--}
|\zeta^n-z^n|\geq C^{-1}\rho^{2}\sigma
\end{equation}
for $z\in M_{(1-\sigma)\rho}$ and $\zeta\in\partial M_\rho$. Using \re{1fx2+},
  we get   for $\zeta\in\partial M_\rho$
\eq{nsrr++}
\|(\zeta^n-z^n(\cdot))^{-a}\|_{(1-\sigma)\rho,\alpha}\leq
C(\rho^2\sigma)^{-a-\alpha},\quad a\geq1, \quad 0\leq\alpha\leq1.
\end{equation}
Note that $B^{\alpha\beta J}_{Is}$ has the form $\pd_*^2r$.   We have  \aln
\pd_x^Kl (\xi,x)&=\sum_{a+b+d+c=|K|}\ \ \sum_{c_1+|L|=c}\ \ \sum_{|L|=0,1}
\f{\pd_*^{2+c_1} r  \cdot r_{\zeta^\beta}\cdot\pd_x^{L}
(r_{\zeta^j}-r_{z^j}) }
{(\zeta^n-z^n)^{1+d}(N_0^{a_0-1+a}S_0^{b_0+b})(\zeta,z)}.
\end{align*}
Set $|K|=k$ and $|L|=c_2$. We  now estimate the $\cC^\alpha$-norm in the $x$ variables.
 Fix $\zeta\in \pd D_{\rho}$.  Using \re{nsrr--}-\re{rzbe}, \re{nsrr0}
for three terms in the  numerator and  \re{nsrr},
 \re{nsrr++--}-\re{nsrr++} for the reciprocals of  three terms in
 the denominator,
we obtain
\al \label{ytes}
&\|l (\xi,\cdot)\|_{{(1-\sigma)\rho}, k+\alpha}\leq
 C_{k}\sum_{a+b+d+c_1+c_2=k}\ \sum_{ c_2=0,1}\
(\rho\sigma)^{-2(a_0+a+b_0+b+d)}\sigma^{1+d}\rho\cdot
\\
\nonumber
&\hspace{1.5em} \cdot \Bigl\{
(\rho\sigma)^{1-c_2}\rho^{-c_1-\alpha}\|r\|_{ {\rho},2+c_1+\alpha}
 +\rho^{(1-c_2)(1-\alpha)} \rho^{-c_1-c_2\alpha}
 \|r\|_{\rho,2+c_1+c_2\alpha}\\
 \nonumber &\hspace{9.5em}+ \rho^{-c_1}(\rho\sigma)^{1-c_2 }
 ((\rho^{2}\sigma)^{-\alpha}+2(\rho\sigma)^{-2\alpha})
 \|r\|_{\rho,2+c_1}
 \Bigr\}.
 \end{align}
 The worst term in terms of   powers
of $\rho,\sigma$ occurs when $c_1=c_2=d=0, a+b=k$.
This shows that for each $\zeta\in\ov M_\rho
 \setminus M_{(1-\yt\sigma)\rho}$, we have
\al\label{kabj}
\|l (\xi,\cdot)\|_{{(1-\sigma)\rho}, k+\alpha}\leq C_{k}
(\rho\sigma)^{-s_2} 
 \|r\|_{ {\rho},k+2+\alpha}
 \end{align}
 with  $s_2=2(n+k-1+\alpha)$.

 By \re{ares} and \re{kabj},
  we estimate the boundary term by
\eq{p1va}
\|P_1'\varphi\|_{ {(1-\sigma)\rho},k+\alpha}
\leq C_k\rho^{2n-2}(\rho\sigma)^{-s_2}
\|r\|_{ {\rho},k+\alpha+2}\|\varphi\|_{ \rho,0}.
\end{equation}
 Estimating the  cutoff term by \re{ikabj} and \re{arvo},
 we obtain
 \eq{p0v0}
  \|P_0'\varphi_1\|_{ {(1-\sigma)\rho},k+\alpha}
\leq C_a\rho^{2n-1-s_1}\sigma^{1-s_1}
\|r\|_{ {\rho},k+2+\alpha}\|\varphi\|_{ \rho,0}.
\end{equation}
Define $s\df\max\{s_1-1, s_2,k+\alpha\}$ and
$s_*\df\max\{s_1-2n+1,s_2-2n+2, k+\alpha-1\}$.
Thus  for $a=k+\alpha$, we get
\eq{sval}
  s=2(a+n-1),\quad  s_*= 2a.
 \end{equation}
Combining \re{fk0}, \re{p1va}-\re{p0v0},  we get the following.
\pr{dkck} Let $n\geq4$,
and let $a\geq0$ be   a real number. Let $M\colon
y^n= |z'|^2+\hat r(z',x^n)$ satisfy \rea{bhroo}.
Let $P'$ be either of $P',Q'$ in the \hf\
  $\varphi=\dbm P'\varphi+Q'\dbm\varphi$
  on $M_\rho$. Assume that $0<\rho\leq\rho_0\leq3$. Then for a tangential form $\varphi$
\gan
 \|P'\varphi\|_{D_{(1-\sigma)\rho},a}
\leq
C_a\rho^{-s_*}\sigma^{-s}(\|\varphi\|_{D_\rho,a} + \|\varphi\|_{D_\rho,0}
 \|\hat r\|_{D_\rho, a+2}),
\end{gather*}
where $0<\sigma<1$, $s, s_*$ are given by \rea{sval} .
\end{prop}

\setcounter{thm}{0}\setcounter{equation}{0}\section{Reduction of
 $\cC^{k+\yt}$-estimates to $\cC^\yt$-estimate. Summary}
\label{sec10}

We want to prove
the following $\cC^{k+\yt}$ estimates.
\pr{dkck+} Let $k\geq0$ be an integer. Let $M\colon
y^n= |z'|^2+\hat r(z',x^n)$ satisfy \rea{bhroo}.
Let $P'$ be one of $P',Q'$ in the \hf\
  $\varphi=\dbm P'\varphi+Q'\dbm\varphi$  on $M_\rho$.
Then for $0<\rho\leq\rho_0\leq3$, $0<\sigma<1$, and a tangential form $\varphi$
\gan
 \|P'\varphi\|_{D_{(1-\sigma)\rho},k+\yt}\leq
\f{C_k}{\rho^{2k+1}\sigma^{2n+2k-1}}\bigl
( \|  r\|_{D_\rho,\f{5}{2}}\|\varphi\|_{D_\rho,k}+
\|\hat r\|_{D_\rho, k+\f{5}{2}}\|\varphi\|_{D_\rho,0}\bigr),\\
\|P\varphi\|_{D_{(1-\sigma)\rho},\yt}\leq
C  \rho^{-1}\sigma^{1-2n} \|\varphi\|_{D_\rho,0},\quad q=1.
\end{gather*}
\end{prop}
\begin{proof}  Fix
   a positive integer $k$. To estimate the $\cC^{k+\yt}$-norm of
$P'\varphi=(P'_0+P'_1)\varphi$, we first recall estimates  \re{p1va} and \re{p0v0} for the boundary and
cutoff terms
\eq{cuto}
\|(P'_0\varphi_1,  P'_1\varphi)
\|_{ {(1-\sigma)\rho},k+\yt}
\leq C_a\rho^{-s_*}\sigma^{-s}
\|r\|_{ {\rho},k+\f{5}{2}}\|\varphi\|_{ \rho,0}.
\end{equation}
Here  $s=2n+2k-1$ and $s_*=2k+1$ are  computed by \re{sval}
for $a=k+\yt$.
It remains to
estimate the $\cC^{k+\yt}$-norm of  $P'_0\varphi_1$, where
  $\varphi_1=\chi\varphi$ and $\chi$
has  compact support in $D_{\rho}$ and  $\|\chi\|_{\rho,a}\leq C_a(\rho\sigma)^{-a}$.
The proof   will be completed later.
For the rest of proof, 
we reduce it
to the special case of $k=0$.

\

  The second formula \re{for2}   says that  the coefficients
  of $\pd^kP'_0\varphi_0$ are sums of
 \ga
\mathcal Ku(x)=\int_{D_\rho}
 u (\xi,x)
k  (\xi,x) \,dV \label{kIJ-} \end{gather}
with functions $u(\xi,x)$ and kernels $k(\xi,x)$ of the form
\ga
\label {speu}
u(\xi,x)=
\pd^E\chi (\xi)\pd^F\varphi_{\ov J}(\xi)\pd_*^{2+l}
r(\xi,x), \quad |E|+|F|+l=k,\\
\label{kIJ}
k (\xi,x) \df k_{ab}^I(\zeta,z)=\f{(\zeta'-z',\ov{\zeta'-z'},
\IM (r_z\cdot(\zeta-z)))^{I}}{(r_{\zeta}\cdot(\zeta-z))^{a}(r_z\cdot(\zeta-z))^{b}},
\quad  \zeta, z\in M_\rho.\end{gather}
Moreover,  by \re{ipij}, the non-negative integers $a,b,$
  $I=(i_1,\ldots,i_{2n-1})$  satisfy
\eq{kIJ+}
|I|-2a-2b\geq 1-2n, \quad |I|\leq 4k+2,\quad a+b\leq n+k.
\end{equation}
By \re{fk012} in \rp{cptv+},
we have \al\label{u12}
\|u(\xi,\cdot)\|_{
 \rho,\yt}
 \leq C_k(\rho\sigma)^{-k-\yt}  \bigl(\|\varphi\|_{ \rho, k} \|  r\|_{ \rho,\f{5}{2}}
 +\|\varphi\|_{ \rho,0}\|\hat r\|_{ \rho,k+\f{5}{2}}
\bigr).
\end{align}
In section~\ref{sec14}
 we will prove that if  $u\in L^\infty(D_\rho\times D_\rho)$   then
\ga\label{yte}
\|\mathcal Ku\|_{D_{\rho},{1/2}}
 \leq C
 \sup_{\xi\in D_{\rho}}\|u(\xi,\cdot)\|_{D_{\rho},{1/2}}.
\end{gather}
Combining it with \re{cuto} and \re{u12} yields
the first estimate in the proposition.

\

We now consider the case that  $\varphi$ is a tangential   $(0,1)$-form.
We return to \re{hftk}-\re{hftk+}
and  look at a special  property of the kernels of
$P'\varphi $.   Recall that in this case 
\gan
  P'_0\varphi(x)= \sum_{1\leq \gamma <n }\ \ \sum_{1\leq j\leq n}
\ \int_{D_\rho}
A^{j\ov \gamma}(\xi,x)
\f{\varphi_{\ov \gamma}(\xi) (r_{\zeta^j}-r_{z^j})}{(N_0^{n-1}S_0 )(\zeta,z)} \, dV(\xi),\\
P'_1\varphi(x)=   \sum_{1\leq\alpha,\beta,\gamma< n}
 \ \ \sum_{1\leq s<  2n}
 \ \int_{\partial D_\rho}
\f{B^{\alpha\beta \ov \gamma}_{s}(\xi,x) \varphi_{\ov \gamma}(\xi)(r_{\zeta^\alpha}-
r_{z^\alpha})r_{\zeta^\beta}}{(\zeta^n-z^n)(N_0^{n-2}S_0 )(\zeta,z)} \,
dV^s( \xi).
\end{gather*}
Here $\zeta, z\in M_\rho$. 
Also, 
$ A^{j\gamma}$ and $  B^{\alpha\beta \gamma}_{s}$ are polynomials in
$  r_{\zeta },   r_{ \ov\zeta},
 r_{\zeta\ov\zeta}$. In particular, they     are independent of $z$.
Moreover, in the kernels there are only  the first-order derivatives
$r_{z^j}$ in the $z$ variable. Then all norms of $r$
in \re{ytes}-\re{p1va}, in which $k=0$,
 can be replaced by $\|r\|_{\rho,2}<C$ and
the estimate \re{p1va} for the boundary term becomes
\gan
\|P'_1\varphi\|_{{(1-\sigma)\rho},\alpha}
\leq C\rho^{-2\alpha}\sigma^{2-2n-2\alpha} \|\varphi \|_{ \rho,0},\quad 0\leq\alpha\leq1.
\end{gather*}
Absorb  $A_{I}^{j\gamma}(\xi)$ into $\varphi_\gamma(\xi)$,
With  $\zeta,z\in M_\rho$ the kernels of interior integral $P'_0\varphi $ have the form
\eq{ker2}
k (\xi,x)=\f{r_{\zeta_j}-r_{z^j}}
{(r_{\zeta}\cdot(\zeta-z))^{a}(r_z\cdot(\zeta-z))^{b}},
\quad a+b=n.
\end{equation}
We will show that \re{yte} holds
for this new kernel. Applying it  to
$$u(\xi,x)=A^{j\ov \gamma}(\xi)\varphi_{\ov \gamma}(\xi),\quad \zeta\in M$$ we obtain
$
\|P'_0\varphi_{0,1}\|_{\rho,1/2}\leq C\|\varphi_{0,1}\|_{\rho,0}.
$
This shows the second estimate of the proposition.

The proof of \rp{dkck+} is thus complete, by assuming  \re{yte} in which  $\mathcal Ku$,
 given by \re{kIJ-},  has  a kernel of the form
 \re{kIJ}-\re{kIJ+} or \re{ker2}.
\end{proof}

We want to unify  the two kernels  \re{kIJ}, which  satisfies
\re{kIJ+}, and   \re{ker2}. For \re{kIJ}, we write
$$\IM (r_z\cdot(\zeta-z))=\f{1}{2i}
r_z\cdot(\zeta-z)-\f{1}{2i}
\ov{r_z\cdot(\zeta-z)}.$$  Then  kernel \re{kIJ} is a linear combination of
\eq{nker}
k(\xi,x)\df 
\f{(\zeta-z,\ov{\zeta-z},r_\zeta-r_z)^{I}}
{(r_\zeta\cdot(\zeta-z))^{a}(r_z\cdot(\zeta-z))^{b}
\ov{(r_\zeta\cdot(\zeta-z))^{c}}\, \ov{(r_z\cdot(\zeta-z))^{d}}},
\end{equation}
where $\zeta,z\in M$,
 $a, b, c, d$ are now  possibly negative  integers, $I$ is a $(3n)$-tuple
of nonnegative integers,
and
\eq{nker+}
|I|-2(a+b+c+d)\geq 1-2n,\  |I|\leq 4k+2,\  |a|+|b|+|c|+|d|\leq n+5k+2.
\end{equation}
Indeed,  the $I=(i_1,\ldots, i_{2n-1})$ in \re{kIJ+} is a multi-index
 of nonnegative integers.
 Hence, $\re{kIJ+}$ and $i_{2n-1}\geq0$ implies \re{nker+}.
  Obviously,
\re{ker2} is of the form \re{nker}. In
section~\ref{sec14}, for the new kernel \re{nker} with the condition \re{nker+}  we will prove \re{yte}, i.e.
\ga\label{yte+}
\Bigl\|\mathcal \int_{D_\rho}u(\xi,\cdot)k(\xi,\cdot)\, dV(\xi)\Bigr\|_{D_{\rho},{1/2}}
 \leq C
 \sup_{\xi\in D_{\rho}}\|u(\xi,\cdot)\|_{D_{\rho},{1/2}},
 \end{gather}
where $u\in L^\infty(D_\rho\times D_\rho)$.

\

\noindent{\bf Summary.}
We would like to    summarize our   observations on $\cC^{k+\yt}$-estimates
to explain how the general estimate \re{yte+} gives us  an  actual $\yt$-gain
in special cases.

\medskip

a) $k\ge1$ and $q\geq1$.
 In our applications of \re{yte+} to the  $\cC^{k+\yt}$-estimates,
 $u(\xi,x)$, arising from the $k$-th order derivatives
 of $P'_0\varphi$ after applying the approximate Heisenberg
 transformation,
 is of the form $\pd^{k-j}\varphi_{\ov J}(\xi)\pd^{j+2}_*r(\xi,x)$
 where $\varphi$ has compact support in $D_\rho$;
 see
 the second formula of the derivatives in section~\ref{sec8-}.
  In particular, the $\cC^{\yt}$-norm
 of $u(\xi,x)$ in $x$ variables involves only $\|\varphi\|_{\rho,k}$ (and
 $\|r\|_{\rho,k+\f{5}{2}}$). Thus the estimate for
 $\|P_0'\varphi\|_{\rho,k+\yt}$ gains $\yt$ in H\"older exponent  from
  $\|\varphi\|_{\rho,k}$ at the expense of   two extra derivatives
  in $\|r\|_{\rho,k+\f{5}{2}}$.


\medskip

 b) $k=0$ and $q=1$.  Apply \re{yte+} to $(0,1)$ forms. In this case,    we will not need to
 differentiate  the   kernels and apply cutoff. Both  simplify the estimate
 considerably.
 However, it is crucial that when \re{yte+} is applied, $u(\xi,x)$
 has the form
  $A^{j\gamma}(\xi)\varphi_\gamma(\xi)$; in particular, it is {\it independent}
  of $x$.
Therefore, the $ \yt$-estimate
 of $u$ in the $x$-variables and hence $ \yt$-estimate of $P'_0\varphi$ only
 require   $r\in\cC^2$.

\medskip

  c) While the $ \yt$-estimate of $P'_0\varphi$
  is on the original   domain (\rp{k12+} below), the
   estimate  for the boundary term $P'_1\varphi$, which only requires  $r\in \cC^2$ for
   the same reason as in b),
  is  on   shrinking  domains.

\setcounter{thm}{0}\setcounter{equation}{0}\section{$\yt$-estimate}
\label{sec14}

We will derive the
 $\yt$-estimate, by modifying a standard argument for H\"older estimates.
This will complete the proof of \rp{dkck+} via \re{u12} and \re{yte+}.
We restate the latter as the following.

\pr{k12+}
 Let $M$  be a real hypersurface of class $\cC^2$ and
 satisfy \rea{bhroo}.  Assume that $0<\rho\leq\rho_0\leq3$.
 Let $a,b,c,d$ be (possibly negative)  integers
and  $I$ be a $(3n)$-tuple
of nonnegative integers. Suppose that
 $|I|- 2(a+b+c+d)\geq1-2n$. Let
\gan
k(\xi,x)=\f{(\zeta-z,\ov{\zeta-z},r_\zeta-r_z)^{I}}
{(r_\zeta\cdot(\zeta-z))^{a}(r_z\cdot(\zeta-z))^{b}
\ov{(r_\zeta\cdot(\zeta-z))^{c}}\, \ov{(r_z\cdot(\zeta-z))^{d}}},\quad\zeta,z\in M,\\
\mathcal Ku(x)=\int_{M_\rho} u(\xi,x)k(\xi,x)\,dV.
\end{gather*}
 Then for $u\in L^\infty (D_\rho\times D_\rho)$,
  $$
\|\mathcal Ku\|_{D_\rho,\yt}\leq C
\Bigl( \rho\sup_{\xi\in D_\rho }
\|u(\xi,\cdot)\|_{D_\rho,\yt}+\sup_{\xi\in D_\rho }
\|u(\xi,\cdot)\|_{D_\rho,0}\Bigr).
$$
\end{prop}
\begin{proof}
Recall that $D_\rho$ is convex. It will be convenient
to regard $Ku$ as a function on $M$. The $\yt$-H\"older ratios
on $M_\rho$ and $\pi M_\rho$ are equivalence, since
$|\pi(z_2-z_1)|\geq |z_2-z_1|/C$ for $z_1,z_2\in M_\rho$.
By abuse of notation, we write $u(\xi,x)$
as $u(\zeta,z)$ with $\zeta,z\in M$ and apply the same change
of notation to  $k(\xi,x)$, $\mathcal
 Ku(x)$.

  Fix  $z_1,z_2$ in $M_\rho$.  We have
\aln
\mathcal Ku(z_2)-\mathcal Ku(z_1)&=\int_{M_\rho}  u(\zeta, z_2)
\bigl(k(\zeta, z_2)-k(\zeta, z_1)\bigr)\,dV\\
&\quad +
\int_{M_\rho} \bigl(u(\zeta, z_2)-u(\zeta, z_1)\bigr)k(\zeta, z_1)\,
dV=\mathcal I+\mathcal I'.
\end{align*}
To estimate $\mathcal I'$, we need to estimate
 $\int_{M_\rho}|k(\xi,z_1)|\, dV(\xi)$. So we
apply the approximating Heisenberg transformation
 $\zs=\psi_{z_1}(\zeta)$.
(However, we want to refrain from use of   Taylor remainder expansions, such as \re{ty12-}-\re{ty12},
 for the kernel. This
  avoids the requirement of $r\in \cC^{5/2}$.)
 Let $\zeta\in M_\rho$.
 Using the fundamental theorem of calculus, we get
 $\eta^n-y_1^n= A(\xi,x_1)\cdot(\xi-x_1)$ with $|A|<C$. By \re{xinx}, we
 have $\xi^n-x_1^n=p(\xi,x_1)\cdot\xi_*$ with $|p|<C$.
 Therefore,
 $$
 |\zeta-z_1|\leq C|\xi_*|.
 $$
   By
\rl{m2h}, we have
$$|r_\zeta-r_{z_1}|\leq C|\xi_*|, \quad
1/C\leq\f{|r_{z_1}\cdot(\zeta-z_1)|}{|\zs'|^2+i\xi_*^n|}\leq C, \quad
1/C\leq\f{|r_{\zeta}\cdot(\zeta-z_1)|}{|\zs'|^2+i\xi_*^n|}\leq C.
$$
  Therefore, for $\zeta\in M_\rho$
$$
|k(\zeta,z_1)|\leq
\f{C|\xi_*|^{|I|}}{|\zs'|^2+i\xi_*^n|^{a+b+c+d}}\leq\sum_{|J|=|I|}\f{ C|\xi_*^J|
 }{||\zs'|^2+i\xi_*^n|^{a+b+c+d}}.$$
By \rl{psii},   $\pi(\psi_{z_1}(M_\rho))\subset B_{9\rho}$.
 Applying \rl{outl}
with $\rho_1=\rho_0=9\rho$ and $\beta\geq0$, we get
   \eq{12rho1}
 \int_{M_\rho}|k(\zeta,z_1)|\,dV
\leq \sum_{|J|=|I|}\ \int_{B_{9\rho}}
\f{C|\xi_*^{J}|}{|\xi_*^n+i|\zs'|^2|^{a+b+c+d}}\,dV
\leq C'\rho.
   \end{equation}
Since $u(\zeta,z)$ is of class $\cC^{1/2}$ in the $z$-variable, one
gets
$$
|\mathcal I|\leq
C|z_2-z_1|^{1/2} \rho\sup_{\zeta\in D_\rho}\|u(\zeta,\cdot)
\|_{\rho,1/2}.
$$
To estimate the integral $\mathcal I$, it suffices to show that
$$
\int_{M_\rho} |k(\zeta, z_2)-k(\zeta, z_1)|\,dV\leq C \delta^{1/2},
\quad \delta=|z_2-z_1|.
$$

\

As mentioned in the introduction we  decompose $M_\rho$ into
 a {\it cylinder} and its complement.
 Consider the cylinder $ M_\rho\cap\{|\zeta'-z_1'|<\rho_1\}$,
where $\rho_1=C_*|z_2-z_1|^{1/2}$ with $C_*$ to be determined.
  Notice that the radius of cylinder is about
 $|z_2-z_1|^{1/2}$, which is much large than $|z_2-z_1|$,
  the euclidean distance
 between $z_1,z_2$.

We have
\gan
\int_{M_\rho} |k(\zeta, z_2)-k(\zeta, z_1)|\,dV\leq \mathcal I_1+\mathcal I_2+\mathcal I_3,\\
\mathcal I_1=\int_{M_\rho\cap\{|\zeta'-z_1'|<\rho_1\}}
|k(\zeta, z_1) |\,dV, \quad \mathcal I_2=\int_{M_\rho\cap\{|\zeta'-z_1'|<\rho_1\}}
|k(\zeta, z_2) |\,dV,\\
\mathcal I_3=\int_{M_\rho\cap\{|\zeta'-z_1'|\geq\rho_1\}}
|k(\zeta, z_2)-k(\zeta, z_1)|\,dV.
\end{gather*}
By an analogy of    \re{12rho1}, we have
$$ 
\mathcal I_1\leq \sum_{|J|=|I|}C\int_{|\zs'|<\rho_1,|\zs|<9\rho}
\f{|\xi_*^{J}|}{|\xi_*^n+i|\zs'|^2|^{a+b+c+d}}\,dV.
 $$ 
By \rl{outl} with $\beta\geq0$ and $\rho_0=9\rho$,   we get
$\mathcal I_1\leq C\rho_1=CC_*|z_2-z_1|^{1/2}$. Note that
$$
M_\rho\cap\{|\zeta'-z_1'|<\rho_1\}\subset M_{\rho}\cap\{|\zeta'-z_2'|<\rho_1+|z_2-z_1|\}.
$$
Applying the estimate  for $\mathcal I_1$, we get
 $\mathcal I_2\leq C(\rho_1+|z_2-z_1|)\leq C'C_* |z_2-z_1|^{1/2}.$

\

To estimate $\mathcal I_3$, we ignore any non-isotropic distance and
connect $(z_1',x_1^n)$ and $ (z_2',x_2^n)$ by a line segment
 in the convex domain $D_\rho$.
  Let $z(t)$ be the graph of the line segment
 in $M$.
Let $\zeta$ be a point of $M_\rho$ which is not in the cylinder. Then
$$|\zs'|=|\zeta'-z_1'|\geq \rho_1=C_*|z_2-z_1|^{1/2}.$$
 We have
\gan
|r_{z(t)}\cdot(\zeta-z(t))-r_{z_1}\cdot(\zeta-z_1)|\leq C_0|z_2-z_1|,\\
|r_{z_1}\cdot(\zeta-z_1)|\geq C_1|\zeta-z_1|^2\geq C_1|\zeta'-z_1'|^2\geq C_1C_*^2|z_2-z_1|,
\end{gather*}
where the second inequality comes from \rl{distb}.
Hence for  $C_1C_*^2>2C_0$, we obtain
$$
|r_{z_1}\cdot(\zeta-z_1)|/2\leq
|r_{z(t)}\cdot(\zeta-z(t))|\leq 2|r_{z_1}\cdot(\zeta-z_1)|.
$$
Recall that $\zeta_*=\psi_{z_1}(\zeta)$ and
$
C^{-1}||\zs'|^2+i\xi_*^n|\leq |r_{z_1}\cdot(\zeta-z_1)|\leq C||\zs'|^2+i\xi_*^n|.
$
Using $|r_{z(t)}\cdot(\zeta-z(t))|/C\leq
|r_\zeta\cdot(\zeta-z(t))|\leq C|r_{z(t)}\cdot(\zeta-z(t))|$, we get
\ga\label{esq1}
C^{-1}||\zs'|^2+i\xi_*^n|\leq
|r_\zeta\cdot(\zeta-z(t))|\leq C||\zs'|^2+i\xi_*^n|,\\
C^{-1}||\zs'|^2+i\xi_*^n|\leq
|r_{z(t)}\cdot(\zeta-z(t))|\leq C||\zs'|^2+i\xi_*^n|.\label{esq2}
\end{gather}

Write
$
k(\zeta,z)=\frac{p(\zeta,z)}{q(\zeta,z)}$ with
$
p(\zeta,z)=(\zeta'-z',\ov{\zeta'-  z'},r_{\zeta}-r_z)^I$ and
$$
q(\zeta,z)=(r_\zeta\cdot(\zeta-z))^{a}(r_z\cdot(\zeta-z))^{b}
\ov{(r_\zeta\cdot(\zeta-z))^{c}}\, \ov{(r_z\cdot(\zeta-z))^{d}}.
$$
For $\zeta\in M_\rho$ we have
\aln
|\zeta-z(t)|&\leq|\zeta-z_1|+|z(t)-z_1|\leq|\zeta-z_1|+C|z_2-z_1|\\
&\leq
|\zeta-z_1|+CC_*^{-2}|\zeta'-z_1'|^2\leq
C'|\zeta-z_1|. 
\end{align*}
By \re{xinx}, $|\xi^n-x_1^n|\leq C|\xi_*|$. Thus $|\zeta-z_1|\leq C|\xi_*|.$ Therefore,
$|\zeta-z(t)|\leq C|\xis|$ and
$$
| p(\zeta,z(t))|\leq C|\zeta-z(t)|^{|I|} \leq C'|\xi_*|^{|I|}.
$$
By \re{esq1}-\re{esq2}, we have
$$
|\{q(\zeta,z(t))\}^{-1}|\leq \f{C}{||\zs'|^2+i\xi_*^n|^{a+b+c+d}}.
$$
It is easy to see that $|z'(t)|\leq C|x_2-x_1|\leq C|z_2-z_1|$.
Then $|\pd_tq(\zeta,z(t))^{-1}|$ does not exceed the  sum of
$$
\f{C|z_2-z_1|}{|(r_\zeta\cdot(\zeta-z(t)))^{a'}(r_z\cdot(\zeta-z(t)))^{b'}
\ov{(r_\zeta\cdot(\zeta-z(t)))^{c'}}\, \ov{(r_{z(t)}\cdot(\zeta-z(t)))^{d'}}|},$$
where $a'+b'+c'+d'=a+b+c+d+1$. Applying the product rule, we get
$$|\pd_tp(\zeta,z(t))|\leq C|z_2-z_1|\sum_{|I'|=|I|-1}
|(\zeta-z(t),\ov{\zeta-z(t)},r_{\zeta}-r_{z(t)})^{I'}|.
$$
Therefore,
$$|\pd_tp(\zeta,z(t))|\leq C|z_2-z_1| |(\zs',\xi_*^n)|^{|I|-1},\  \bigl|\pd_t
 \f{1}{q(\zeta,z(t))}\bigr |\leq \f{C|z_2-z_1|}{||\zs'|^2+i\xi_*^n|^{a+b+c+d+1}}.
$$
By the mean-value-theorem, we get for $\zeta\in M_\rho\cap\{|\zeta'-z'|>\rho_1\}$
\aln
&|k(\zeta,z_2)-k(\zeta,z_1)|\leq C|z_2-z_1| \Bigl\{
\f{|\xi_*|^{|I|}}{||\zs'|^2+i\xi_*^n|^{a+b+c+d+1}}
 +
\f{|\xi_*|^{|I|-1}}{||\zs'|^2+i\xi_*^n|^{a+b+c+d}}
\Bigr\}\\
&\hspace{5em}\leq C'|z_2-z_1|
\f{|\xi_*|^{|I|}}{||\zs'|^2+i\xi_*^n|^{a+b+c+d+1}}
\leq \sum_{|J|=|I|}
\f{C'|z_2-z_1||\xi_*^J|}{||\zs'|^2+i\xi_*^n|^{a+b+c+d+1}}.
\end{align*}
Applying \rl{outl}  with $\rho_1=C_*|z_2-z_1|^{1/2}$, $\rho_0=9\rho$ and $\beta\geq-2$,
we get
\begin{equation*}
\mathcal I_3\leq C|z_2-z_1|\rho_1^{-1}\leq C'|z_2-z_1|^{1/2}.   \qedhere
\end{equation*}
 \end{proof}

\setcounter{thm}{0}\setcounter{equation}{0}\section{Proof of \rt{1}}
\label{sec19}

We now prove \rt{1}, following a KAM argument
in~\ci{GWzese}.
We restate the theorem.

\medskip

\noindent{\bf Theorem.\, }{\em
Let $M\colon r=0$ be a strongly pseudoconvex real hypersurface of class $\cC^2$
in $\cc^n$ with $n\geq4$.
Let $\omega$
 be a continuous $r\times r$   
  matrix
of $(0,1)$-forms on $M$ satisfying the integrability condition
$\dbb\omega=\omega\wedge\omega\mod\db r$. Near each point of $M$ there exists a non-singular
matrix $A\in \cC^{1/2}(M)$ such that $\dbb A=-A\omega\mod\db r$.
Moreover, if $a$ is a positive real number, $M$ is of class $\cC^{a+2}$ and
$\omega\in \cC^{a}(M)$   there is a solution
$A\in \cC^{a}(M)$;  if $k$ is a positive integer,  $\omega\in \cC^k$ and $M\in
\cC^{k+\f{5}{2}}$,  there is a solution $A\in
\cC^{k+\yt}(M)$.
}

\medskip

Note that not all solutions have the same regularity.  If $u$ is a
continuous CR function vanishing nowhere on $M$,
then $uA$ is still  a solution.

\

\newcommand{\mh}[1]{\noindent{\bf #1.}}

\mh{Non-isotropic dilations}
The non-isotropic dilation $T_\delta\colon z\to(\sqrt\delta z',\delta z^n)$ with $\delta>0$
does
not preserve the real hypersurface $M$.
However, it is obvious that it sends   $M^\delta=T_\delta^{-1}M$ onto $M$.
By abuse of notation, we will denote $T_\delta$ its restriction  to $\cc^{n-1}\times\rr$.

We want to show that $M$ can be put in the form \re{bhroo-}-\re{bhroo} after
a second order normalization and a non-isotropic dilation.

Let
$M$ be a graph  $y^n=|z'|^2+\hat r(z',x^n)$ over $D$.
Recall that $M_\rho=M\cap\{(x^n)^2+y^n<\rho^2\}$ and
$$D_\rho=\pi(M_\rho)=\left\{(z',x^n)\in D\colon |z'|^2+|x^n|^2+\hat r(z',x^n)
<\rho^2\right\}.$$
Set $r(z)=-y^n+\hat r(z',x^n)$.
Define
$$\hat r^\delta(z)=\delta^{-2} \hat r( \delta z',\delta^2 x^n),
\quad r^\delta(z)=\delta^{-2} r(  \delta z',\delta^2 x^n).
$$
Then $M^\delta=T_\delta^{-1} M$ is the graph  over
 $D^{\delta}=T_\delta^{-1}D$, given by
$$
y^n=|z'|^2+\hat r^\delta(z',x^n).
$$
For $ 0<\delta<1$, we have
$$
M^\delta_\rho=M^\delta\cap\{|x^n|^2+y^n<\rho^2\}
\subset M^\delta\cap\{\delta^2|x^n|^2+y^n<\rho^2\}.
$$
So $D^\delta_\rho\df\pi(M_\rho^\delta)\subset T_\delta^{-1}D_\rho$ for $0<\delta<1$.

In \rt{1}, we need a local solution $A$. By a change of
local holomorphic coordinates, we may assume that $\hat r(0)=\pd\hat r(0)=
\pd^2\hat r(0)=0$. Therefore
$$
 \|\hat r^\delta\|_{D^\delta_1,2}<\epsilon,
\quad \ov{D_1^\delta}\subset D^\delta,
$$
if $\delta$ is sufficiently small.

\

\mh{Integrability conditions} We will find our solution through
a sequence of frame changes.  We also need a small norm of initial $\omega$ via
dilation.
Therefore, we need to verify that the integrability
condition is preserved under dilation and frame changes.

\medskip

Recall that for  $\varphi=
\sum_{|I|=q}\varphi_{\ov I}d\ov{z'}^I$,
we define $\dbm\varphi=\sum_{|I|=q,1\leq\alpha<n} X_{\ov\alpha}
\varphi_{\ov I}d\ov{z^\alpha}\wedge d\ov{z'}^I$
for $X_{\ov\alpha}=\partial_{\ov{z^\alpha}}-{r_{\ov{z^\alpha}}}/{r_{\ov{z^n}}}\partial_{\ov{z^n}}$.
A direct computation shows that $\delta dT_\delta^{-1}X_{\ov\alpha}=
X_{\ov\alpha}^\delta\equiv \partial_{\ov{z^\alpha}}-
{r^\delta_{\ov{z^\alpha}}}/{r^\delta_{\ov{z^n}}}\partial_{\ov{z^n}}$.
Thus,  $T_\delta^*\ov\pd_M=\ov\pd_{M^\delta}T_\delta^*$.
This shows that the formal integrability condition
is   invariant under dilation, i.e.
$\dbar_{M^\delta}\omega^\delta=\omega^\delta
\wedge\omega^\delta$.
If $\omega=\sum\varphi_{\ov J}\, d\ov z^J$ is a tangential $(0,1)$-form
on $M$, then  $T_\delta^*\omega$
is a tangential $(0,1)$-form on $M^\delta$. Moreover, if $\omega\in \cC^a(M)$, then
\eq{tdsw}
T_\delta^*\omega=\delta \sum_{\alpha=1}^{n-1}\varphi_{\ov\alpha}\circ T_\delta\, d\ov{z^\alpha},\quad
\lim_{\delta\to0}\|\omega^\delta\|_{\cC^a(D^\delta_1)}=0.
\end{equation}
In other words, we have achieved the smallness of $\omega$ via dilation alone.

We now consider the integrability condition under a frame change.
We are given an $r\times r$ matrix of
continuous $(0,1)$-forms $\omega$
 on  $M$.
We assume that $\dbb\omega=\omega\wedge\omega \mod{\ov\partial r}$.
   Without loss of
generality, we may assume that
$\omega$ are tangential.
The
formal integrability condition is that as
currents
\eq{intc}\dbm\omega=\omega\wedge\omega.
\end{equation}
Recall that our goal is to  find a non-singular matrix $A$ which
solves
\eq{dbma}
\dbm A+A\omega=0.
\end{equation}
We will consider only the solution $A $ such that both $A$ and
$\dbm A$ are continuous.
For any such $A$,  the transformation
$\omega\to \tilde\omega=(\dbm A+A\omega)A^{-1}$ preserves the integrability condition
\re{intc}.    
Indeed,  differentiating  $\tilde\omega A=\dbm A+A\omega$  and then
using $\omega=A^{-1}\tilde\omega A-A^{-1}\dbm A$,
$\dbm\omega=\omega\wedge\omega$, we verify $\dbm\tilde\omega=\tilde\omega \wedge\tilde\omega$ by
\aln
(\dbm\tilde\omega) A&-\tilde\omega\wedge \dbm A
=\dbm A\wedge\omega+A\dbm\omega=\dbm A\wedge (A^{-1}\tilde\omega A-A^{-1}\dbm A)\\
&+(\tilde \omega A-\dbm A)\wedge(A^{-1}\tilde
\omega A-A^{-1}\dbm A)=\tilde\omega\wedge\tilde\omega A-\tilde\omega\wedge\dbm A.
\end{align*}

Assume that the    matrix $\omega$
  is of class
$\cC^a$.
In what follows, all constants, including $\delta$,
 will depend on $a$.  For simplicity this dependence
will not be indicated sometimes. However,
  constants $C_0, \delta_*$ and $\epsilon$ do not depend on $a$.

\

\mh{Proof of \rt{1}}
We need to find a non-singular matrix $A=I+B$,    
defined
near the origin of $M$,   
such that
$$
\dbm B+\omega+B\omega=0.
$$
 It suffices to find a $\delta>0$ and a
non-singular matrix $A^\delta$ defined near $0\in M^\delta$
such that $\ov\partial_{M^\delta}\omega^\delta
+A^\delta\omega^\delta=0$. Then $A^\delta\circ T_\delta^{-1}$ is a solution
to the original equation.

\

 Take $\rho_0=1$, $\sigma_j=2^{-j-1}$ and $\rho_{j+1}=(1-\sigma_j)\rho_j$.
 Then $\rho_\infty=\lim_{j\to\infty}\rho_j>0$. We will assume that
 $0<\delta\leq\delta_*$.
We want to apply our estimates for $P'$ and $Q'$.
So we choose $\delta_*\in(0,1]$ such that the \hf\ holds on $M^\delta_\rho$ for $
\rho\in(0, 1]$ and $\delta\in(0,\delta_*]$. For $\rho_\infty<\rho\leq 1$ we have
\eq{estj}
\|P'\varphi\|_{\cC ^{a}(M^\delta_{(1-\sigma)\rho})} \leq
C \sigma^{-s}
\|\varphi\|_{\cC ^{a+2}(M^\delta_{\rho})},
\end{equation}
where $P'$ is either of operators $P',Q'$ in the \hf\ on $M^\delta_\rho$.
We emphasize that the constant $C_{a}$ is independent of $\delta\in(0,\rho_*)$ and
$\rho\in(\rho_\infty,1]$. We have also absorbed
$\|r^\delta\|_{\cC^{a+2}(D_{\rho_j}^\delta)}$ into $C_a$,
since $\hat r^\delta(z',x^n)=\delta^{-2}\hat r(\delta z',\delta^2 x^n)$
and $\hat r(0)=
 \pd\hat r(0)=0$ imply that
$$
\|\hat r^\delta\|_{\cC^{a+2}(D_{\rho_j}^\delta)}\leq C\|\hat r\|_{\cC^2(D_{\delta\rho_0})}
+\delta^a\|\hat r\|_{\cC^{a+2}(D_{\delta\rho_0})}<C_a, \quad 0<\delta<\delta_*.
$$

Set $M_j=M^\delta_{\rho_j}$
and
$\|\cdot\|_{\rho_{j+1},a}=\|\cdot\|_{\cC^{a} (M_{ {j+1}})}$.
We have $M_{j+1} \subset M_0$.
Let $\omega_0=T_\delta^*\omega$, restricted on $M_0$.
On $M_j$  we have the \hf\  $\varphi=
\db_{M^\delta} P'_j\varphi+Q'_j\db_{M^\delta}\varphi$, where
 $P'_j=P'_{M_j}$ and $Q'_j=Q'_{M_j}$.

 Using $\db_{M^\delta}\omega_0
=\omega_0\wedge\omega_0$, we arrive at the equation
$$
\db_{M^\delta} (B_0+P'_0\omega_0)+Q'_0(\omega_0\wedge\omega_0)+B_0\omega_0=0
$$
where $P'_0, Q'_0$ are applied  entrywise to the matrices. We use the approximate
solution
$$B_0=-P'_0\omega_0.$$ Assume that $A_0=I+B_0$ is invertible.
We repeat this procedure and
get $B_j=-P'_j\omega_j$ and
\al
\label{toqof}\omega_{j+1}&=(\dbm A_j+A_j\omega_j)A_j^{-1} =
 \{Q'_j(\omega_j\wedge\omega_j)-(P'_j\omega_j)\omega_j\}
(I-P'_j\omega_j )^{-1}.  
\end{align}
Here, we need
all $A_{j}=I+B_{j}$ to be non-singular on $M_j$. We want to show that
 $\lim_{j\to\infty}
 A_jA_{j-1}\cdots A_0$
is a solution.

We now estimate   $\|B_j\|_{\rho_{j+1},a}$
and $\|\omega_{j+1}\|_{\rho_{j+1},a}$.

For an $r\times r$
matrix $B=(b_i^j)$ of functions on $M_\rho$   we define $\|B\|_{\rho,a}
=\max\{\|b_{i}^j\|_{\rho,a}\}$.
If $B,D$ are two such matrices, we have
$$
\|DB\|_{\rho,0}\leq r^{2}\|D\|_{\rho,0}\|B\|_{\rho,0},\quad
\|B^l\|_{\rho,0}\leq r^{l}\|B\|_{\rho,0}^{l}.
$$
Assume that $\rho_\infty<\rho\leq\rho_0$.
We want to show that if $\|B\|_{ \rho,0}\leq\f{1}{2r}$, then
\eq{invm}
\|D\|_{\rho,a}
\leq c_{a}\|B\|_{\rho,a}, \quad 0\leq a<\infty,
\end{equation}
where $c_a>1$ depends on $a, r,n$.
 Let $A^{-1}=I+D$. We
know that $D=\sum_{l\geq1}(-1)^{l}B^l$ and
\gan
\|D\|_{ \rho,0}\leq
\f{r\|B\|_{ \rho,0}}{1-r\|B\|_{ \rho,0}} \leq
2r\|B\|_{ \rho,0},
\end{gather*}
which is \re{invm} with $a=0$. Since $I$ is constant, then
\aln
D(z_2)-D(z_1)&=A(z_2)^{-1}-A(z_1)^{-1}=A(z_1)^{-1} (A(z_1) -A(z_2) )A(z_2)^{-1}\\
&=A(z_1)^{-1}(B(z_1)-B(z_2))A(z_2)^{-1}.
\end{align*}
Using \re{invm} with $a=0$, we get
$\|D\|_{\rho,a}\leq C \|B\|_{\rho,a}$ if $0<a\leq 1.$
Assume that \re{invm} holds
when $[a]<k$. Let $[a]=k\geq1$.
Applying the product rule to $(I+B)D=-B$
and multiplying from left by $(I+B)^{-1}=I+D$, we get
 \aln
\pd D
&=(I+D) (\pd B)D - (I+D) \pd B.
\end{align*}
By  \rp{hoin}  and the induction assumption, we obtain
$$
\|D\|_{\rho,a}\leq C_a\{(1+\|B\|_{\rho,a-1})
\|B\|_{\rho,1}+\|B\|_{\rho,a}\}\leq C_a'\|B\|_{\rho,a},
$$
which proves \re{invm}.

Since $B_j=-P'_j\omega_j$, by \re{estj} we have
\eq{bwj}
\|B_j\|_{\rho_{j+1},a} \leq
c_a^*\sigma_j^{-s_a}\|\omega_j\|_{\rho_j,a}, \quad c_a^*>1.
\end{equation}
We want to achieve $\|B_j\|_{ \rho_{j+1},0} \leq\f{1}{2r}$.
So it suffices to obtain
\eq{bjj}
\|\omega_j\|_{\rho_j,a} \leq\f{\sigma_j^{s_a}}{2rc_a^* c_a}=b_j, \quad j=0, 1,2, \ldots.
\end{equation}
By \re{invm}-\re{bjj}, we have $\|(I+B)^{-1}\|_{\rho_j,a}\leq 1+\|D\|_{\rho_j,a}\leq2$.
Using \re{estj}-\re{toqof} and the estimates on
   matrix products,  we get
\al\label{wjj1} \|\omega_{j+1}\|_{\rho_{j+1},a } &\leq
2\cdot r^2
\bigl\{\|Q'_j(\omega_j\wedge\omega_j)\|_{\rho_{j+1},a }
\bigr.\\ &\quad
\bigl.+r^2\|B_j\|_{\rho_{j+1},a }\|\omega_j\|_{\rho_{j+1},a }\bigr\}
\leq { C_a^*}{\sigma_j}^{-s_a}\|\omega_j\|_{\rho_j,a }^2.\nonumber
\end{align}
Assume that $\|\omega_0\|_{\rho_0,a}\neq0$. Otherwise the theorem
holds trivially.
Define
$$
\hat b_{j+1}= { C_a^*}{\sigma_j}^{-s_a}\hat b_j ^2,\quad \hat b_0=
\|\omega_0\|_{\rho_0,a }.
$$
  By \re{tdsw},
we choose   a dilation $T_\delta$ with $\delta\in(0,\delta_*]$ such that $\omega_0=
T_\delta\omega$ satisfies
$$
\hat b_0=\|\omega_0\|_{\rho_0,a}\leq b_0.
$$
Then  $r_j=\hat b_{j+1}/{\hat b_j}$ satisfies
\gan
r_0=  { C_a^*}{\sigma_0}^{-s_a}\hat b_0, \quad r_{j}=
2^{s_a}(r_{j-1})^2,\\
r_1/r_0=2^{s_a}r_0, \quad r_{j+1}/r_j=(r_j/r_{j-1})^2.
\end{gather*}
This shows that $r_{j+1}/r_j$, and hence $r_j, \hat b_j$,
converge rapidly, if $\hat b_0$ is sufficiently small. Specifically,
$$
r_{j+1}/r_j=(r_1/r_0)^{2^j}, \quad r_j=r_0 (r_1/r_0)^{2^j-1}, \quad
\hat b_j=\hat b_0 r_0^j(r_1/r_0)^{2^j-j-1}.
$$
Recall that $\sigma_j=2^{-j-1}$.
Clearly, $\|\omega_j\|_{\rho_j,a}\leq\hat b_j\leq \hat b_0r_0^j\leq
\f{2^{-s_a}}{2rc_a^*c_a}(2^{-s_a})^j=b_j$, provided
$$
\hat b_0\leq \f{2^{-s_a}}{2rc_a^*c_a}, \quad r_0\leq 2^{-s_a},
\quad r_1/r_0=2^{s_a}r_0\leq1.
$$
Therefore, we have shown that if $\omega\in \cC^a$ and
$M\in \cC^{2+a}$ then   $(I+B_j)\cdots (I+B_0)$
converges in $\cC^a$-norm on $\cap D_{\rho_j}^\delta$
 to an invertible matrix $A_\infty$.
When $a=k$ is an integer and $M\in \cC^{k+5/2}$,  or when
$M\in \cC^{2}$
for $k=0$, we have
$$
\|B_j\|_{\rho_{j+1},k+1/2} \leq
c_{k} \sigma_j^{-s_{k}}\|\omega_j\|_{\rho_j,k}.
$$
Since $\|\omega_j\|_{\rho_j,k}$ tends to zero rapidly, it is obvious
that $B_j$ converges to $0$ rapidly in $\cC^{k+1/2}$ on
$\cap D_{\rho_j}^j$. This shows that $A^\infty$ is of class $\cC^{k+1/2}$.

To complete the proof,  we  check that
$\db_{M^\delta} A^\infty+A^\infty\omega_0=0$. Recall that for $A_j=I+B_j$,
we have
$\db_{M^\delta} A_j+A_j\omega_j=w_{j+1}A_j$.
Thus
\aln
\omega_{j+1}A_jA_{j-1}&=(\db_{M^\delta} A_j)A_{j-1}+A_j\omega_jA_{j-1}\\
&=(\db_{M^\delta} A_j)A_{j-1}+A_j(\db_{M^\delta} A_{j-1}+A_{j-1}\omega_{j-1})\\
&=\db_{M^\delta}(A_jA_{j-1})+A_jA_{j-1}\omega_{j-1}.
\end{align*}
 Inductively,
we get $\db_{M^\delta}(A_j\cdots A_0)+A_j\cdots A_0\omega_0=\omega_{j+1}A_j\cdots A_0.$
Taking the  limits, we get $\db_{M^\delta} A^\infty+A^\infty\omega_0=0$ on
$M^\delta_{\rho_\infty}$.
When $k=0$,  the   derivatives in the sense of currents
are   continuous and the above computation is valid as
currents. \qed

\

In the above argument,
we obtain the rapid convergence of $B_j, \omega_j$ in $\cC^a$ norm in one step
when $a$ is finite. One can also establish a rapid convergence of $\omega_j, B_j$
first in $\cC^0$-norm and then in higher order derivatives. See~\ci{Wenion},
\ci{GWzese}
for  details.

\appendix

\setcounter{thm}{0}\setcounter{equation}{0}\section{H\"older inequalities}
\label{hoine}
The main purpose of this appendix is to present some H\"older
inequalities on domains in $\rr^m$.
We do not claim any originality in deriving these inequalities. In fact, we will
just modify formulation and proofs  of H\"ormander~\ci{Hosesi}.
The inequalities  in~\ci{Hosesi} are for a fixed convex domain.
In our applications, we need to allow
the domain $D_\rho$ to vary. Therefore,
we will derive them in full details
and omit   simple repetitions only.

\

We say that a domain $D$ in $\rr^m$ has the
{\it cone  property} if the following
hold:
(i)
 Given two points $p_0,p_1$ in $D$ there exists
a piecewise  $\cC^1$   curve $\gamma(t)$ in $D$
such that $\gamma(0)=p_0$ and $\gamma(1)=p_1$,
$|\gamma'(t)|\leq   C_*\|p_1-p_0\|$ for all $t$
except  finitely many   values. The diameter of $D$ is
less than   $C_*$.
(ii)
For each point $x\in \ov D$, $D$ contains a cone $V$  with vertex $x$,
opening $\theta>C_*^{-1}$ and height $h>C_*^{-1}$.

We will denote $C_*(D)$ a  constant $C_*>1$  satisfying
 (i) and (ii).

In this appendix,  by  a   cone   $V=V(\theta,h,v)$   with vertex at the origin,
opening $\theta>0$
and height $h>0$, and centered at positive $v$ axis where   $v$ is
a unit vector, we mean 
\eq{ttmt}
V=\{t\in\rr^m\colon v\cdot t>\theta^{-1}\|t-(v\cdot t)v\|,   v\cdot t<h\},
\end{equation}
where $\|t\|=(|t^1|^2+\cdots+|t^{m}|^2)^{\yt}$.
Note that $x+V$
is  a cone with vertex at $x$.
If $V_m\subset\rr^m$ is the cone $1\geq x^m>\|x'\|$,
a simple computation shows that
$V_m\times V_n\subset\rr^{m+n}$ contains
the cone
$$
1\geq x^m+y^n>(|x^m-y^n|^2+4\|x'\|^2+4\|y'\|^2)^{1/2}.
$$
Thus  if $D_1$, $D_2$ are domains of \ps, then $D_1\times D_2$ has
\ps\ with a constant $C_*(D_1\times D_2)$ depending only on $C_*(D_i)$.

\

For the rest of the appendix, until \rp{vbhi},
we assume that the domain  $D$
has \ps\ unless stated otherwise.
The constants in all H\"older inequalities will
depend on $m$ and $C_*(D)$.  For simplicity this dependence
will not be expressed sometimes.

\

Let $k\geq0$ be an integer.
For a complex-valued function  $u$    on $D\subset \rr^{m}$, define
\gan
\|\partial^ku\|_{D,0}= \sup_{x\in D,|I|=k} |\partial^Iu(x)|,
\quad \|u\|_{D,k}=\max_{0\leq j\leq k}\|\partial^ju\|_{D,0},\\
|u|_{D,\alpha}= \sup_{x,y\in D}\frac{|u(x)-u(y)|}{|x-y|^\alpha}, \quad 0< \alpha\leq1,\\
\|u\|_{D,k+\alpha}=\max\{
 \|u\|_{D,k}, |\partial^Iu|_{D,\alpha}\colon |I|=k\}, \quad 0<\alpha<1.
\end{gather*}
If $A=(a_i^j)$ is a 
matrix of functions on $D$, we define $\|A\|_{D,k+\alpha}
=\max\{\|a_{i}^j\|_{D,k+\alpha}\}$. 

We also define
$$
|u|_{D,k+\alpha}=\max_{|I|=k}|\partial^Iu|_{D,\alpha},\quad 0<\alpha\leq1.
$$
By (i) of \ps\ and
 the fundamental theorem of calculus
\eq{ftcn}
C_0^{-1} |u|_{D,k}\leq \|\partial^ku\|_{D,0}\leq |u|_{D,k}, \quad k=1,2,\ldots,
\end{equation}
  provided that $u\in \cC^k(D)$.
It will be convenient to use both $|u|_{D,k}$ and $\|\partial^ku\|_{D,0}$.

\le{c2ps} Let
 $D\subset \rr^m$ satisfy (i) of the cone property.
Let $f$ be a $\cC^1$ map from  $D$ into $\rr^m$.
 There exists a constant $C_0>1$
such that  if
  $|f'-I|<C_0^{-1}$ on $D$,
then $f$ is a $\cC^1$ diffeomorphism from
$D$ onto $D'$. Moreover,
 $\|f^{-1}-I\|_{D',1}\leq C\|f-I\|_{D,1}$.   \end{lemma}
  \pf Take
two points $p_0,p_1$ in $D$. By   assumption
there exists a piecewise $C^1$ curve $\gamma$ in $D$ such
that $\gamma(0)=p_0$, $\gamma(1)=p_1$, and
$|\gamma'(t)|\leq C_*\|p_1-p_0\|$.
 Write $f=I+\tilde f$.  We have
 $$
  f(p_1)-f(p_0)= p_1-p_0 - \int_0^1\nabla \tilde f(\gamma(t))
  \cdot\gamma'(t)\, dt.
  $$
 Since $\| \tilde f'\|_{D,0}<C_0^{-1}$, then for $C_0$ sufficiently
 large we get
\al\label{fp1p0}
\yt \|p_1-p_0\|\leq
\|f(p_1)-f(p_0)\|\leq2\|p_1-p_0\|.
\end{align}
   Now it is obvious that $f$ is a $\cC^1$ diffeomorphism from $D$ onto $D'$ and that
$D'$ satisfies (i) of \ps.
Let $\tilde g=f^{-1}-I$. Then $\tilde g\circ f=-\tilde f$. In particular,
$\|\tilde g\|_{D',0}
\leq \|\tilde f\|_{D,0}$, and by \re{fp1p0}, $|\tilde g|_{D',1}
\leq C|\tilde f|_{D,1}$.
 \end{proof}


\

We are ready to derive   H\"older inequalities.
We start with two lemmas in~\ci{Hosesi} for our domains.
\le{a1}
Let $0<a<b$. Let $D$  satisfy   \ps.
If $|u|_{D,a}\leq 1$ and $|u|_{D,b}\leq1$ then $|u|_{D,c}
\leq C$ for $a<c<b$, where $C$ depends only on $a,b$ and $C_*(D)$.
\end{lemma}
\pf
For simplicity, denote $|u|_{D,a}$ by $|u|_a$.
Since the diameter of $D$ is bounded by some constant $C$,
it is obvious that $|u|_c\leq C'|u|_{c'}$ if $k<c<c'\leq k+1$.
 Therefore, it suffices to prove the inequality when $c$ is an integer.

Let $V\subset D$ be a cone  as stated  in   \ps. We may assume that $0\in V$.

If $a$ is an integer,  we get $|\partial^a u|\leq1$ on $V$. If
$a$ is not an integer, let $P$ be its Taylor polynomial of degree $[a]$ at $0$.
Set $v=u-P$. Since $\partial^IP$ are constants   for all $|I|\geq [a]$,
then $|v|_{c'}=|u|_{c'}$ for all $c'>[a]$.
Therefore, we may assume that
 $\partial^Iu(0)=0$ for all $|I|=[a]$. Now
$|u|_a\leq1$ implies that $|\partial^Iu|<C_0$ for all $|I|=[a]$.

We want to use the mean-value-theorem repeatedly. Let us first look
at the one-variable case.
When $f$ is $\cC^k$ on $[0,1]$ and $\|f\|_0\equiv\|f\|_{[0,1],0}< C_0$,
there is a point $t$
in $[0,1]$ such that $|f'(t)|\leq 2C_0$. One can divide $[0,1]$ into a
sufficient number
of  equal parts and find
a point $t_j\in[0,1]$ such that $|f^{(j)}(t_j)|<C'$ for $j=1,\dots, k$.
If $|f|_{a}\leq1$ for some $a\in(k,k+1]$,
we obtain further that $\|f^{(k)}\|_0<C_k$. Consequently, $\|f^{(j)}\|_0<C_j$ hold
 for $j=k, k-1,\dots,  0$.

Return to our case. Fix a polydisc $\Delta^m$ in $V$ with side larger than $C^{-1}$.
Fix $I$ with $|I|=[a]$. Using $|\pd^Iu|<C_0$ for $|I|=[a]$ and $|u|_b\leq1$,
by the   one-variable argument we obtain $ |\partial^j\partial^Iu |
<C$ on $\Delta^m$
for $j=[b]-[a],\dots, c-[a]$. So $ |\partial^{j}u | <C$ on $\Delta^m$ for
$j=c,\dots, [b]$.  If   $[b]<b$,
  $|u|_b\leq1$ implies that $ |\partial^{[b]}u |<C'$ on $D$.
  If $b=[b]$ the assumption and \re{ftcn} implies
   $|\pd^{b}u|\leq1$. Using a path connecting
 a point in $D$ to $\Delta^m$ and $|\pd^ju|<C$ on $\Delta^m$ for $j=[b],\ldots, c$, we
  get   $|\pd^ju|\leq C''$ on $D$ for
 $j=[b]-1,\ldots,c$.  Using \re{ftcn} again, we get
  $|u|_c\leq C$.
\end{proof}

\noindent{\bf Example.} Let $f(x)=x+x^3$ and $D=[0,\epsilon]$. Then
$f'(0)=f'(x)-\int_0^xf''(t)\, dt$
and
$$
\epsilon f'(0)=\int_0^\epsilon f'(x)\, dx-\int_0^\epsilon\int_0^xf''(t)\, dt\, dx
=f(\epsilon)-f(0)-\int_0^\epsilon\int_0^xf''(t)\, dt\,  dx.
$$
This  shows that
$
1=|f'(0)|\leq C_1\|f\|_0+C_2\|f''\|_0\leq 2\epsilon C_1+6\epsilon C_2.
$
However, $|C_1|+|C_2|$ tends to $\infty$ as $\epsilon\to0$.
  This example demonstrates that in some
 inequalities derived in this appendix,
   the constants indeed depend on $C_*(D)$.
For special domains used in this paper, we will find these
constants by   dilation; see \rp{vbhi}.

\le{a2} Let $D\subset\rr^m$ satisfy the cone property. Set $\|\cdot\|_{D,a}
=\|\cdot\|_a$.
If $0<a<b$, $c=\lambda a+(1-\lambda)b$ and $0<\lambda<1$, then
$
|u|_{c}\leq C|u|_{a}^\lambda(|u|_{a}+|u|_{b})^{1-\lambda},
$ 
where $C$ depends only on $a,b$.
\end{lemma}
\pf
The case $|u|_a\geq |u|_c$ is obvious and we may assume that $|u|_a\leq|u|_c$.
If $|u|_a=0$, then $u$
is a polynomial of degree $<a$. Then $|u|_c=0$ too for $c>a$, and the inequality
holds.
We may assume that $|u|_a\neq0$.  Without loss of generality,
we may assume that $|u|_a=1<|u|_c$.
If $|u|_b\leq1$, \rl{a1} implies that $|u|_c\leq C$ and the inequality holds.
Therefore, it suffices to verify
 \eq{uccu}
|u|_c\leq C|u|_b^{1-\lambda}, \quad
\text{if $|u|_a=1<|u|_b$.}
 \end{equation}
 We first assume the inequality for integer $c$ and  verify it
for non-integer $c$. Set $[c]=k$.   We have $\lambda=\f{b-c}{b-a}$
and $1-\lambda=\f{c-a}{b-a}$.
Depending on whether $a,b$ are in $[k,k+1]$,  we have
the following cases.

{\bf Case i)
 $k\leq a<c<b\leq k+1$.}
 Since $c=\lambda a+(1-\lambda)b$ then $c-k=\lambda(a-k)+(1-\lambda)(b-k)$.
 Consider first the case
  $a=k$. Then  $c-k=(1-\lambda)(b-k)$. Let $v=\partial^Iu$ with $|I|=k$.
 Hence
 \aln
 \f{| v(y)- v(x)|}{|y-x|^{c-k}}&=| v(y)- v(x)|^\lambda
 \left| \f{v(y)-v(x)}{|y-x|^{b-k}}\right|^{1-\lambda}\leq 2^\lambda
 \| v\|_0^\lambda | v|_{b-k}^{1-\lambda}\\
 &
 \leq C'|u|_k^\lambda |\pd^Iu|_{b-k}^{1-\lambda}\leq C'|u|_k^\lambda |u|_{b}^{1-\lambda},
 \end{align*}
where the  second
 last inequality is obtained by \re{ftcn} using $k=a>0$.  Therefore,
$|u|_c\leq C |u|_a^\lambda|u|_b^{1-\lambda}$.
  Assume now that $a-k>0$. Then $a-k,b-k,c-k$ are in $(0,1]$.
  Computing the H\"older
 ratio
 gives us $| v|_{c-k}\leq | v|_{a-k}^\lambda| v|_{b-k}^{1-\lambda}$, i.e.
 $|u|_c\leq |u|_a^\lambda|u|_b^{1-\lambda}$.
 We emphasize that we have proved
 $|u|_c\leq C|u|_a^\lambda|u|_b^{1-\lambda}$ (and
 hence  \re{uccu} for case i))
 without using any
  condition on $|u|_a,|u|_b,|u|_c$ other than  that on $a,b,c,k$.

{\bf Case ii) $a<k<c<b\leq k+1$.}
We get
$
|u|_c\leq C|u|_k^{(b-c)/{(b-k)} 
}|u|_b^{(c-k)/{(b-k)} 
}$, by case i).   
By the assumption for the integer case  (applied to triple
$a<k<b$) we obtain
\eq{int1}
|u|_{k}\leq C
|u|_b^{(k-a)/{(b-a)}
}.
\end{equation}
Eliminating $|u|_{k}$ from two inequalities gives us \re{uccu}; indeed
$$
\f{k-a}{b-a}\cdot\f{b-c}{b-k}+\f{c-k}{b-k}=
\f{c-a}{b-a}.
$$

{\bf Case iii) $k\leq a<c<k+1<b$.}
We have
$
|u|_c\leq C
|u|_{k+1}^{(c-a)/{(k+1-a)} 
}$,
by case i).  
The assumption on the integer case (applied to  triple
$a<k+1<b$) gives us
\eq{int2}
|u|_{k+1}\leq C
|u|_b^{(k+1-a)/{(b-a)} 
}.
\end{equation}
Eliminating $|u|_{k+1}$  gives us \re{uccu}.

{\bf Case iv) $a<k<c<k+1<b$.} Then \re{int1} and \re{int2}  are  still valid.
By i), we have
$
|u|_c\leq C|u|_k^{k+1-c}|u|_{k+1}^{c-k}.
$  
Using \re{int1}-\re{int2} and eliminating $|u|_k$, $|u|_{k+1}$
gives us \re{uccu}.

\

Finally, we  prove \re{uccu} when $c$ is a positive integer by repeating an
argument in~\ci{Hosesi}.

Fix $x_0\in D$ and let $V$ be a cone in $D$ with vertex $x_0$,  height   and
opening $1/C_*$.
Since $V$ is convex,  for
$x\in V$ and
$0<\epsilon<1$ we can define
$$
u_x^\epsilon(y)=u((1-\epsilon)x+\epsilon y), \quad y\in V.
$$
Then $|u_x^\epsilon|_{V,b}\leq\epsilon^b|u|_b$
and $|u_x^\epsilon|_{V,a}\leq\epsilon^a$. Since $|u|_b>1$, there is
an $\epsilon\in(0,1)$ so that $\epsilon^a=\mu=\epsilon^b|u|_b$.
Now,  apply \rl{a1}   to  the domain $V$ and the function $\mu^{-1}u_x^\epsilon$.
For any multiindex $I$ with $|I|=c$, we have
$$
\epsilon^c|\partial^Iu(x)|=|\partial^Iu_x^\epsilon(x)|\leq C_0\mu=C_0(\epsilon^a)^\lambda
(\epsilon^b|u|_b)^{1-\lambda}.
$$
Canceling $\epsilon$'s  shows $|\partial^Iu(x)|\leq C_0|u|_b^{1-\lambda}$
for $x\in \ov V$.  Since $C_0$ does not depend on $x_0\in\ov V$,
we get $\|\pd^cu\|_0\leq C_0|u|_b^{1-\lambda}$ on $D$; by
\re{ftcn}, $|u|_c\leq C\|\pd^cu\|_0$
and
\re{uccu} follows.
\end{proof}

\pr{hoin} Let $D\subset\rr^m$  have \ps\
and denote $\|\cdot\|_{D,a}$ by $\|\cdot\|_a$. Let $a,b,a_j,b_j$ be nonnegative real
numbers.
\bppp
\item
$
\|u\|_{\lambda a+(1-\lambda)b}\leq
C_{a,b}\|u\|_{a}^\lambda\|u\|_{b}^{1-\lambda}$ for $0<\lambda<1$.
\item
$\|uv\|_{a}\leq C_a(\|u\|_{0}\|v\|_{a}+\|u\|_{a}\|v\|_{0})$ .
\item
Suppose that $D_j\subset\rr^{n_j}$ has \ps\ for $j=1,\ldots, k$.
Let  $a_j,c_j$ be non-negative real
numbers,
and let $(b_1,\ldots, b_k)$ be in the convex hull of $(a_1,\ldots, a_k)$
and $(c_1,\ldots, c_k)$. Then
$$
\prod_{j=1}^k\|u_j\|_{D_j,b_j}\leq C_{k+|a|+|b|+|c|}
\Bigl(\prod_{j=1}^k\|u\|_{D_j,a_j}+\prod_{j=1}^k\|u\|_{D_j,c_j}\Bigr).
$$
\item Let $f$ be a map from $D$ into $D'\subset\rr^{n}$.
Assume that $D'$ has \ps.     Then
\gan
\|u\circ f\|_{a}\leq C_a(\|u'\|_{D',a-1}\|f'\|_{0}^a
 +\|u'\|_{D',0}\|f'\|_{a-1})+\|u\|_{D',0},\quad a\geq1,\\
\|u\circ f\|_{a}\leq C\min(\|u'\|_{D',0}|f|_{a},
\|u\|_{D',a}\|f'\|_{0}^a)+\|u\|_{D',0},\quad 0\leq a\leq 1.
\end{gather*}
\item  Let $f=I+\tilde f$ be a $C^a$ map from $ D$ into $\rr^m$.
There exists $C_0>1$
such that if  $|\tilde f'|<C_0^{-1}$ then
\gan
 \|f^{-1}-I\|_{f(D),a}\leq C_a\|\tilde f\|_{a}, \quad a\geq0.
\end{gather*}
\item   Let $f=I+\tilde f$ be a $C^1$ map  from $D$ into $ D'\subset\rr^m$ with $|
f'|<C_0$.
 Assume that
 $D\cup D'$ is contained in a   {\bf convex}  domain
  $D''$ of the cone property. Then
\gan
\|u\circ f-u\|_{a}\leq C_a(\|u'\|_{D'',a}\|\tilde f\|_{0}+
\|u'\|_{D'',0}\|\tilde f\|_{a}), \quad a\geq0.
\end{gather*}
\eppp
\end{prop}
\pf (i)-(ii) are proved in~\ci{Hosesi}.
The proof for (iii) is for $D_j$ being the same domain.
In fact, by (i), one has $\log\|u_j\|_{D_j,b_j}\leq
\lambda\log\|u_j\|_{D_j,a_j}+(1-\lambda)
\log\|u_j\|_{D_j,c_j}+\log C$ for all $j$.
Sum  over $j=1,\dots, k$ and take exponential on both
sides. The convexity of $e^x$ yields (iii).
For   $0\leq a\leq1$,  two inequalities in (iv) are verified directly.
Assume that (iv) holds when $a$ is replaced by $a-1$.
Assume that $a>1$. Then
$
(u\circ f)'= u'(f)f'.
$
For both cases $0< a-1\leq 1$ and $a-1>1$ we obtain
\aln
\|(u\circ f)'\|_{a-1}&\leq C(\|u'(f)\|_{a-1}\|f'\|_0+\|u'(f)\|_0\|f'\|_{a-1})\\
&\leq
C'(\|u'\|_{D',a-1}\|f'\|_0^{a}
+\|u'\|_0\|f'\|_{a-1}),
\end{align*}
where $\|u'\|_{a-1}\equiv \|u'\|_{D',a-1}$. This gives us
(iv) as $\|u\circ f\|_a\leq\|u \|_0+\|(u\circ f)'\|_{a-1}$.

\

(v).  It follows immediately from \rl{c2ps}
when $0\leq a\leq1$. We now prove  it  for $a>1$ by
using
a variant of
counting scheme in section~\ref{sec2}.
Define
\eq{shar}
\widehat\pd ^{1+k}v=\sum_{j_1+\cdots+j_l\leq k} p( \tilde f')\pd^{J_1}v\cdots
\pd^{J_l}v, \quad j_i=|J_i|-1\geq0,
\end{equation}
where $p(\tilde f')$ is a   polynomial in $(I+\tilde f')^{-1} $
and $\tilde f'$ and it might be different when it reoccurs.
Let $g=f^{-1}$ and $\tilde g=f^{-1}-I$. We have
$\tilde g'=-(({\mathbf 1}+\tilde f')^{-1}\tilde f'
)
\circ g$, i.e.
$\pd^I\tilde g=(\widehat\pd ^1\tilde f)\circ g$ for $|I|=1$.
Inductively,
\eq{shar+}
\pd^{J}\tilde g=(\widehat\pd ^{|J|} \tilde f)\circ g,\quad
|J|\geq1.
\end{equation}
By \rl{c2ps}, we have
$\yt|y_1-y_0|\leq|g(y_1)-g(y_0)|\leq 2|y_1-y_0|$.
Let $k=[a]-1$ and $\alpha=a-k-1$. Thus $\|g\|_{D',1+k+\alpha}\leq C\|f\|_{0}+
C\sum_{j\leq k}\|\widehat\pd ^{1+j}\tilde f\|_{\alpha}$. Now
\aln
 \|\widehat\pd ^{1+k}\tilde f\|_{\alpha}
&\leq\sum_{j_1+\ldots+j_l\leq k} \ C'\Bigl\{\|\tilde f\|_{ 1+\alpha}\|\tilde f\|_{1+j_1}\cdots
\|\tilde f\|_{1+j_l}\\
&\qquad + \sum_{1\leq i\leq l}\|\tilde f\|_{1+j_1}\cdots\|\tilde f\|_{1+j_i+\alpha}
\cdots\|\tilde f\|_{1+j_l}\Bigr\}.
\end{align*}
  By (iii) we obtain $\|\tilde g\|_{D',1+k+\alpha}
\leq \sum C(\|\tilde f\|_1^{l}+\|\tilde f\|_1^{l-1})\|\tilde f\|_{1+k+\alpha}\leq
C'\|\tilde f\|_{a}$.

 \

(vi).
By convexity of $D''$, we have
$ 
u(x+\tilde f(x))-u(x)=\tilde f(x)\cdot\int_0^1u'(x+t\tilde f(x))dt.
$ 
Consider case $0\leq a<1$.
We have
\aln
\Bigl|\int_0^1\{u'(y+t\tilde f(y))-u'(x+t\tilde f(x))\}dt\Bigr|&\leq
\max_t\|u'\|_{D'',a}|y-x+t(f(y)-f(x))|^a\\ &
\leq C\|u'\|_{D'',a}(1+\|\tilde f\|_{D,1})|y-x|^\alpha.
\end{align*}
 By $\|\tilde f\|_{1}\leq C$,
one sees easily that
$$
\|u(I+\tilde f)-u\|_{a}\leq C(\|u'\|_{D'',a}\|\tilde f\|_{0}+
\|u'\|_{D'',0}\|\tilde f\|_{a}),
$$
which is (vi) for $0\leq a<1$. Assume that the above inequality  holds
when $a$ is replaced by $a-1$. Assume that $a\geq1$. We  need
to estimate the $\|\cdot\|_{a-1}$ norm of
\begin{equation}\label{dxiu}
\pd_{x^i}(u\circ(I+\tilde f)-u)=(\pd_{y^i}u)\circ(I+\tilde f)-\pd_{x^i}  u
+\sum_{1\leq j\leq m}(\pd_{y^j}u)\circ(I+\tilde f)\cdot\pd_{x^i}\tilde f^j.
\end{equation}
By the induction assumption, we have
$$\|(\pd_{y^i}u)\circ(I+\tilde f)-\pd_{x^i}  u\|_{a-1}\leq
C(\|\pd^2u\|_{D'',a-1}\|\tilde f\|_{0}+\|\pd^2 u\|_{D'',0}\|\tilde f\|_{a-1}).$$
We need to  put 
 $\|\pd^2u\|_{D'',l} \|\tilde f\|_{a-1-l}\leq \|u'\|_{D'',l+1}\|\tilde f\|_{a-1-l}$
 into the desired form. By (iii), we get
\eq{udaf}
\|u'\|_{D'',l+1} \|\tilde f\|_{a-1-l}\leq
C( \|\tilde f\|_{0}\|u'\|_{D'',a} +\|u'\|_{D'',0} \|\tilde f\|_{a}
  ).
  \end{equation}
We now treat the term in the sum  of \re{dxiu}.
Set $\|u'\|_{b}\equiv\|u'\|_{D'',b}$.
By (iv)  we have $\|(\pd_{y^j}u)(I+\tilde f)\|_{a-1} \leq C(\|u'\|_{ a-1}+
\|u'\|_{1}\|\tilde f\|_{a-1})+\|u'\|_0$. Thus
\aln
\|(\pd_{y^j}u)(I+\tilde f)\cdot\pd_{x^i}\tilde f^j\|_{a-1}&\leq C((\|u'\|_{ a-1}+
\|u'\|_{1} \|\tilde f\|_{a-1} )\|\tilde f\|_{1}+\|u'\|_{ 0}\|\tilde f\|_{a})
\\
&\leq C'(\|u'\|_{ a-1}\|\tilde f\|_{1}+
\|u'\|_{1}\|\tilde f\|_{a-1}+
\|u'\|_{ 0}\|\tilde f\|_{a}).
\end{align*}
We can put the first two  terms in  the desired form by \re{udaf}.
\end{proof}


Two inequalities in the next proposition   are  used in estimating
 $P$ and $Q$.
 \pr{vbhi} Let $D$ be a convex domain in $\rr^m$ satisfying
 $$
 B_{\rho/c_0}\subset D\subset B_{c_0\rho}, \quad 0< \rho\leq3.
 $$
 Let $a_*=0$ for $0\leq a\leq1$ and $a_*=a$ for $a>1$. Let
 $\|\cdot\|_{a}=\|\cdot\|_{\cC^a(D)}$. Then
\aln
\rho^{a_*}\Bigl\|\prod_{j=1}^mu_j\Bigr\|_{ a}&\leq C_{a,c_0}\sum_{j=1}^m\|u_j\|_{ a}\prod_{i\neq j}
\|u_i\|_{ 0},
\\
\rho^{e}\prod_{j=1}^m\|u\|_{ d_j+b_j} &\leq C_{a,b,c,c_0}
\Bigl(\prod_{j=1}^m\|u_j\|_{ d_j+a_j}+\prod_{j=1}^m\|u_j\|_{ d_j+c_j}\Bigr),
\end{align*}
where
 $e=(b_1+d_1-[d_1])+\cdots+(b_m+d_m-[d_m])$,
 $a_j,c_j, d_j$ are non-negative real
numbers,  and $(b_1,\ldots, b_m)$
 is in the convex hull of $(a_1,\ldots, a_m)$ and $(c_1,\ldots, c_m)$.
  \end{prop}
 \begin{proof}
 The first estimate    is trivial if $a\leq1$. So   assume that $a>1$.
  We know that $D$ is convex and $B_{\rho/{c_0}}\subset
 D\subset B_{c_0\rho}$. Thus $D$ has \ps\
 if $1\leq\rho\leq\rho_0$, where the cone of fixed size with vertex at $x\in\ov D$
 can be found from the convex hull of $x$ and $B_{\rho/2}$.
 Assume now that $0<\rho<1$. Consider the isotropic
 dilation
 $S_\rho(x)=\rho x$.
 Then, $D_*=S_\rho^{-1}D$ has
 \ps.  Since $0<\rho<1$, we have
 \eq{ctw}
\rho^{a}\|u\|_{ a}\leq
\|u\circ S_\rho\|_{D_*,a}\leq\|u\|_{ a}.
\end{equation}
By \rp{hoin} (i) we obtain
\aln
 \rho^{a}\Bigl\|\prod_{j=1}^mu_j\Bigr\|_{ a}&\leq
\Bigl\|\prod_{j=1}^mu_j\circ S_\rho\Bigr\|_{D_*,a}
\leq
C \sum_{j=1}^m  \| u_j\circ S_\rho \|_{D_*,a}\prod_{i\neq j}
\|u_i\circ S_\rho\|_{D_*,0}.
\end{align*}
Using \re{ctw}, we get the first inequality easily.
Let $l_j\leq[d_j]$ be any non-negative integers.
By \re{ctw} and \rp{hoin} (iii), we get
\aln
& \rho^{e}\prod_{j=1}^m\|\pd^{l_j}u_j\|_{ b_j+d_j-[d_j]}\leq
\prod_{j=1}^m\|(\pd^{l_j}u_j)\circ S_\rho\|_{D_*,b_j+d_j-[d_j]}  \\
&\hspace{8ex} \leq   C\Bigl(\prod_{j=1}^m\|(\pd^{l_j}u_j)\circ S_\rho\|_{D_*,a_j+d_j-[d_j]}+
\prod_{j=1}^m\|(\pd^{l_j}u_j)\circ S_\rho\|_{D_*,c_j+d_j-[d_j]}\Bigr) \\
&\hspace{8ex} \leq C \Bigl(\prod_{j=1}^m\| u_j \|_{ l_j+a_j+d_j-[d_j]}+
\prod_{j=1}^m\| u_j \|_{ l_j+c_j+d_j-[d_j]}\Bigr)\\
&\hspace{8ex} \leq C  \Bigl(\prod_{j=1}^m\| u_j \|_{ d_j+a_j }+
\prod_{j=1}^m\| u_j \|_{ d_j+c_j }\Bigr).
\end{align*}
Summing   over all non-negative integers $l_j\leq [d_j]$ gives us the second inequality.
\end{proof}

One can also obtain other inequalities via  dilation. The inequalities below
are not directly used in this paper.
\pr{hoin+}   Let $a_*$ be as in \rp{vbhi} .
Let $\rho, D, \|\cdot\|_{a}$ be as in \rp{vbhi}.
 Let $a,b,a_j,b_j$ be nonnegative real
numbers.
\bppp
\item
$ \rho^{c_*}\|u\|_{\lambda a+(1-\lambda)b}\leq C_{a,b}
\|u\|_{a}^\lambda\|u\|_{b}^{1-\lambda}$ for $0<\lambda<1$, where
$c_* =\lambda a+(1-\lambda)b$.
 \item Let $f$ be a map from $D $ into $D'\subset\rr^{n}$.
Assume that $D'$  is convex and $B_{c_0^{-1}\rho}\subset D'\subset B_{c_0\rho}$.
Then
\gan
\rho^{a_*}\|u\circ f\|_{a}\leq C_{a,c_0}\rho(\|u'\|_{D',a-1}\|f'\|_{0}^a
 + \|u'\|_{D',0}\|f'\|_{a-1})+\|u\|_{D',0},\quad a\geq1,\\
 \|u\circ f\|_{a}\leq C_{c_0}\min(\|u\|_{D',1}\|f\|_{a},
\|u\|_{D',a}\|f\|_{1}^a)+\|u\|_{D',0},\quad 0\leq a\leq 1.
\end{gather*}
\item  Let $a_{**}=0$ for $0\leq a\leq2$ and $a_{**}=a-1$ for $a>2$.
 Let $f$ be a $\cC^a$ map from $ D$ into $\rr^m$.
There exists $C_0>1$
such that if  $|  f'-I|<1/C_0$ then
\gan
 \rho^{a_{**}} \|f^{-1}-I\|_{f(D ),a}\leq C_a \|  f-I\|_{a}.
  \end{gather*}
\item Let $f=I+\tilde f$ be a $\cC^1$ map  from $D$ into $ D'\subset\rr^m$ with $|f'
|<C_0$.
 Assume that
 $D\cup D'$ is contained in a   {\bf convex}  domain
  $D''$  satisfying $B_{\rho/c_0}\subset D''\subset B_{c_0\rho}$. Then
\gan
\rho^{a_*}\|u\circ f-u\|_{a}\leq C_{a,c_0}(\|u'\|_{D'',a}\|\tilde f\|_{0}+
\|u'\|_{D'',0}\|\tilde f\|_{a}).
\end{gather*}
\eppp
\end{prop}
\begin{proof}  We may assume that $0<\rho<1$.
Let $u_*(x)=u(\rho x)$,  $ f^*(x)=\rho^{-1}f(\rho x)=x+\tilde f^*(x)$,
and $D_*=\rho^{-1}D$.
 Then $D_*$ has 
 \ps.

(i) is immediate, by applying \rp{hoin} to $u_*$ and by using \re{ctw}.

(ii). The case $0\leq a\leq1$ is verified directly.
Assume that $a>1$. 
Applying \rp{hoin} to $u_*\circ f^*$, we get
\gan
\|u_*\circ f^*\|_{D_*,a}\leq C_{a}(\|\pd u_*\|_{D'_*,a-1}\|f^{*\prime}\|_{D_*,0}^a
 +\|\pd u_*\|_{D'_*,0}\|f^{*\prime}\|_{D_*,a-1})+\|u_*\|_{D'_*,0}.
\end{gather*}
Since $f^{*\prime}(x)=f'(\rho x)$, then
$\|f^{*\prime}\|_{D_*,b}\leq \|f'\|_{D,b}
$  for $b\geq0$.
By \re{ctw}, we also have
\gan
\|u_*\circ f^*\|_{D_*,a}=\|(u\circ f)\circ S_\rho\|_{D_*,a}\geq\rho^a\|u\circ f\|_{ a},
\\ \|\pd u_*\|_{D_*',a-1}=  \|\rho (\pd u)_*\|_{D_*',a-1}\leq \rho \|\pd u\|_{a-1},
 \quad a\geq1.
 \end{gather*}
Simplifying gives us (ii).

 (iv). Let $\tilde f=f-I$. The case $0\leq a\leq1$ is verified directly. Assume that $a>1$.
Then
$$
\|u_*\circ f^*-u_*\|_{D_*,a}\leq C
(\|\tilde f^*\|_{a}\|\pd u_*\|_{D''_*,0}+\|\tilde f^*\|_{0}\|\pd u_*\|_{D''_*,a}).
$$
As in (ii),  we can  get (iv) by \re{ctw} and 
$$
 \|\tilde
 f^*\|_{D_*',a}=\rho^{-1} \|\tilde f\circ S_\rho\|_{D_*',a}\leq  \rho^{-1}\|\tilde f\|_{a}.
 $$

(iii). (A dilation would give us $a_{**}=a$.)
Let  $ g=f^{-1}=I+\tilde g$. We have $\tilde g=-\tilde f\circ g$
and $\tilde g'=-\{\tilde f'(I+\tilde f')^{-1}\}\circ g$.
By \rl{c2ps},  $C^{-1}|x'-x|\leq |g(x')-g(x)|\leq C|x'-x|$. We get immediately
$$
\|\tilde g\|_{f(D),a}\leq C\|\tilde f\|_a, \quad 0\leq a\leq 2.
$$
Assume that $k\geq1$. Let $\pd^j$ denote a derivative
of order $j$.  Recall that by \re{shar}-\re{shar+},
$$
\pd^{K}\tilde g=\sum_{
j_1+\cdots+j_l\leq k } \{p( \tilde f')\pd^{J_1}\tilde f\cdots
\pd^{J_l}\tilde f\}\circ g, \quad k=|K|-1\geq0, \  j_i=|J|_i-1\geq0.
$$
Set $1+k=[a]$ and $\alpha=a-k-1$.
Computing the H\"older  ratio of \re{shar+} and applying \rp{vbhi}    with $d_j=1$, we obtain
(iii) from
\begin{equation*}
\|\pd^{K}\tilde g\|_\alpha\leq C\sum_{
j_1+\cdots+j_l\leq k }\rho^{-j_1-\cdots-j_l-\alpha}\|\tilde f\|_{1+k+\alpha}
\leq C'\rho^{-k-\alpha}\|\tilde f\|_{1+k+\alpha}. \qedhere
\end{equation*}
\end{proof}

%
%
%

\begin{rem} If all norms in \rp{hoin+} are replaced by the scalar-invariant norm $\|\cdot\|^*$,
defined as
$$
\|u\|_{D,a}^*=\|u\circ S_d\|_{S_d^{-1}D,\alpha}
$$
with $d$ being the diameter of $D$,  we can take
 $c_*=a_*=a_{**}=0$  for the proposition.
 For the use of scalar-invariant norms to derive estimates on the
 Bochner-Martinelli-Koppelman formula for balls in $\cc^n$, see \ci{Weeini}.
  \end{rem}

\setcounter{thm}{0}\setcounter{equation}{0}\section{The Henkin homotopy formula
}
\label{residue}
Recall notation $z'=(z_1,\ldots, z_{n-1})$, $z=(z',z^n)$ and $x=\pi(z)=(\RE z,\IM z')$.

In this appendix, we will derive the following version of
Henkin's homotopy formula.
\th{hhf} Let $M\subset\cc^n$ be a graph  $y^n=|z'|^2+\hat r(x)$ over $D\subset\rr^{2n-1}$. Assume
that $0<\rho<\rho_0\leq3$ and
 \ga
\label{bhrooap}\ov{D_{\rho_0}}\subset D,
\quad
 \hat r(0)=0,\quad \pd \hat r(0)=0,\quad \|\hat r\|_{ \rho_0,2}=\|\hat r\|_{ D_{\rho_0},2}<1/C_0
\end{gather}
with $C_0$ sufficiently
large.  Assume that $0<\rho<\rho_0\leq3$. Let $\varphi$ be a continuous
 tangential $(0,q)$-form on $\ov M_\rho$. Assume
that $\dbm\varphi$ is continuous as currents on $M_\rho$
 and admits a continuous extension on $\ov M_\rho$.
 If $0<q<n-2$, then on $M_\rho$ and as currents
 $$
 \varphi=\dbm (P_0+P_1)\varphi+(Q_0+Q_1)\dbm\varphi\mod\db r,
 $$
 where $P_0,P_1,Q_0,Q_1$ are defined by \rea{fhc0+}.
\end{thm}
Recall that $M_{\rho }=M\cap
\{z\colon |x^n|^2+y^n<\rho ^2\}$ and $D_{\rho }=\{x\in D
\colon |x^n|^2+y^n<\rho^2\}$.
When $M$ is strictly convex,
see Henkin~\ci{Hesese} for the proof. Our proof will follow~\ci{Weeinia}
via  Stokes' theorem.


Set $r=-y^n+|z'|^2+\hat r(x)$ and define
\ga\label{levip}
F(\zeta,z)=r_z\cdot(\zeta-z)+\yt\sum_{1\leq j,k\leq n}r_{z^jz^k}
(\zeta^j-z^j)(\zeta^k-z^k),\\
S_t(z)=\pi\{\zeta\in M\colon |F(\zeta,z)|=t\}, \quad S'_t(\zeta)=
\pi\{z\in M\colon |F(\zeta,z)|=t\}.\nonumber
\end{gather}
Note that $r_{z_jz_k}=\hat r_{z_jz_k}$.

\le{ap0} Let $n\geq2$. Let $M$ be as in \rt{hhf}.
There exists $C_1>1$ satisfying the following.
\bppp\item
If $0<t<(\rho_0-\rho)^2/C_1$ and $z \in M_\rho$,
then $S_t(z)$ and $S'_t(z)$ are
compact subsets of $D_{\rho_0}$.
\item
Assume that $\hat r\in \cC^3(\ov D_{\rho_0})$, $z\in M_\rho$ and
$0<t<C_1^{-1}\min\{(\rho_0-\rho)^2,
(1+ \|\hat r\|_{{\rho_0},3})^{-1}\}$. Then
$S_t(z)$ and $S'_t(z)
$
are  smooth and of classes $\cC^3$ and $
\cC^1$, respectively.\eppp
\end{lemma}
\begin{proof} Set $\epsilon=\|\hat r\|_{ \rho_0,2}$.
(i)
By \rl{mdb}, $D_{\rho_0}$ is convex. On $M_{\rho_0}\times M_{\rho_0}$,
 $|r_z\cdot(\zeta-z)|\geq |\zeta-z|^2/C$
 by \rl{distb} and $|\hat r_{z^jz^k}|\leq C\epsilon$. Hence,    (i) follows from
 \eq{Fzzg}
 |F(\zeta,z)|\geq |\zeta-z|^2/{(2C)}.
 \end{equation}

 (ii)
Fix  $z\in M_{\rho_0}$. Let $\zs^n=-2ir_z(\zeta-z)$.  By \rl{m2h},
 $\xis=(\RE\zeta_*,\IM\zeta'_*)$ are coordinates of $M_{\rho_0}$ and
\al
\label{anxi}
\xi_*^n&=\RE(-2ir_z\cdot(\zeta-z))=\xi^n-x^n+2\IM(\ov{z'}\cdot\zeta')+\IM\{2\hat r_z\cdot(\zeta-z)\},\\
\eta^n_*&=\IM(-2ir_z\cdot(\zeta-z))=|\zs'|^2+
\sum_{|I|= 2}\trd_*^2\hat r(\xi,x)\xis^I,
\label{anxi+}\end{align}
where $\trd_*^i$ is defined by \re{cont}.
Write $u=\RE(-2i F(\zeta,z))$ and $v=\IM(
-2iF(\zeta,z))$.  We have
\al\nonumber
-2i F(\zeta,z)&=-2ir_z\cdot(\zeta-z)-i\sum_{1\leq j,k\leq n}r_{z^jz^k}
(\zeta^j-z^j)(\zeta^k-z^k),\\
\label{anu}
u(\xis)&=\xi_*^n+\tilde u(\xis), \quad \tilde u(\xis)=
\sum_{|I|=2}a_I\circ\Psi^{-1}(\xis,x)\xis^I,\\
v(\xis)&= |\zs'|^2+\tilde v(\xis),\quad \tilde v(\xis)=
\label{anu+}\sum_{|I|=2}b_I\circ\Psi^{-1}(\xis,x)\xis^I,
\end{align}
where $|  a_I |+| b_I |\leq C\|\hat
r\|_{\rho_0,2}\leq C'\epsilon$
and  $|\pd^1 a_I(\zeta,z)|+|\pd^1b_I(\zeta,z)|\leq C\|\hat
r\|_{\rho_0,3}$. Suppose
that
$
|F(\zeta,z)|=t$
and
\eq{htsq}
t^{1/2} \|\hat r\|_{\rho_0,3}<1/{C_0}
\end{equation}
and $C_0$ is sufficiently large.
To show that $S_t(z)$ is smooth, we
need to verify that $d_{\zs',\xi_*^n}(u^2+v^2)\neq0$ when $u^2+v^2=t^2$.
By \re{Fzzg}-\re{anxi} we know that
$$|\zeta-z|/C\leq |\xis|\leq C|\zeta-z|\leq C't^{1/2}.
$$
By \re{Fzzg} and \re{htsq}, we obtain $|\zeta-z|\|\hat r\|_{\rho_0,3}\leq 1/C_0 $ and
$$
|\pd^1\tilde u|+|\pd^1\tilde v|\leq C(\|\hat r\|_{\rho_0,3}|\zeta-z|^2
+\epsilon|\zeta-z|)\leq C'(C_0^{-1}+\epsilon)|\xis|.
$$
Assume that $d_{\zs',\xi_*^n}(u^2+v^2)=0$.
 Then
$$
(\xi_*^n+\tilde u)(1+\tilde u_{\xi_*^n})+ vv_{  \xi_*^n} =0, \quad
u\tilde u_{\zs'}+v (\ov{\zs'}+\tilde v_{\zs'})=0.
$$
From the first identity we get
$|\xi_*^n|\leq C(|\tilde u|+|v|)\leq C'(|\zs'|^2+|\xi_*^n|^2)$.
Therefore, for $t<1/C$,
$|\xi_*^n|\leq C|\zs'|^2$. Now, $|v|\geq|\zs'|^2-C\epsilon(|\zs'|^2+|\xi_*^n|^2)\geq
\yt|\zs'|^2$ and $|u|\leq C|\zs'|^2$. By \rl{psii}, $|\pd^1_{\xi_* } \xi |\leq C$,
and by  \re{anu+} and \re{htsq}, 
\gan
|\tilde v_{\zs'}|  \leq |\zs'|/2,
\quad
|\ov{\zs'}+\tilde v_{\zs'}| \geq  |\zs'|/2,\\
 |\zs'|^3/4\leq|v(\ov{\zs'}+\tilde v_{\zs'})|=|u\tilde u_{\zs'}|\leq
C|\zs'|^2(\epsilon +C_0^{-1})|\zs'|.
\end{gather*}
Hence $\zs'=\xi_*^n=0$, when $C_0^{-1},\epsilon$
are sufficiently small. This shows that $t=0$, a contradiction.
It is clear that for a fixed $z$,
$u^2+v^2$ is a function
of class $\cC^3$ in $(\zeta',\xi^n)$. This shows that $S_t(z)$ is smooth
and of class $\cC^3$.

For a fixed $\zeta\in M_\rho$,   $ (\zs',\xi_*^n)$ still
form a   coordinate system of class $\cC^2$ for $M$. We have the
same formulae \re{anxi}-\re{anxi+}, only to vary $z\in M_{\rho_0}$; see
\rrem{62}. Thus   \re{anu}-\re{anu+}
 are still valid, where $\xi$ is fixed and $z$ varies in $ M_{\rho_0}$.
The same argument
shows that   $S'_t(\zeta)$ is smooth and of class $\cC^1$.
\end{proof}

Let us recall from section~\ref{sec3}
$$
\Omega^{+-}_{0,q}=\f{(r_\zeta-r_z)\cdot d\zeta \wedge
r_z\cdot d\zeta\wedge(\db_\zeta \pd_\zeta r)^{n-2-q}
\wedge(\db_z r_z\wedge d\zeta)^q}
{(r_\zeta\cdot(\zeta-z))^{n-1-q}(r_z\cdot(\zeta-z))^{q+1}}.
$$
The following computation
 is essential in Folland-Stein~\ci{FSsefo}. See also Romero~\ci{Ronion}.
\le{ap1} Let $n\geq2$ and $0\leq q\leq n-1$. Let $M$  satisfy \rea{bhrooap}
and $\hat r\in \cC^3$. Let $F(\zeta,z)$ be the Levi polynomial of $r$
about $z\in M$.
Let
$\varphi$ be a continuous  tangential $(0,q)$-form  on $M$. Then
 \eq{resi}
\lim_{t\to0}\int_{ |F(\zeta,z)|=t,\zeta\in M}\varphi(\zeta)\wedge\Omega_{0,q}^{+-}(\zeta,z)
= \f{(2\pi i)^n}{2} \varphi(z) \mod{\db r(z)},
\end{equation}
  where
 the convergence is uniformly in $z$ on each compact subset of $M$.
\end{lemma}


Here is an outline.
Following~\ci{FSsefo},
we  will use the Levi polynomial $F(\zeta,z)$ to define
new coordinates of $M$ near $z$ and compute each term in the kernel.
On  $T^{1,0}_z M$ there are two
  quadratic forms  $h=\sum_{1\leq j,k\leq n}r_{z^j\ov{z^k}}t^j\ov t^k$,
$A=\sum_{1\leq j,k\leq n}r_{z^j z^k}t^jt^k$.
We will express the kernel   in    $h$ and $A$, with error terms.
When we compute the residue,  the error terms can be removed
 via   non-isotropic dilation. This gives us a limit kernel  (see~\re{keyob})
 on a non-isotropic sphere. Roughly speaking,
  the limit kernel is expressed in   $h$ and
  $A$, but not in   $\ov A$. For latter purpose
 we will use   $-2i  F(\zeta,z)$, instead of its linear part
 $-2ir_z\cdot(\zeta-z)$, as part of coordinates.
Without $\ov A$,
  we  remove $A$ in the limit kernel by averaging. By a     linear transformation, we reduce
$h$ to the identity, and   compute the residue.

\medskip

CHANGE NOTATION. Let $\zs'=\zeta'-z'$ and
\eq{-2if}
\zs^n=-2i F(\zeta,z)=-2ir_z\cdot(\zeta-z)-i\sum_{1\leq j,k\leq n}\hat r_{z^jz^k}
(\zeta^j-z^j)(\zeta^k-z^k).
\end{equation}
Then  near $z\in M_{\rho_0}$, $\zeta_*',\xi_*^n$ form
coordinates of $M_{\rho_0}$. We will modify  some earlier computations where
the approximate Heisenberg transformation is used.

Applying a remainder formula on convex domain $D_{\rho_0}\times\rr$ and then
letting $\zeta,z\in M_{\rho_0}$,
 we obtain 
\gan
r(\zeta)-r(z)=2\RE F(\zeta,z)
+\sum_{1\leq j,k\leq n}r_{z^j\ov{z^k}}
(\zeta^j-z^j)(\ov{\zeta^k - z^k})+E(\zeta,z),\\
E(\zeta,z)=o(|(\zeta'-z',\xi^n-x^n)|^2),\quad
|\pd^1E(\zeta,z)|\leq C\|\hat r\|_{\rho_0,3}|(\zeta'-z',\xi^n-x^n)|^2.
\end{gather*}
With $z\in M_{\rho_0}$ the equation $r(\zeta)=0$ becomes
$$
-\eta_*^n+\sum_{1\leq j,k\leq n}r_{z^j\ov{z^k}}(\zeta^j-z^j)(\ov{\zeta^k- z^k})
 +E(\zeta,z)=0,
$$
where $\zeta_n-z_n$, with $\zeta'-z'=\zs$, is
a local solution to \re{-2if}. Replace   $\zeta-z$ by the solution
expressed in  $\zs$.  Solving
for $\eta_*^n$  from the new equation shows that for $(\zs',\xi_*^n)\in \pi\psi_z(M_{\rho_0})
$
\eq{etss}
\eta_*^n=\sum_{1\leq\alpha,\beta<n}h_{\alpha\ov \beta}\zs^\alpha\ov\zs^\beta+
o(|\zs'|^2)+O(|\xi_*^n|
|\zs'|+|\xi_*^n|^2).
\end{equation}

Write
$A(\zs',\xi_*^n)=o(k)$ if $ t^{-k}A(t\zs',t^2\xi_*^n)$ tends to zero uniformly for $|\zs'|+|\xi_*^n|<1$ and
 small $|t|$, where  $A $ is a  function or differential form.
 Write
$A(\zs',\xi_*^n)=O(k)$
if $|t^{-k}A(t\zs',t^2\xi_*^n)|\leq C$.
Thus $ \zeta'-z'=O(1)$,  $\zs^n=O(2)$,
$d\xi_*'=O(1)$ and $d\xi_*^n=O(2)$.
By \re{-2if}-\re{etss}, $\zeta^n-z^n= -\frac{r_{z'}}{r_{z^n}}\cdot\zs'
+O( 2)$.
The Levi matrix $(h_{\alpha\ov\beta})$ at $z$ is determined by
\eq{derh}
\sum_{1\leq j,k\leq n}r_{z^j\ov{z^k}}(\zeta^j-z^j)(\ov{\zeta^k - z^k})=
\sum_{1\leq\alpha,\beta<n}h_{\alpha\ov \beta}\zs^\alpha\ov\zs^\beta+
o( 2).
\end{equation}
Explicitly, we have
\al\label{forh}
 h_{\alpha\ov \beta} =
 r_{z^\alpha\ov {z^\beta}} -r_{z^\alpha \ov{z^n}}\frac{r_{ \ov {z^\beta}}}{r_{\ov{z^n}}}
 -\frac{r_{z^\alpha}}{r_{z^n}} r_{z^n\ov {z^\beta}}+r_{z^n\ov{z^n}}
 \frac{r_{z^\alpha}}{r_{z^n}}
 \frac{r_{\ov {z^\beta}}}{r_{\ov{z^n}}}.
\end{align}
We also have the holomorphic quadratic form
\eq{derh+}
\sum_{1\leq j,k\leq n}r_{z^j  z^k}(\zeta^j-z^j)( \zeta^k-  z^k)=
\sum_{1\leq\alpha,\beta<n}A_{\alpha \beta}\zs^\alpha \zs^\beta+ o(2).
\end{equation}
Set
 $h=(h_{\alpha\ov\beta})$ and $A=(A_{\alpha\beta})$.   To simplify the notation, we define
$$
h(\zs',\ov{\zs'})=\sum_{1\leq\alpha,\beta<n}h_{\alpha\ov\beta}\zs ^\alpha\ov{\zs}^\beta,
\quad A(\zs',\zs')=\sum_{1\leq\alpha,\beta<n}A_{\alpha\beta}\zs^\alpha\zs ^\beta.
$$
We can write
$$
\ov\pd_{\zs}h(\zs',\ov{\zs'})=h(\zs',d\ov{\zs'}), \quad
\pd_{\zs}A(\zs',\zs')=2A(d\zs',\zs').
$$
We will need to
express the kernel in
 $h=(h_{\alpha\ov\beta})$ and $A=(A_{\alpha\beta})$,
  but not in $(\ov A_{\alpha\beta})$.
  This will play a crucial  role,  when we eliminate $A$ later via   averaging.

\

For the denominator of the kernel, using \re{-2if} we get
\al\label{2irz}
-2ir_z\cdot(\zeta-z)&=\zs^n+i\sum_{1\leq j,k\leq n}r_{z^j  z^k}(\zeta^j-z^j)( \zeta^k-  z^k)\\
& =
\xi_*^n+i (h(\zs',\ov{\zs'})+A(\zs',\zs'))+ o( 2).\nonumber
\end{align}
Using $r_\zeta\cdot(\zeta-z)= r_z\cdot(\zeta-z)+(r_{\zeta}-r_z)\cdot(\zeta-z)$, we get
\aln
&-2ir_\zeta\cdot(\zeta-z)
=-2ir_z\cdot(\zeta-z)-2i\sum_{1\leq j,k\leq n}r_{z^jz^k}(\zeta^j-z^j)( \zeta^k-  z^k)
\\& \qquad \qquad -2i\sum_{1\leq j,k\leq n}r_{z^j\ov {z^k}}(\zeta^j-z^j)
(\ov{\zeta^k - z^k})+o(2)
\\
&\qquad \quad=-2ir_z\cdot(\zeta-z)-2i(h(\zs',\ov{\zs'})+A(\zs',\zs'))+
o( 2)\quad \text{(by \re{derh}-\re{derh+})}\\
&\qquad \quad=\xi_*^n-i (h(\zs',\ov{\zs'})+A(\zs',\zs'))+ o(2).
\quad \quad \text{(by \re{2irz})}
\end{align*}
We arrive at the basic relations
\al
-2ir_\zeta\cdot(\zeta-z)&=\xi_*^n-i (h(\zs',\ov{\zs'})+A(\zs',\zs'))+ o(2),\label{2irzz}\\
-2ir_z\cdot(\zeta-z)&=\xi_*^n+
i (h(\zs',\ov{\zs'})+A(\zs',\zs'))+ o(2).\label{2irzz+}
\end{align}


We now computer the numerator of the kernel. Using the first-order expansion of $r_\zeta$
 about $z$, we get
\aln
&(r_\zeta-r_z)\cdot d\zeta=\sum_{1\leq j,k\leq n}\bigl\{r_{z^jz^k}(\zeta^k-z^k)
+ r_{z^j\ov{z^k}}(\ov{\zeta^k - z^k})+O(2) 
\bigr\}d \zeta^j
\\
&\qquad \quad=\pd_\zeta \sum_{1\leq j,k\leq n}\Bigl\{\yt r_{z^jz^k}( \zeta^j-z^j)(\zeta^k-z^k)
+ r_{z^j\ov{z^k}}( \zeta^j-z^j)(\ov{\zeta^k- z^k})\Bigr\}+o(2). 
\end{align*}
Note that for fixed $z$, $\zeta\to\zs$ is holomorphic. So we can switch the above
$\pd_\zeta$ to
$\pd_{\zs}$ and then restrict  it to $M$, which gives us
\eq{}(r_\zeta-r_z)\cdot d\zeta= h(d\zs',\ov{\zs'})+
A(d\zs',\zs')+o(2).
\end{equation}
Recall that
$
\zs^n=-2ir_z\cdot(\zeta-z)-i\sum_{1\leq j,k\leq n}r_{z^jz^k}
(\zeta^j-z^j)(\zeta^k-z^k).
$
Applying $\pd_\zeta$  gives us
\aln
-2ir_z\cdot d\zeta&=\pd_{\zeta}\Bigl\{\zs^n+i\sum_{1\leq j,k\leq n}r_{z^jz^k}
(\zeta^j-z^j)(\zeta^k-z^k)\Bigr\}\\
&=\pd_{\zeta}\Bigl\{\zs^n
+ i\sum_{1\leq\alpha,\beta<n}A_{\alpha \beta}\zs^\alpha \zs^\beta+ o(2)\Bigr\}
\qquad\text{by \re{derh+}}\\
&=d\zs^n+2i
A(d\zs',\zs')+o(2).
\end{align*}
Recall that $z$ is fixed. So on $M$ we have
\aln
\ov\pd_\zeta \pd_\zeta r&=\ov\pd\pd r(\zeta)=\ov\pd\pd\Bigl\{\sum_{1\leq j,k\leq n}
r_{z^j\ov{z^k}}(\zeta^j-z^j)(\ov{\zeta^k- z^k})+o(2)\Bigr\}\\
&=\ov\pd\pd \Bigl\{\sum_{1\leq\alpha,\beta<n}
h_{\alpha\ov \beta}\zs^\alpha\ov\zs^\beta
 \Bigr\}+o(2)\qquad \text{by \re{derh}}
\\
&=- h(d\zs',d\ov{\zs'})+o(2).
\end{align*}
On $M$, 
$\ov\pd r_{z^j}=\ov\pd_M r_{z^j}\!\mod{\db r(z)}$, and
$
\ov\pd r_{z^j}=
\sum_{\beta=1}^{n-1}\bigl(r_{z^j\ov z^\beta}-\frac{r_{\ov z^\beta}}{r_{\ov z^n}}r_{z^j
\ov z^n }\bigr)d\ov z^\beta\mod{\db r(z)}.
$
Recall that
  $\zeta^n-z^n=-r_{z^n}^{-1}r_{z'}\cdot\zs'+O(2)$. Thus
 $d\zeta^n= -\frac{r_{z'}}{r_{z^n}}\cdot d\zs'
+O(2). 
$
Therefore,
\aln
(\ov\pd r_z)\wedge d\zeta&=\sum_{1\leq\alpha,\beta<n}
\bigl(r_{z^\alpha\ov z^\beta}-\frac{r_{\ov z^\beta}}{r_{\ov z^n }}r_{z^\alpha
\ov z^n }\bigr)d\ov z^\beta \wedge d\zs^\alpha+O(2)
\\
&\quad-\sum_{1\leq\alpha,\beta<n}
 \bigl(r_{z^n\ov z^\beta}- \frac{r_{\ov z^\beta}}{r_{\ov z^n }}r_{z^n
\ov z^n }\bigr)d\ov z^\beta\wedge \frac{r_{z^\alpha}}{r_{z^n}}\,
d\zs^\alpha\mod{\db r(z)}.
\end{align*}
Looking at \re{forh}, we see that
$$
(\pd r_z) {\wedge}d\zeta=-h(d\zs',d\ov{z'})
+O(2)\mod\db r(z).
$$
We apply the non-isotropy dilation $T_t(\zs',\xi_*^n)=(  t\zs',t^2\xi_*^n)$
and summarize the above
as
\aln
t^{-2}T_t^*\{-2ir_\zeta\cdot(\zeta-z)\}&=\xi_*^n-
i (h(\zs',\ov{\zs'})+A(\zs',\zs'))+o(t^0),\\
t^{-2}T_t^*\{-2ir_z\cdot(\zeta-z)\}&=\xi_*^n+
i (h(\zs',\ov{\zs'})+A(\zs',\zs'))+o(t^0),\\
t^{-2}T_t^*\{(r_\zeta-r_z)\cdot d\zeta\}&=h(d\zs',\ov{\zs'})+
A(d\zs',\zs')+o(t^0),
\\
t^{-2}T_t^*\{-2ir_z\cdot d\zeta\}&=
d\zs^n+2i
A(d\zs',\zs')+
o(t^0),
\\
t^{-2}T_t^*\{\ov\pd_\zeta r_\zeta\dot{\wedge} d\zeta\}&=
-h(d\zs',d\ov{\zs'}) +
 o(t^0),\\
t^{-1}T_t^*\{ \ov\pd_z r_z\dot{\wedge}d\zeta\}&=-h(d\zs',d\ov{z'})
+O(t)\mod{\db r(z)}.
\end{align*}

By \re{-2if} and \re{2irz},    the
 subset of $M$ defined by $2|F(\zeta,z)|=t^2$   has the form
$$S_t(z)= \left\{(\zs',\xi_*^n)\colon |\xi_*^n+ih(\zs',\ov{\zs'})+o(2)|=t^2\right\}.$$
Then we have a non-isotropic sphere
\gan
\lim_{t\to0}T_t^{-1}(S_t)=\left\{(\zs',\xi_*^n)\colon
|\xi_*^n+ih(\zs',\ov{\zs'})|=1\right\}\df S.
\end{gather*}
Set $N_*(\zs',\xi_*^n)=\xi_*^n +ih(\zs',\ov{\zs'})$.
As $t$ tends to $0$,  $t^{-q}T_t^*\Omega_{(0,q)}^{+-}\mod\db r(z)$ converges uniformly
in a neighborhood of $  S$  to   $\widetilde\Omega$ given by
\eq{keyob}
\frac{(h(d\zs',\ov{\zs'})+A(d\zs',\zs'))\wedge
(d\zs^n+2i A(d\zs',\zs'))\wedge h(d\zs',d\ov{\zs'})^{n-q-2}\wedge
h(d\zs',d\ov{z'})^q}{-(2i)^{1-n}\ov N_*^{n-q-1}N_*^{q+1}
{(1-i\ov N_*^{-1} \zs' A \zs'^T)}^{n-q-1}{(1+i N_*^{-1} \zs' A \zs'^T)}^{q+1}}.
\end{equation}
Here $N_*=N_*(\zs',\xis^n)$. Note that except at the origin, the above form is smooth  in the
 $(\zs',\xi_*^n)$-space.
Therefore, for $\varphi(\zeta)=\sum_{|I|=q}\varphi_{\ov I}(\zeta)d\ov{\zeta'}^I$
\aln
& \int_{2|F(\zeta,z)|=t^2}\varphi(\zeta)\wedge\Omega_{0,q}^{+-}(\zeta,z)
= \int_{(\zs',\xi_*^n)\in
S_t(z)}\varphi(\zeta(\zs',\xi_*^n))\wedge\Omega_{0,q}^{+-}(\zeta(\zs',\xi_*^n),z)\\
&\hspace{3em}  =
\int_{(\zs',\xi_*^n)\in T_t^{-1}S_t(z)}
T_t^*\left\{\varphi(\zeta(\zs',\xi_*^n))\wedge\Omega_{0,q}^{+-}(\zeta(\zs',\xi_*^n),z)\right\}\\
&\hspace{3em} =\sum_{|I|=q}\varphi_{\ov I}(z)
\int_{(\zs',\xi_*^n)\in T_t^{-1}S_t(z)}
T_t^* \left\{d \ov{\zs'}^I\wedge
\Omega_{0,q}^{+-}(\zeta(\zs',\xi_*^n),z)\right\}\mod\db r(z).
\end{align*}
 We obtain \aln &
\lim_{t\to0}
\int_{2|F(\zeta,z)|=t^2}\varphi(\zeta)\wedge\Omega_{0,q}^{+-}(\zeta,z)
=\sum_{|I|=q}\varphi_{\ov I}(z)
\int_{  S}
 d\ov{\zs'}^I\wedge
\widetilde\Omega (\zs',\xi_*^n,z)\mod\db r(z).
\end{align*}

Return to $\widetilde\Omega$, defined by \re{keyob}. We are ready  to remove $A$ via an averaging.
By \re{bhrooap},  $|A|$ is small. Express
$
 (1-iN_*^{-1} \zs' A\zs'^T)^{-(n-q-1)}(1+i\ov N_*^{-1} \zs' A\zs'^T)^{-(q+1)}
$
 as a convergent power series
in $N_*^{-1} \zs' A\zs'^T$ and $\ov N_*^{-1} \zs' A\zs'^T$.
Note that $N(\xis)$ is invariant under the rotation
$e_\theta\colon(\zs',\xi^n)\to(e^{i\theta}\zs',\xi_*^n)$.
Using $|I|=q$, we can verify that
\aln
 &\int_{(\zs',\xi_*^n)\in  S}d\ov{\zs'}^{I}\wedge\widetilde\Omega(\zs',\xi_*^n,z)=
\frac{1}{2\pi}\int_{\theta=0}^{2\pi}d\theta
\int_{(\zs',\xi_*^n)\in e_\theta   S}d\ov{\zs'}^{I}\wedge\widetilde\Omega(\zs',\xi_*^n,z) \\
&\hspace{2em}=
\frac{1}{2\pi}\int_{\theta=0}^{2\pi}d\theta
\int_{(\zs',\xi_*^n)\in  S}e_\theta^*\{d\ov{\zs'}^{I}\wedge\widetilde\Omega \}
=
\frac{1}{2\pi}
\int_{(\zs',\xi_*^n)\in  S}\int_{\theta=0}^{2\pi}
e_\theta^*\{d\ov{\zs'}^{I}\wedge\widetilde\Omega \}d\theta.
\end{align*}
Therefore
$\int_{(\zs',\xi_*^n)\in  S}d\ov{\zs'}^{I}\wedge\widetilde\Omega(\zs',\xi_*^n,z)=
\int_{(\zs',\xi_*^n)\in  S}d\ov{\zs'}^{I}\wedge \Omega'(\zs',\xi_*^n,z),$
where
$$
 \Omega' (\zs',\xi_*^n,z)
=\frac{d \zs' h\ov{\zs'}^T\wedge
d\zs^n \wedge(d \zs'\wedge hd\ov {\zs'}^T)^{n-q-2}\wedge(d \zs' \wedge
hd \ov{z'}^T)^q}{-(2i)^{1-n}\ov {N_*^{n-q-1}(\xis)}N_*^{q+1}(\xis) }.
$$
In eliminating $A$,  we have used the fact that
there are no  terms  $\ov\zs^\alpha\ov\zs^\beta$ in
\re{2irzz}-\re{2irzz+}.

\


Take a linear transformation
$\hat\zeta'= U(\zs'),  \hat\xi^n=\xi_*^n$
such that
$$
h(\zs',\ov\zs)=|\hat\zeta'|^2,\quad \hat \zeta^n=\hat\xi^n+i|\hat\zeta'|^2.
$$
Let $\hat z'=U(z')$.
Under the new coordinates and on $|\hat\zeta^n|=1$,  $ \Omega' (\zs',\xi_*^n,z)
$ becomes
\aln
 \widehat
\Omega(\hat\zeta',\hat\xi^n,\hat z')
&=\frac{ (\ov{\hat \zeta'}\cdot d \hat\zeta')\wedge d\hat\zeta^n\wedge
 (d\hat\zeta' \dot{\wedge}   d\ov{\hat \zeta'})^{n-q-2}\wedge (d \hat\zeta' \dot{\wedge}
 d\ov{\hat z'})^q}{-(2i)^{1-n}(\ov {\hat \zeta^n})^{n-q-1}(\hat\zeta^n)^{q+1}}\\
 &=(\ov{\hat \zeta'}\cdot d \hat\zeta')\wedge d\mu(\hat\zeta^n)\wedge
 (d\hat\zeta' \dot{\wedge}   d\ov{\hat \zeta'})^{n-q-2}\wedge
 (d \hat\zeta' \dot{\wedge}
 d\ov{\hat z'})^q.\end{align*}
 Here
 \aln
 \mu(\hat\zeta^n)&=\begin{cases}
 -(2i)^{n-1} (n-2q-1)^{-1}(\hat\zeta^n)^{n-2q-1},&  n-2q>1,\\
 -(2i)^{n-1} \log\hat \zeta^n,& n-2q=1,\\
-(2i)^{n-1} (n-2q-1)^{-1}(\ov{\hat\zeta^n})^{1+2q-n}, & n-2q<1.
\end{cases}
\end{align*}
Note that $\IM\hat\zeta^n\geq0$ and
$\mu(\hat\zeta^n)$ is smooth on $|\hat\zeta^n|<1$.
By Stokes' theorem, we get
\aln
c_qd\ov{\hat z'}^I&\df\int_{|\hat\zeta^n|=1}d\ov{\hat\zeta'}^I
\wedge(\ov{\hat \zeta'}\cdot d \hat\zeta')\wedge d\mu(\hat\zeta^n)\wedge
 (d\hat\zeta' \dot{\wedge}   d\ov{\hat \zeta'})^{n-q-2}\wedge (d \hat\zeta' \dot{\wedge}
 d\ov{\hat z'})^q\\
& 
= \int_{|\hat\zeta^n|<1} d\mu(\hat\zeta^n)\wedge
 (d\hat\zeta' \dot{\wedge}   d\ov{\hat \zeta'})^{n-q-1}\wedge d\ov{\hat\zeta'}^I
\wedge(d \hat\zeta' \dot{\wedge}
 d\ov{\hat z'})^q.
\end{align*}
To rewrite the integrand, introduce  variables $\xi=(\xi_1,\cdots,\xi_q)$,    $\eta$, $x$, $y$   so that
$$
d\ov{\hat\zeta'}^I=d\ov\xi^I, \quad
d\hat\zeta' \dot{\wedge}   d\ov{\hat \zeta'}=d\xi \dot{\wedge}d\ov\xi+
d\eta \dot{\wedge}d\ov\eta,\quad
d \hat\zeta' \dot{\wedge}
 d\ov{\hat z'}=d\xi \dot{\wedge}d\ov x+d\eta\dot{\wedge}d\ov y. $$
Then $(d\xi \dot{\wedge}d\ov\xi+
d\eta \dot{\wedge}d\ov\eta)^{n-q-1}\wedge d\ov{\xi}^I
\wedge(d\xi \dot{\wedge}d\ov x+d\eta\dot{\wedge}d\ov y)^q$ equals
 \aln
& 
(d\eta \dot{\wedge}d\ov\eta)^{n-q-1}\wedge d\ov{\xi}^I
\wedge(d\xi \dot{\wedge}d\ov x+ d\eta\dot{\wedge}d\ov y)^q
=
(d\eta \dot{\wedge}d\ov\eta)^{n-q-1}\wedge d\ov{\xi}^I
\wedge(d\xi \dot{\wedge}d\ov x)^q
\\ & \hspace{6em}=(-1)^q
(d\eta \dot{\wedge}d\ov\eta)^{n-q-1}
\wedge(d\xi \dot{\wedge}d\ov{\xi})^q\wedge d\ov x^I.
 \end{align*}
The last term equals $
(-1)^q\binom{n-1}{q}^{-1}
(d\xi \dot{\wedge}d\ov\xi+d\eta \dot{\wedge}d\ov\eta)^{n-1} \wedge d\ov x^I$. Therefore
$$
c_qd\ov{\hat z'}^I=  (-1)^q\binom{n-1}{q}^{-1}\int_{|\hat\zeta^n|<1}
  d\mu(\hat\zeta^n)\wedge
  (d\hat\zeta' \dot{\wedge}   d\ov{\hat \zeta'})^{n-1}
  \wedge d\ov{\hat z'}^I. 
$$
We now compute $c_q$. By the definition of $\mu(\hat\zeta^n)$, we get
\aln
 c_q&=
 { (-1)^{q+1}4^{n-1}}{(n-1-q)!q!} \int_{| \zeta^n|<1}  (\zeta^n)^{n-2q-2}dV(\xi), \quad n\geq 2q+2,\\
 c_q&= {(-1)^{q}4^{n-1}}{(n-1-q)!q!}
 \int_{| \zeta^n|<1}  (\ov\zeta^n)^{2q-n}dV(\xi), \quad
 n<2q+2.
 \end{align*}
 We first consider the case $n\geq 2q+2$.
Define $\sigma_{n-1}= \f{2\pi^{n-1}}{(n-2)!}$ and $\sigma_{n-1,q}=
   \f{\sigma_{n-1}}{2(n-2q-1)}$. Set $t=\xi^n$ and $s=|\zeta'|$.
Using polar coordinates, we get
\aln
&\mathcal J=\int_{| \zeta^n|<1}( \zeta^n)^{n-2q-2}
  dV(\xi) =\sigma_{n-1}
 \int_{
 -1}^1\int_{s=0}^{(1-t^2)^{{1}/{4}}}
 (t+is^2)^{n-2q-2}s^{2n-3}\, dsdt\\
 & =\f{\sigma_{n-1}}{n-2q-1}\int_{
 0}^1\{((1-s^4)^{\yt}+is^2)^{n-2q-1}-(-
 (1-s^4)^{\yt}+is^2)^{n-2q-1}\}
 s^{2n-3}\, ds\\
 & =\sigma_{n-1,q}\int_{
 -1}^1((1-s^2)^{\yt}+is)^{n-2q-1}
 s^{n-2}\, ds.
 \end{align*}
 Set  $s=\sin \theta$. We get $\mathcal J=\sigma_{n-1,q}\int_{
  - {\pi}/{2}}^{ {\pi}/{2}}
 (\cos\theta+i\sin\theta)^{n-2q-1}\sin^{(n-2)}\theta\,\cos\theta\, d\theta$. Hence
 \aln
 & \mathcal J
 = \f{\sigma_{n-1,q}}{2 (2i)^{n-1}}\int_{|z|=1 }z^{n-2q-1}(z-\ov z)^{n-2}(z+\ov z)\frac{dz}{z}
\\ &
= \f{\sigma_{n-1,q}\pi i (-1)^{n-q-1}}{ (2i)^{n-1}}
\left\{\binom{n-2}{n-q-1}-\binom{n-2}{n-q-2}\right\}
=
 \f{\sigma_{n-1}\pi i(-1)^{n-q}(n-2)! }{2(2i)^{n-1}(n-q-1)!q!},
 \end{align*}
where $\binom{n-2}{n-1}=\binom{n-2}{ -1}=0$. We get
$
c_q=\yt(2\pi i)^n
$
for $n\geq2q+2$.
For $n<2q+2$, we can verify that
$c_q=(-1)^n\ov c_{n-1-q} =\yt(2\pi i)^n$.
The proof of \rl{ap1} is complete.
%


\

We need the following lemma from~\ci{Hesese}, where
$|F(\zeta,z)|=t$ is replaced by $|\zeta-z|=t$
and only $r\in \cC^2$ is needed.
\le{b3} Let $1\leq q\leq n-1$. Let $M, \rho_0$ be as in  \rl{ap1}. Assume that $r$
is of class $\cC^3$.
 Let $\varphi$ be a continuous
$(0,q)$-form   on $M_{\rho_0}$. If $0<\rho<\rho_0$ and $t>0$
is sufficiently small,
then on $M_\rho$, in the sense of currents, and modulo $\db r(z)$, 
$$
 \int_{ F(\zeta,z)>t}\ \varphi(\zeta)\wedge
\db_z\Omega^{+-}_{(0,q-1)}(\zeta,z)=(-1)^{q-1}\db_b\int_{ F(\zeta,z)>t}\ \varphi(\zeta)\wedge
 \Omega^{+-}_{(0,q-1)}(\zeta,z).
$$
\end{lemma}
\begin{proof}
  For the convenience
of the reader, we  reproduce the proof in~\ci{Hesese}.
  Since $r\in \cC^3$, for $t$   sufficiently small, both
integrals are  continuous on $M_\rho$. To verify the identity
in the sense of  currents,
let $\psi$ be a smooth $(n,n-1-q)$-form  with compact support in $M_\rho$.
Interchanging the order of integration,
we have
\aln
\mathcal I(z)&= \int_{z\in M}\int_{ F(\zeta,z)>t}\ \varphi(\zeta)\wedge
\db_z\Omega^{+-}_{(0,q-1)}(\zeta,z)\wedge\psi(z)
 \\ &\quad\quad =-\int_{\zeta\in M}\int_{ F(\zeta,z)>t}\ \varphi(\zeta)\wedge
 \db_z\Omega^{+-}_{(0,q-1)}(\zeta,z)\wedge\psi(z)\quad \text{(by Fubini)}
 \\ &\quad\quad=\int_{\zeta\in M}\int_{ F(\zeta,z)>t}\ \psi(z)\wedge
  \db_z\Omega^{+-}_{(0,q-1)}(\zeta,z)\wedge\varphi(\zeta).
 \end{align*}
  Applying Stokes' theorem 
  yields
 \aln
  &\mathcal I(z) =(-1)^{q}\int_{\zeta\in M}\int_{ F(\zeta,z)>t}\ \db_z\psi(z)\wedge
  \Omega^{+-}_{(0,q-1)}(\zeta,z)\wedge\varphi(\zeta)
  \\ &\quad\quad \quad+(-1)^{q+1}\int_{\zeta\in M}\int_{ F(\zeta,z)=t}\  \psi(z)\wedge
  \Omega^{+-}_{(0,q-1)}(\zeta,z)\wedge\varphi(\zeta) 
  \\ &\quad\quad = \int_{\zeta\in M}\int_{ F(\zeta,z)>t}\ \varphi(\zeta)\wedge
  \Omega^{+-}_{(0,q-1)}(\zeta,z)\wedge\db_z\psi(z)
  \\ &\quad\quad \quad+(-1)^{q+1}\int_{\zeta\in M}\int_{ F(\zeta,z)=t}\ \varphi(\zeta) \wedge
 \Omega^{+-}_{(0,q-1)}(\zeta,z)\wedge\psi(z).
 \end{align*}
Interchanging the  order of integration in both terms   yields
  \aln
&\mathcal I(z) = -\int_{z\in M}\int_{ F(\zeta,z)>t}\ \varphi(\zeta)\wedge
\Omega^{+-}_{(0,q-1)}(\zeta,z)\wedge\db_z\psi(z)
\\ &\quad\quad \quad+(-1)^{q+1}\int_{z\in M}\int_{ F(\zeta,z)=t}\ \varphi(\zeta) \wedge
\Omega^{+-}_{(0,q-1)}(\zeta,z)\wedge\psi(z)
\\ &\quad\quad = -\int_{z\in M}\int_{ F(\zeta,z)>t}\ \varphi(\zeta)\wedge
\Omega^{+-}_{(0,q-1)}(\zeta,z)\wedge\db_z\psi(z).
\end{align*}
Here the second last term vanishes  by counting total degree in $\zeta$,  which equals $q+2n-q-1>2n-2$.
\end{proof}


As in
 \ci{Weeinia}, the above two lemmas
can be used to derive
the Henkin homotopy formula for
$(0,q)$-forms on $\ov M_\rho$  via Stokes' theorem,  when $r\in \cC^3$
satisfies \re{bhrooap}. Here are the details. Let $M^t(z)=M_{\rho_0}
\cap\{\zeta\colon
|F(\zeta,z)|>t\}$. We will use
$
\ov\pd_\zeta\Omega_{(0,q)}^{+-}+\ov\pd_z\Omega_{(0,q)}^{+-}=0$ for $
\zeta\neq z$ and $1\leq q\leq n-2$; see~\ci{Weeinia}.
Then on $M_\rho$
\aln
 &\int_{\pd M^t}\  \varphi(\zeta)\wedge\Omega_{(0,q)}^{+-}(\zeta,z)=
\int_{M^t(z)}\ \ov\pd_\zeta\varphi \wedge\Omega_{(0,q)}^{+-}
+(-1)^q\int_{M^t(z)}\ \varphi(\zeta)\wedge\ov\pd_\zeta\Omega_{(0,q)}^{+-} \\
&\hspace{10ex}=\int_{M^t(z)}\ \ov\pd_\zeta\varphi \wedge\Omega_{(0,q)}^{+-}(\zeta,z)
+(-1)^{q-1}
\int_{M^t(z)}\ \varphi(\zeta)\wedge\ov\pd_z\Omega_{(0,q-1)}^{+-}(\zeta,z).
\end{align*}
 Using Lemmas \ref{ap1}-\ref{b3} and letting
  $t\to0$, we obtain, modulo $\ov\pd r(z)$,
$$
\int_{\pd M_{\rho_0}}\ \varphi(\zeta)\wedge\Omega_{(0,q)}^{+-}(\zeta,z)
=c_0\varphi(z)+\int_{M_{\rho_0}}\ \ov\pd_\zeta\varphi\wedge\Omega_{(0,q)}^{+-}
+\db_b\int_{M_{\rho_0}}
\varphi(\zeta)\wedge \Omega_{(0,q-1)}^{+-}.
$$
Now assume that
$0< q< n-2$. Then
$-\Omega_{(0,q)}^{+-}=\db_\zeta\Omega_{(0,q)}^{0+-}+\db_z\Omega_{(0,q-1)}^{0+-}$;
see~\ci{Weeinia}.
We get that modulo $\ov\pd r(z)$
\aln
&\int_{\pd M_{\rho_0}}\varphi(\zeta)\wedge\Omega_{(0,q)}^{+-}(\zeta,z)
=-\int_{\pd M_{\rho_0}}
\varphi(\zeta)\wedge\ov\pd_\zeta\Omega_{(0,q)}^{0+-}
-\int_{\pd M_{\rho_0}}
\varphi(\zeta)\wedge\ov\pd_z\Omega_{(0,q-1)}^{0+-}\\
&\hspace{16ex}=(-1)^{q}\int_{\pd M_{\rho_0}}
 \dbb\varphi(\zeta)\wedge\Omega_{(0,q)}^{0+-}
-(-1)^{q} \dbb\int_{\pd M_{\rho_0}}
\varphi(\zeta)\wedge\Omega_{(0,q-1)}^{0+-}.
\end{align*}
Therefore, modulo $\db r(z)$ and as currents,
\aln
c_0\varphi(z)&=
\ov\pd_b\Bigl\{-\int_{M_{\rho_0}}\varphi(\zeta)\wedge\Omega_{(0,q-1)}^{+-}(\zeta,z)-(-1)^{q}
\int_{\pd M_{\rho_0}}\varphi(\zeta)\wedge\Omega_{(0,q-1)}^{0+-}(\zeta,z)\Bigr\}\\
&\quad +\Bigl\{-\int_{M_{\rho_0}}\db_b\varphi(\zeta)\wedge\Omega_{(0,q)}^{+-}(\zeta,z)+(-1)^{q}
\int_{\pd M_{\rho_0}}\db_b\varphi(\zeta)\wedge\Omega_{(0,q)}^{0+-}(\zeta,z)\Bigr\}\\
&=-
\ov\pd_b\Bigl\{\int_{M_{\rho_0}}\Omega_{(0,q-1)}^{+-}(\zeta,z)\wedge\varphi(\zeta)+
\int_{\pd M_{\rho_0}}\Omega_{(0,q-1)}^{0+-}(\zeta,z)\wedge\varphi(\zeta)\Bigr\}\\
&\quad -\Bigl\{\int_{M_{\rho_0}}\Omega_{(0,q)}^{+-}(\zeta,z)\wedge\db_b\varphi(\zeta)+
\int_{\pd M_{\rho_0}}\Omega_{(0,q)}^{0+-}(\zeta,z)\wedge\db_b\varphi(\zeta)\Bigr\}.
\end{align*}

We now derive the homotopy formula by reducing the regularity condition.
  Notice that when $M$
 is $\cC^2$ and strictly convex
 it
is  proved by Henkin~\ci{Weeinia}.

\medskip

\noindent
{\bf Case a) $r\in \cC^2$ and $\varphi\in \cC^1(\ov M_{\rho_0})$.}
Let $0< q<n-2$. Fix a tangential $(0,q)$-form $\varphi$
of class $\cC^1$ on $M_{\rho_0}$.

 Write $x=(\RE z,\IM z')$ and $\xi=(\RE\zeta,\IM\zeta')$ with   $\zeta,z\in M_{\rho_0}$.
We first consider the case
when $\varphi$ has compact support, say in $M_{\rho_1}$ with $\rho_1<\rho_0$. Then $P'_{M_{\rho_0}}
\varphi $ has no boundary term.
Recall that
$ \Psi(\xi,x)=(\tilde\psi_x(\xi),x)=(\xis,x)$ is defined by relations
  $\zs'=\zeta'-z',
\zs^n=-2ir_z(\zeta-z)$ and $\zeta,z\in M_{\rho_0}$.
 By \re{t1=1}-\re{t2=1} and \re{firstfa}, on $D_{\rho_0}$  
 \gan
 P'_{0, M_{\rho_0}} \varphi(x)=\sum_{|I|=q-1}
  \mathcal I_{\ov I} (x)\, d\ov{z'}^I,\\
 \mathcal I_{\ov I} (x)=
\sum_{|K|=1}\sum_{|J|=q} \int_{B_{9\rho_0}}
\left((\varphi_{\ov J} \tilde A_{\ov IK}^{\ov  J} )\circ\Psi^{-1}\cdot {\hat T_1^{-a}}\cdot
 {\hat T_2^{-b}}\right)(\xis , x)\cdot\hat k_{ab}^K(\xis)\, dV(\xis),\\
 \varphi_{\ov J}(\xi,x)=\varphi_{\ov J}(\xi),\quad \hat T_j(\xi_*,x)=
 1+ \sum_{|L| =2} C_{jL}\circ\Psi^{-1}(\xi_*,x)
 {\xis^L}{N_*^{-1}(\xis)}.
\end{gather*}
Here $\tilde  A^{\ov J}_{\ov IK}, C_{jL}$ are functions of the form $\pd_*^2r$, $\pd_*^2\hat r$,
respectively,   defined in section~\ref{sec8-}.
These functions depend only on derivatives of $\hat r$ of order at most two.
We   take a sequence of $\cC^\infty$
functions $\hat r^m $  converging to $\hat r$ in $\cC^2$-norm on $\ov D_{\rho_0}$.
Subtracting $\hat r^m$ by its Taylor polynomial of order $1$ about the origin, we may assume that
 $\hat r^m(0)=0$ and $\pd\hat r^m(0)=0$.
In what follows we use  the letter $m$   to  indicate   dependence   on
$\hat r^m$.
Let $M^m$ be the graph $y^n=|z'|^2+\hat r^m(x)$ over $D_{\rho_0}$.
 Let  $\rho\in(\rho_1,\rho_0)$. There exists $m_* $ such that
 for $m>m_*$
$$
D_{\rho_1}\subset
  \ov D^m_{\rho}\subset D_{\rho_0},\quad  D^m_{\rho}\df
 \{x\in D_{\rho_0}\colon|x|^2+\hat r^m(x)<
 \rho^2\}.
 $$
 Assume that $m>m_*$.
Then by $
\|\hat r^m \|_{ D_{\rho_0},2}<1/C_0$,  $D_{\rho}^m$ is strictly convex. Note that $
 \varphi$ is still
 a tangential $(0,q)$-form of $M^m$
   on $D_{\rho_0}$.
 By the \hf\ for the smooth $M^m$,
 \eq{hfmm}
\varphi=\db_{M^m}P'_{M^m_{\rho}}\varphi+Q'_{M^m_{\rho}}\db_{M^m}
\varphi. 
\end{equation}
 By the formula of $\hat T_j$, it is clear that $|\hat
T_{j,m}|\geq1/2$  for  $m$ sufficiently large and  $0\neq \xis\in
B_{9\rho}\cap\tilde\psi^m_x(D_{\rho})$. Since $\varphi_{\ov J}$ has
compact support, we have
$$|\{(\varphi_{\ov J} \tilde  A^{Jm}_{I K})\circ\Psi_m^{-1} \cdot {\hat T_{1,m}^{-a}}\cdot
 {\hat T_{2,m}^{-b}}\}(\xis , x)|\leq C, \quad\xis\neq0.
$$
Since $\hat k_{ab}^I\in L^1_{loc}$, by the dominated convergence theorem,
$\mathcal I^m_{I J}$ converges to $\mathcal I_{I J}$ pointwise on $\ov D_{\rho_1}$  as $m\to\infty$.
Thus,
 $P'_{M_\rho^m}\varphi$ converges to $P'_{M_\rho}\varphi$ pointwise on $\ov D_{\rho_1}$.
 By the $\cC^0$-estimate on $P_0'$ (see \re{fk0+}), we know that $|P'_{M^m_{\rho}}\varphi|<C$ on $\ov D_{\rho_1}\subset
\cap_mD_{\rho}^m$.
Since $\varphi$ is of class $\cC^1$, $\db_{M^m}\varphi$  converges to $\db_{M}\varphi$
uniformly on $\ov D_{\rho_0}$.
Next, we fix an $(n,n-q-1)$-form $\psi=\sum_{|J|=n-q-1}\psi_{\ov J}dz^1\wedge
\cdots\wedge dz^{n-1}\wedge\pd r(z)\wedge d\ov{z'}^J$ of class $\cC^1$ and   compact support
in $D_{\rho_1}$.  Then $\psi^m=\sum_{|J|=n-q-1}\psi_{\ov J}dz^1\wedge
\cdots\wedge dz^{n-1}\wedge\pd r^m(z)\wedge d\ov{z'}^J$ is an $(n,n-1-q)$-form
  of $M^m$ on $D_\rho^m$. Obviously,
 $\psi^m$ converges to $\psi$ in $\cC^1$-norm
uniformly on $\rr^{2n-1}$.
Now one can see that
  $\int_{D_{\rho_1}}P'_{0,M_{\rho}^m}\varphi\wedge \db_{M^m}\psi^m$
 converges to $\int_{D_{\rho_1}} P'_{0,M_{\rho}}\varphi\wedge \db_{M}\psi$.
 A similar argument shows that
 $Q'_{0,M_{\rho}^m}\db_{M^m}\varphi$ are uniformly bounded and
 converge to $Q'_{0,M_{\rho}}\db_{M}\varphi$ pointwise  on $D_{\rho_1}$.
Thus $ \int_{D_{\rho_1}}( Q'_{0,M^m_{\rho}}\db_{M^m} \varphi) \wedge\psi^m$
converges to $ \int_{D_{\rho_1}}( Q'_{0,M_{\rho}}\db_{M} \varphi) \wedge\psi$.
 By the \hf\  \re{hfmm},
\eq{hfmm+}
\int_{D_{\rho_1}} \varphi \wedge\psi^m =(-1)^q\int_{D_{\rho_1}} (P'_{0,M^m_{\rho}}\varphi) \wedge\db_{M^m}\psi^m
+\int_{D_{\rho_1}}( Q'_{0,M^m_{\rho}}\db_{M^m} \varphi) \wedge\psi^m.
\end{equation}
Taking limits, we  see that as currents,
$
\varphi=\db_{M}P'_{0,M_{\rho}}\varphi+Q'_{0,M_{\rho}}\db_{M}\varphi
$ holds on $ D_{ \rho_1}.$  Since $\varphi$ has compact support
in $D_{\rho_1}\subset D_\rho$,   we     replace the domain of integration $M_{\rho} $ by $M_{\rho_0}$ and
add boundary integrals.
We get
 $\varphi=\db_{M}P'_{M_{\rho_0}}\varphi+Q'_{M_{\rho_0}}\db_{M}\varphi
$ on $D_{\rho_1}$ as currents,
 and hence on $D_{\rho_0}$ whenever $\varphi$ has compact support in $D_{\rho_0}$.

Return to the general case.
Let $M^m$ be as before.  Take any $\rho_1,\rho_2,\rho$ such that $0<\rho_1<\rho_2<\rho<\rho_0$.
Take a $\cC^\infty$ function $\chi$ which is $1$ on $D_{\rho_2}$   and has
compact support in $D_{\rho_0}$.
Let $\varphi_0=\chi\varphi$ and
$\varphi_1=(1-\chi)\varphi$. We have proved that
$ \varphi_0 =\db_{M}P'_{M_{\rho_0}}\varphi_0+Q'_{M_{\rho_0}}\db_{M}
 \varphi_0$
 as currents  on $D_{\rho_0}$. Let $\psi,\psi^m$ be as before, which  are supported in $D_{\rho_1}$.
  By an analogy of \re{hfmm+}, for $m>m_*$ we have
$$ 
\int_{D_{\rho_1}} \varphi_1 \wedge\psi^m =(-1)^q\int_{D_{\rho_1}} (P'_{M^m_{\rho}}\varphi_1) \wedge\db_{M^m}\psi^m
+\int_{D_{\rho_1}}( Q'_{M^m_{\rho}}\db_{M^m} \varphi_1) \wedge\psi^m.
$$ 
Since $\supp\psi^m\subset D_{\rho_1}$ and $\supp\varphi_1\subset
\ov D_{\rho_0}\setminus D_{\rho_2}$, 
  the dominated convergence theorem implies
the convergence of $
 \int_{D_{\rho_1}} (P'_{0, M^m_{\rho}}\varphi_1) \wedge\db_{M^m}\psi^m $ to $
\int_{D_{\rho_1}} (P'_{0,M_{\rho}}\varphi_1) \wedge\db_{M}\psi.
$
The integrands for the boundary integrals $(P'_{1,M^m_{\rho}}\varphi_1)(x)$
are of class $\cC^1$ in a neighborhood of  $\pd D_\rho$,
when $x\in D_{\rho_2}$.  Note that
the map $d_\rho$
 sending  $x\in\pd D_\rho$ to $r_\rho(x)x\in\pd D^m_\rho$ with $r_\rho>0$
converges to the identity map in the $\cC^2$ norm as $m\to\infty$.
 Thus we can  verify that $P'_{1,
M^m_{\rho}}\varphi_1$ converges uniformly to $
P'_{1,M_{\rho}}\varphi_1$ on $D_{\rho_1}$. Combining with
$\supp \psi^m\subset D_{\rho_1}$, we obtain
 the convergence of $ 
\int_{D_{\rho_1}} (P'_{1, M^m_{\rho}}\varphi_1)
\wedge\db_{M^m}\psi^m $ to $ \int_{D_{\rho_1}} (P'_{1,M_{\rho}}\varphi_1)
\wedge\db_{M}\psi .
$
One can verify that  $ 
\int_{D_{\rho_1}} (Q'_{1, M^m_{\rho}}\db_{M^m}\varphi_1)
\wedge \psi^m $ converges to $\int_{D_{\rho_1}} (Q'_{1,M_{\rho}}\db_{M}\varphi_1)
\wedge \psi.
$
Hence  $\varphi_1 =\db_{M}P_{M_{\rho}}\varphi_1+Q_{M_{\rho}}\db_{M}
 \varphi_1$
  as currents
  on $D_{\rho_1}$. The definition of $\varphi_1$ is independent of $\rho>\rho_2$.
  Letting $\rho\to\rho_0$, we get $\varphi_1 =\db_{M}P_{M_{\rho_0}}\varphi_1+Q_{M_{\rho_0}}\db_{M}
 \varphi_1$  on $D_{\rho_1}$. Add to
 $\varphi_0$. We get $\varphi=\db_{M}P_{M_{\rho_0}}\varphi+Q_{M_{\rho_0}}\db_{M}
 \varphi$ as currents on $D_{\rho_1}$, and hence on $D_{\rho_0}$.
%
%

\medskip

\noindent
{\bf Case b) $r\in \cC^2(D)$
and $\varphi,\dbm\varphi\in \cC^0(\ov M_{\rho_0})$.}
We verify the \hf\
  by the Friedrichs approximation theorem, for which we need
   the   commutator of a smoothing operator $S_t$
  and $\dbm$, applied to tangential $(0,q)$-forms.

 Recall that on $M$
\eq{thpa}
\theta(\xi)=-2i\pd r(\zeta)=a \, d\xi^n
\mod{(d\zeta^\alpha,d\ov \zeta^\beta)},\quad a=1+\hat r_{\xi^n}^2.
\end{equation}
Let $\varphi=\sum_{|I|=q}
  \varphi_{\ov I}(x) d\ov{z'}^I$ be a continuous tangential $(0,q)$ form on $M$. Let
  $\chi$ be a smooth function of compact support in $\rr^{2n-1}$ such
  that $\int\chi dV=1$. Let $\chi_t(x)=t^{1-2n}\chi(t^{-1}x)$ and
  define $
  S_t\varphi=\sum_{|I|=q} \varphi_{\ov I}\ast\chi_t\,d\ov{z'}^I.
  $
Recall that $\ov X_\alpha=\partial_{\ov{z^\alpha}}+b_{\ov \alpha}
\partial_{x^n}$ with $b_{ \ov\alpha}= -  r_{\ov{z^\alpha}}/{(2r_{\ov{z^n}})}$.
\le{dbmst}
Let $M\subset\cc^n$ be a graph of class $\cC^2$ over $D\subset\cc^{n-1}\times\rr$.
 Assume that
$\varphi$ is a
\cont\  tangential $(0,q)$-form on $M_\rho$ such that $\dbm\varphi$
is \cont\ on $M_\rho$.  Let $0<\rho'<\rho$.  Then for $t$ sufficiently small
and on $M_{\rho'}$,
\aln
 [\dbm, S_t]\varphi(x)=
\sum_{|I|=q,1\leq\alpha<n}\ \int&\Bigl\{\varphi_{\ov I}(x-t\xi)t^{-1}(b_{\ov\alpha}(x)-b_{\ov\alpha}(x-t\xi))
(\pd_{\xi^n}\chi)(\xi)
\Bigr.\\
\Bigl.
\quad &  -(\varphi_{\ov I}a^{-1} \ov X_\alpha
a)(x-t\xi)\chi(\xi)\Bigr\}
\,
dV(\xi)\wedge d\ov{z^\alpha}\wedge d\ov{z'}^I.\end{align*}
\end{lemma}
\begin{proof} We follow a computation in~\ci{Frfofo}; see also~\ci{Honize} (p.~121).
 Let $
\varphi=\sum_{|I|=q}\varphi_{\ov I}d\ov{z'}^I$.
By the assumption, we have $\dbm\varphi =\sum_{|J|=q+1} \psi_{\ov J} d\ov{z'}^J$ with
$\psi_{\ov J}\in \cC^0.$  Now
\aln
 \dbm S_t\varphi(x)&=
\sum_{|I|=q}\sum_{1\leq\alpha<n}\ \int\varphi_{\ov I}(\xi)\ov X_\alpha^{(x)}\chi_t(x-\xi)\,
dV(\xi)\wedge d\ov{z^\alpha}\wedge d\ov{z'}^I,\\
S_t \dbm \varphi(x)&=
\sum_{|J|=q+1} \ \int \psi_{\ov J}(\xi)\chi_t(x-\xi)\,
dV(\xi)\wedge
  d\ov{z'}^J.
\end{align*}
Set $\nu=d\zeta^1\wedge\cdots\wedge d\zeta^{n-1}\wedge d\xi^n$.
We may assume that $dV=d\ov {\zeta^1}\wedge\cdots\wedge d\ov {\zeta^{n-1}}
\wedge \nu$.
For each $J=(j_1,\ldots, j_{q+1})$, 
there exist increasing indices $J^*=(j_1^*,\ldots, j^*_{n-2-q})$ and $\epsilon^J=\pm1$
such that
$
d\ov {\zeta^1}\wedge\cdots\wedge
 d\ov {\zeta^{n-1}}
 =\epsilon^Jd\ov{\zeta'}^J \wedge d\ov{\zeta'}^{J^*}.
$
Now,
\aln
S_t \dbm \varphi(x)&=
 \sum_{|J|=q+1} \epsilon^J\int \psi_{\ov J}(\xi) \chi_t(x-\xi)\,d\ov{\zeta'}^J
\wedge d\ov{\zeta'}^{J^*}\wedge
\nu\wedge d\ov{z'}^J\\
& =
 \sum_{|J|=q+1} \epsilon^J\int \dbm^{(\xi)}\varphi  \wedge\chi_t
 (x-\xi)\,  d\ov{\zeta'}^{J^*}\wedge
\nu\wedge d\ov{z'}^J.
\end{align*}
By Stokes' formula,
 \re{dbnq} and \re{thpa}, $S_t \dbm \varphi(x)$ equals
\aln
&
 \sum_{|J|=q+1}\  (-1)^{q+1}\epsilon^J\int  \varphi(\xi)   \wedge
 a(\xi)\dbm^{(\xi)}(a^{-1}(\xi)\chi_t(x-\xi))\wedge d\ov{\zeta'}^{J^*}\wedge
\nu\wedge d\ov{z'}^J\\
& \qquad  =
 -\sum_{|I|=q}\sum_{|J| =q+1}  \epsilon^J\int \varphi_{\ov I}  (\xi)
  a(\xi)\ov X_\alpha^{(\xi)}\f{\chi_t(x-\xi)}{a (\xi)}\,
  d\ov\zeta^\alpha\wedge d\ov{\zeta'}^I\wedge
d\ov{\zeta'}^{J^*}\wedge
\nu\wedge d\ov{z'}^J\\
&\qquad =-
\sum_{|I|=q}\sum_{1\leq\alpha<n}\ \int\varphi_{\ov I}(\xi)a(\xi)\ov X_\alpha^{(\xi)}(a^{-1}(\xi)\chi_t(x-\xi))\,
dV(\xi)\wedge d\ov{z^\alpha}\wedge d\ov{z'}^I.
\end{align*}
Set $E_\alpha(\xi,x) =  a(\xi)\ov X_\alpha^{(\xi)}
(a^{-1}(\xi)\chi_t(x-\xi))
+\ov X_\alpha^{(x)}(\chi_t(x-\xi)).$ We get
$$ [\dbm, S_t]\varphi(x)=
\sum_{|I|=q}\ \sum_{1\leq\alpha<n} \ \int\varphi_{\ov I}(\xi)E_\alpha(\xi,x)
\,
dV(\xi)\wedge d\ov{z^\alpha}\wedge d\ov{z'}^I,$$
which can be put into the  form in the lemma.
\end{proof}
We now derive the \hf\  for  the case b)  via smoothing.  We  have
\aln
\dbm S_t\varphi-\dbm\varphi
&=S_t\dbm\varphi-\dbm\varphi+[\dbm, S_t](\varphi-S_{t'}\varphi)+[\dbm, S_t] S_{t'}\varphi.
\end{align*}
Fix $\rho<\rho_1<\rho_0$ and $\epsilon>0$.
Fix
$t'>0$ sufficiently small such that $|S_{t'}\varphi-\varphi|<\epsilon$ on $M_{\rho_1}$.
By \rl{dbmst}, for all small $t$
\eq{smrd}
\sup_{M_\rho}\bigl|[\dbm,S_t](\varphi-S_{t'}\varphi)\bigr|\leq C
\sup_{M_{\rho_1}}|\varphi-S_{t'}\varphi|\leq C\epsilon.
\end{equation}
Let $ S_{t'}\varphi=\sum_{|I|=q}\psi_{\ov I}\, d\ov{z'}^I$. Then $[\dbm,S_t]=\sum
[b_{\ov\alpha}\pd_{x^n},S_t]$
and
$$
([\dbm,S_t]\psi)(x)
=\sum_{1\leq\alpha<n}\sum_{|I|=q}\
\int\chi_{t}(y)(b_{\ov\alpha}(x)-b_{\ov\alpha}(x-y))\pd_{x^n}\psi_{\ov I}(x-y)\, dV(y)\, d\ov {z'}^I.
$$
This shows that $\lim_{t\to0}[\dbm,S_t]S_{t'}\varphi=0$ on $M_{\rho}$.
We also have $\lim_{t\to0}S_{t}\dbm\varphi=\dbm\varphi$ on $M_{\rho}$.
Therefore, \re{smrd} implies that
$\lim_{t\to0}\dbm S_t\varphi=\dbm\varphi$ on $M_\rho$.
Now by case a),
$ \varphi_t=\dbm
P'_{M_{\rho}}\varphi_t+Q'_{M_{\rho}}
\dbm\varphi_t $ holds
on $M_{\rho}$  as currents. Letting $t\to0$ and then
$\rho\to\rho_0$, we obtain    $\varphi=\dbm P'_{M_0}\varphi+Q'_{M_0}\dbm \varphi$
on $M_{\rho_0}$  in the sense of currents.

\newcommand{\Wenion}{\bibitem{Wenion}
S.M. Webster, {\em The integrability problem for CR vector bundles},
Several complex variables and complex geometry, Part 3 (Santa
Cruz, CA, 1989),  355--368, Proc. Sympos. Pure Math., 52, Part 3,
Amer. Math. Soc., Providence, RI, 1991. }

\newcommand{\Weeini}{\bibitem{Weeini}
S.M. Webster, {\em A new proof of the Newlander-Nirenberg
theorem}, Math. Z. {\bf 201}(1989), no. 3, 303--316. }

\newcommand{\Wezeon}{\bibitem{Wezeon}
S.M. Webster, {\em
SubRiemannian geometry and the eikonal equation},
 Selected topics in Cauchy-Riemann geometry,  351--378, Quad. Mat., 9,
Dept. Math., Seconda Univ. Napoli, Caserta, 2001.
}

\newcommand{\Weeinia}{\bibitem{Weeinia}
S.M. Webster, {\em On the local solution of the tangential
Cauchy-Riemann equations},  Ann. Inst. H. Poincar\'e Anal. Non
Lin\'eaire  {\bf 6}(1989),  no. 3, 167--182. }

\newcommand{\Weeinib}{\bibitem{Weeinib}
S.M. Webster, {\em On the proof of Kuranishi's embedding theorem},
Ann. Inst. H. Poincar\'e Anal. Non Lin\'eaire  {\bf 6}(1989),  no. 3,
183--207. }

\newcommand{\MMnifi}{\bibitem{MMnifi} J. Michel and
 L. Ma,
 {\em On the regularity of CR structures for
almost CR vector bundles},  Math. Z.  {\bf 218}(1995),  no. 1,
135--142. }

\newcommand{\MMnifo}{\bibitem{MMnifo} L. Ma and J. Michel,
 {\em Regularity of local
 embeddings of strictly pseudoconvex CR structures},
 J. Reine Angew. Math.  {\bf 447}(1994), 147--164.
 }

\newcommand{\MMnith}{\bibitem{MMnith} L. Ma and
 J. Michel, {\em
Local regularity for the tangential Cauchy-Riemann complex},
 J. Reine Angew. Math. {\bf 442}(1993), 63--90.
 }

\newcommand{\Kueitw}{\bibitem{Kueitw}M. Kuranishi, {\em
Strongly pseudoconvex CR structures over small balls},
 {\em I. An a
priori estimate}, Ann. of Math. (2) {\bf 115}(1982), no. 3, 451--500;
 {\em II. A regularity theorem},
  {\bf 116}(1982), no. 1, 1--64;
 {\em III. An embedding theorem,}  {\bf 116}(1982), no. 2, 249--330.
 }

\newcommand{\Akeise}{\bibitem{Akeise}
  T. Akahori, {\em A new approach to the local embedding theorem of
  CR-structures for $n\geq 4$
  (the local solvability for the operator $\overline\partial\sb b$ in
  the abstract sense)},
  Mem. Amer. Math. Soc.  {\bf  67}(1987),  no. 366, xvi+257 pp.
  }

\newcommand{\NReini}{\bibitem{NReini}
A. Nagel and J.-P. Rosay, {\em Nonexistence of homotopy formula
for $(0,1)$ forms on hypersurfaces in $C\sp 3$},  Duke Math. J.
{\bf 58}(1989),  no. 3, 823--827. }

\newcommand{\Canifo}{\bibitem{Canifo}
D. Catlin, {\em Sufficient conditions for the extension of CR structures},
 J. Geom. Anal. {\bf 4}(1994),  no. 4, 467--538.
 }

\newcommand{\Shnise}{\bibitem{Shnise}
M.-C. Shaw, {\em Homotopy formulas for $\overline\partial\sb b$ in
CR manifolds with mixed Levi signatures},
 Math. Z.  {\bf 224}(1997),  no. 1, 113--135.
 }

\newcommand{\Shnize}{\bibitem{Shnize}
M.-C. Shaw, {\em
  $L\sp p$ estimates for local solutions of $\overline\partial\sb {\rm b}$
on strongly pseudo-convex CR manifolds},
  Math. Ann.  {\bf 288}(1990),  no. 1, 35--62.}

\newcommand{\FSsefo}{\bibitem{FSsefo}
G.B. Folland and E.M. Stein, {\em
 Estimates for the $\bar \partial \sb{b}$ complex and analysis on the Heisenberg group},
Comm. Pure Appl. Math.  {\bf 27}(1974), 429--522.
}

\newcommand{\FSeitw}{\bibitem{FSeitw}
G.B. Folland and E.M. Stein,  
{\em
Hardy spaces on homogeneous groups},
Mathematical Notes, 28. Princeton
University Press, Princeton, N.J.; University of Tokyo Press, Tokyo, 1982.
}

\newcommand{\Ronion}{\bibitem{Ronion}
C. Romero, {\em Potential theory for the Kohn Laplacian on the Heisenberg
group}, thesis, University of Minnesota, 1991.}

\newcommand{\Frfofo}{\bibitem{Frfofo} K.O. Friedrichs,
{\em The identity of weak and strong extensions of differential operators},
Trans. Amer. Math. Soc. {\bf 55}(1944). 132--151
}

\newcommand{\Honizea}{\bibitem{Honizea}
L. H\"ormander, {\em
The analysis of linear partial differential operators.
 I. Distribution theory and Fourier analysis.}
 Springer-Verlag, Berlin, 1990.}

\newcommand{\Hosesi}{\bibitem{Hosesi}
L. H\"ormander, {\em
The boundary problems of physical geodesy},
Arch. Rational Mech. Anal. {\bf 62} (1976), no. 1, 1--52.
}

\newcommand{\Honize}{\bibitem{Honize}
L. H\"ormander,  
 {\em An introduction to complex analysis in several
variables},
North Holland, Amsterdam, 1990, 3rd edition.
}

\newcommand{\CSzeon}{\bibitem{CSzeon} S.-C. Chen and M.-C. Shaw,
{\it Partial differential equations in several complex variables},
AMS/IP Studies in Advanced Mathematics, 19,
Amer. Math. Soc., Providence, RI, International Press, Boston, MA,
2001.}

\newcommand{\BBeise}{\bibitem{BBeise}
J. Bruna and J.M. Burgu\'es,
{\em
Holomorphic approximation and estimates for the
$\overline\partial$-equation on strictly pseudoconvex
nonsmooth domains},  Duke Math. J.,  {\bf 55}(1987),  no. 3, 539--596.
}

\newcommand{\Sesifo}{\bibitem{Sesifo}
R.T. Seeley, {\em Extension of $C\sp{\infty }$ functions defined in a half space},
Proc. Amer. Math. Soc.  {\bf 15}(1964), 625--626.
}

\newcommand{\Sisefo}{\bibitem{Sisefo}
Y.T. Siu, {\em
The $\dbar$ problem with uniform bounds on derivatives},
Math. Ann. {\bf 207}(1974), 163--176.}

\newcommand{\Hesese}{\bibitem{Hesese}
G.M. Henkin, {\em The Lewy equation and analysis on pseudoconvex
domains}, Russian. Math. Surveys. {\bf 32}(1977) 59-130.}

\newcommand{\HRseon}{\bibitem{HRseon}
G.M. Henkin  and A.V. Romanov, {\em Exact H\"older estimates
of the solutions of the $\db$-equation},
 Izv. Akad. Nauk SSSR Ser. Mat.  {\bf 35}(1971), 1171-1183.}

\newcommand{\LReize}{\bibitem{LReize}
I. Lieb and R.M. Range,
{\it L\"osungsoperatoren f\"ur den Cauchy-Riemann-Komplex mit
${\mathcal C}\sp{k}$-Absch\"atzungen},
Math. Ann. {\bf 253}(1980), no. 2, 145--164.
}

\newcommand{\GLseei}{\bibitem{GLseei}
J.B. Garnett and R.H. Latter, {\it
The atomic decomposition for Hardy spaces in several complex variables},
Duke Math. J. {\bf 45}(1978), no. 4, 815--845.
}

\newcommand{\GWzese}{\bibitem{GWzese}
X. Gong and S.M. Webster,
{\it Regularity for the CR vector bundle problem I}, Pure and Appl.
Math. Quarterly, (to appear).}

\newcommand{\GWzenib}{\bibitem{GWzenib}
X. Gong and S.M. Webster,
{\it Regularity in the local CR embedding problem}, (submitted).}

\newcommand{\Keseon}{\bibitem{Keseon}
N. Kerzman, {\it
H\"older and $L\sp{p}$ estimates for solutions of $\bar \partial u=f$
 in strongly pseudoconvex domains},
 Comm. Pure Appl. Math. {\bf 24}(1971) 301--379.}

\end{document}